\documentclass[a4paper,10pt]{article}

\usepackage{amssymb,amsfonts,amsmath,amsthm}
\usepackage{graphicx,hyperref}
\usepackage{cite,caption}
\usepackage{mathtools}
\usepackage{booktabs}
\usepackage{newtxtext}
\usepackage[varvw]{newtxmath}
\usepackage{tikz}
\usepackage{tikz-3dplot}
\usepackage{pgfplots}
\usepackage{siunitx}
\pgfplotsset{compat = newest}
\usepackage{subcaption}
\usepackage{overpic}
\newcommand{\fz}{\frac}
\newcommand{\prz}[2]{ \frac{\partial{#1}}{\partial{#2}} }

\renewcommand{\Omega}{\varOmega}
\renewcommand{\Gamma}{\varGamma}
\renewcommand{\Psi}{\varPsi}
\renewcommand{\Pi}{\varPi}
\newcommand{\bR}{\mathbb{R}}

\newcommand{\rmD}{{\rm D}}

\renewcommand{\(}{\left(}
\renewcommand{\)}{\right)}
\newcommand{\rmT}{{\rm T}}
\renewcommand{\(}{\left(}
\renewcommand{\)}{\right)}
\newcommand{\defeq}{\vcentcolon=}

\providecommand{\keywords}[1]{\textit{Keywords:} #1}

\newtheorem{Rmk}{Remark}[section]
\newtheorem{example}{Example}
\numberwithin{equation}{section}
\allowdisplaybreaks[4]
\usepackage{pgfplots}
\usepackage{mathrsfs}
\usetikzlibrary{spy}
\usepackage{placeins}
\usepackage{tabularx,ragged2e,booktabs,caption}
\usepackage{multicol}
\usepackage{float}
\tikzset{
	block/.style={
		draw, 
		rectangle, 
		minimum height=1.5cm, 
		minimum width=3cm, align=center
	}, 
	line/.style={->,>=latex'}
}
\usepackage{environ} 
\usepackage{nccmath}
\usepackage{blindtext}
\usepackage{mathtools}
\NewEnviron{NORMAL}{% 
	\scalebox{1.2}{$\BODY$} 
} 
\makeatletter
\newcommand{\mathleft}{\@fleqntrue\@mathmargin0pt}
\newcommand{\mathcenter}{\@fleqnfalse}
\makeatother

\newcommand{\cB}{\color{blue}}

\usepackage{xcolor}
\hypersetup{
    colorlinks=true,
    citecolor=blue,
    linkcolor=red,
    urlcolor=orange,
}
\title{A two-step Lagrange--Galerkin scheme for the shallow water equations with a transmission boundary condition and its application to the Bay of Bengal region. Part~I: Flat bottom topography}
\author{
Md~Mamunur~Rasid$^{1,2}$,
Masato~Kimura$^{3}$,
Md~Masum~Murshed$^{4}$, \\
Erny~Rahayu~Wijayanti$^{5}$,
and
Hirofumi~Notsu$^{3,}$\footnote{Corresponding author}
}
\date{
\small
$^1$Division of Mathematical and Physical Sciences, Kanazawa University, Kakuma, Kanazawa 920-1192, Japan \smallskip\\
$^2$University of Rajshahi, Rajshahi-6205, Bangladesh \smallskip\\
$^3$Faculty of Mathematics and Physics, Kanazawa University, Kakuma, Kanazawa 920-1192, Japan \smallskip\\
$^4$Department of Mathematics, University of Rajshahi, Rajshahi-6205, Bangladesh \smallskip\\
$^5$Department of Mathematics, Gadjah Mada University, Yogyakarta 55281, Indonesia
\medskip\\
{\tt mamun.math@stu.kanazawa-u.ac.jp},\ \ 
{\tt mkimura@se.kanazawa-u.ac.jp},\ \ 
{\tt mmmurshed82@gmail.com},\ \ 
{\tt wijayanti.erny@gmail.com},\ \ 
{\tt notsu@se.kanazawa-u.ac.jp}
}
\begin{document}
\maketitle
\begin{abstract}
This study presents a two-step Lagrange--Galerkin scheme for the shallow water equations with a transmission boundary condition~(TBC). Firstly, the experimental order of convergence of the scheme is shown to see the second-order accuracy in time. Secondly, the effect of the TBC on a simple domain is discussed; the artificial reflections are kept from the Dirichlet boundaries and removed significantly from the transmission boundaries. Thirdly, the scheme is applied to a complex practical domain, i.e., the Bay of Bengal region, which is non-convex and includes islands. The effect of the TBC is discussed again for the complex domain; the artificial reflections are removed significantly from transmission boundaries, which are set on open sea boundaries. Based on the numerical results, it is revealed that the scheme has the following properties; (i)~the same advantages of Lagrange--Galerkin methods~(the CFL-free robustness for convection-dominated problems and the symmetry of the matrices for the system of linear equations); (ii)~second-order accuracy in time; (iii)~mass preservation of the function for the water level from the reference height (until the contact with the transmission boundaries of the wave); and (iv)~no significant artificial reflection from the transmission boundaries. The numerical results by the scheme are presented in this paper for the flat bottom topography of the domain. In the next part of this work, Part~II, the scheme will be applied to rapidly varying bottom surfaces and a real bottom topography of the Bay of Bengal region.
\par
\vspace{0.5em}
\noindent
\end{abstract}

\section{Introduction} 
\begin{figure}[!htbp]
	\centering
	\includegraphics[width=.5\textwidth]{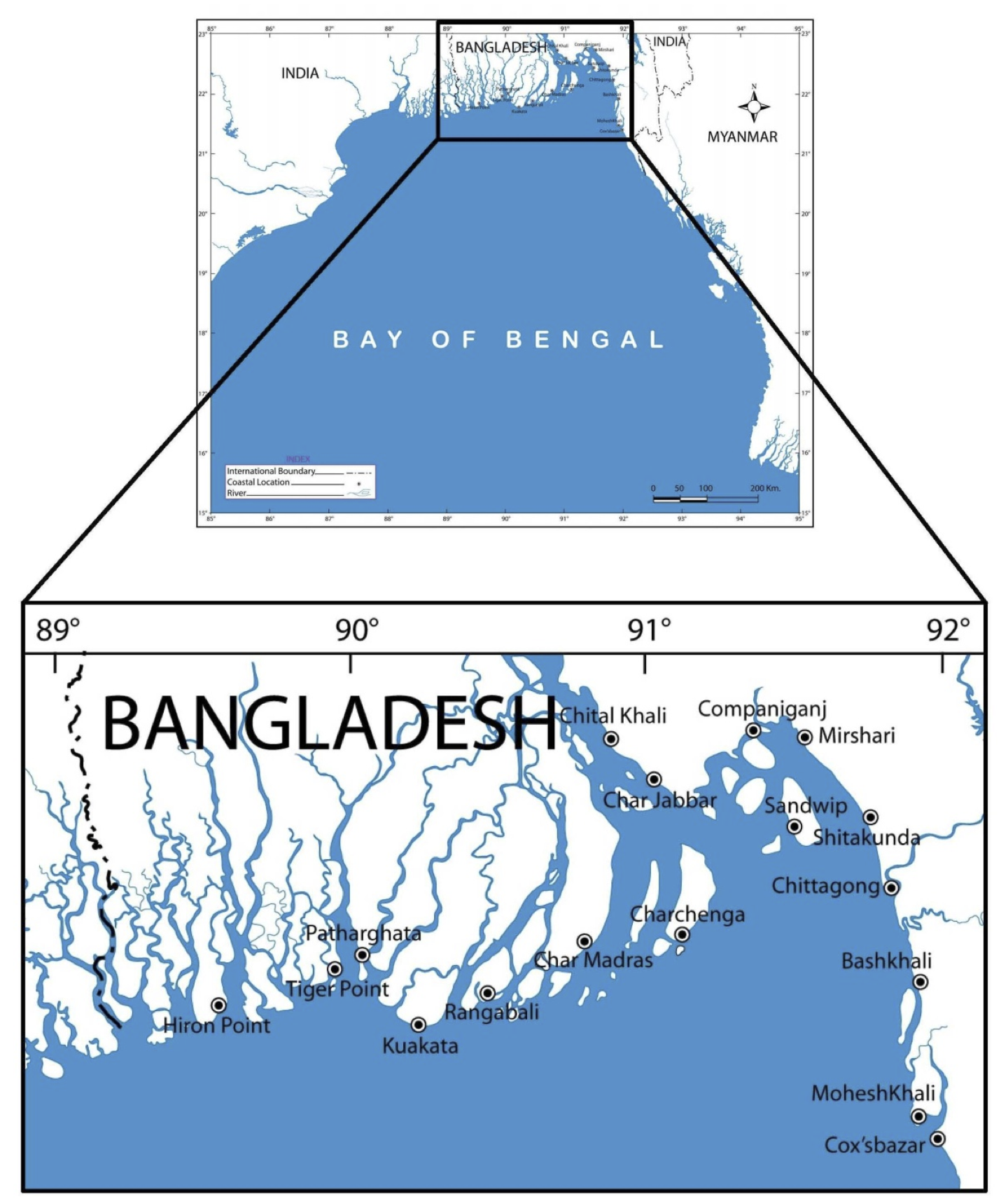}
	\caption{The Bay of Bengal region}\label{Bay}
\end{figure}

The system of the shallow water equations~(SWEs) is one of the most common models for describing fluid flow in rivers, channels, estuaries, and coastal areas and is often used for simulating tsunamis and storm surges in oceanic phenomena.
Natural disasters like tsunamis, cyclones, and storm surges cause a tremendous loss of lives and properties in the coastal areas in several regions.
According to~\cite{debsarma2009}, statistics show that about $5$\,\% of the global tropical cyclones form over the Bay of Bengal, and, on average, five to six storms form in this region every year, but with $80$\,\% of the global casualties.
The significant factors behind the heavy casualties are the shallow coastal water, thickly populated low-lying islands, highly curved coastal and island boundaries, river discharge, high astronomical tidal range, and favorable cyclone track, cf.~\cite{das1972} and Figure~\ref{Bay}.
That is why an effective storm surge prediction model and method are highly desired for the coastal region of Bangladesh to minimize the resulting damage from storm surges.
\par
Studies focusing on the Bay of Bengal region are found in~\cite{das1972, johns1981, roy1999polar, debsarma2009, paul2012tide, paul2013contribution, paul2018storm} and references therein.
For open sea boundaries, almost all the researchers implemented SWEs with a radiation-type boundary condition, which is comparable to a transmission boundary condition~(TBC) employed in~\cite{kanayamadan2016,murshed2021}.
Although for real problems, the finite element method is more suitable than the finite difference method because of the advantages of handling complex physical domains, geometries, or boundary conditions, as far as we know, there is no study to solve SWEs employing a  TBC for the Bay of Bengal region using the finite element method except~\cite{murshed_thesis2019}.
\par
The system of the SWEs consists of two equations, a pure convection equation for the total wave height and a modified Navier--Stokes momentum equation for the velocity derived by taking the average of function values in $x_{3}$-direction, cf.~\cite{kanayamadan2006,murshed2021}, which include the material derivatives in conservative and non-conservative forms, respectively.
For a time step size~$\Delta t > 0$, let $t^n\defeq n\Delta t$.
The so-called Lagrange--Galerkin method is the finite element method combined with the idea of the method of characteristics; the non-conservative and conservative material derivatives are discretized as, for a scalar-valued function~$\phi$ and a velocity~$u$, cf., e.g., \cite{EwiRus-1981, douglasrussell1982, pironneau1982, RuiTab-2010},
\begin{align*}
\Bigl[\prz{\phi}{t} + u\cdot\nabla\phi\Bigr](x,t^n) & = \fz{\phi^n(x)-\phi^{n-1}(x-u^{n}(x)\Delta t)}{\Delta t} + O(\Delta t), \\
\Bigl[\prz{\phi}{t} + \nabla \cdot (u\phi)\Bigr](x,t^n) & = \fz{\phi^n(x)-\phi^{n-1}(x-u^{n}(x)\Delta t)\gamma^n(x)}{\Delta t} + O(\Delta t),
\end{align*}
respectively, which are first-order approximations in time, where $x-u^{n}(x)\Delta t$ is an \textit{upwind} point of~$x$ with respect to~$u^{n}(x)$ and $\gamma^n$ is the Jacobian determinant of the mapping~$x-u^{n}(x)\Delta t$.
In general, the Lagrange--Galerkin method has two advantages; (i)~the CFL-free robustness for convection-dominated problems and (ii)~the symmetry of the resulting coefficient matrices for the system of linear equations.
In addition to the four pioneering works above, many authors have proposed the ideas of this type of approximations in the context of the finite element method, cf.~\cite{EwiRusWhe-1983, suli_1988, Pir-1989, BMMR-1997, AchGue-2000, ruitabata2002, BerNogVaz-2006_part1, BerNogVaz-2006_part2, ChrWal-2008, N-2008-JSCES, PirTab-2010, BenBer-2011, BenBer-2012_part1, BenBer-2012_part2, BerSaa-2012, BerGalSaa-2012, notsuruitabata2013, notsutabata2015, notsu2016error, notsutabata2norder2016, TabUch-2016-CD, LNS-2015, LMNT-2017_Peterlin_Oseen_Part_I, LMNT-2017_Peterlin_Oseen_Part_II, TabUch-2018-NS, Uch-2019, ColCarBer-2020, ColCarBer-2021} and references therein.
When we focus on the SWEs, to the best of our knowledge, Murshed et al.~\cite{murshed2021} and Murshed~\cite{murshed_thesis2019} firstly solved the SWEs with a TBC by a (single-step) Lagrange--Galerkin scheme of first-order in time for a flat bottom topography.
Recently, a two-step mass-preserving Lagrange--Galerkin scheme of second order in time for conservative convection-diffusion problems has been proposed and analyzed with error estimates in~\cite{FutKolNotSuz-2022}.
\par
In this paper, we present a new two-step Lagrange--Galerkin scheme to solve the SWEs together with a TBC, which is of second order in time and maintains the two advantages of the Lagrange--Galerkin methods, i.e., the CFL-free robustness and the symmetry of the resulting matrices.
The two material derivatives are discretized based on the ideas of two-step methods proposed for the non-conservative form in~\cite{EwiRus-1981, EwiRusWhe-1983, BMMR-1997, notsutabata2norder2016} and the conservative form in~\cite{FutKolNotSuz-2022}.
Firstly, preparing an artificial exact solution, we observe our scheme's experimental order of convergence~(EOC) to see the second-order accuracy in time on a simple (square) domain.
Since long (real-)time computations on a mesh refined locally are needed in practical problems, the CFL-free second-order accuracy in time of our scheme is a significant advantage, enabling us to employ a more extensive time increment compared with first-order numerical methods.
Secondly, we observe the effect of the TBC on a simple (square) domain, and the artificial reflections are kept from the Dirichlet boundaries and removed significantly from the transmission boundaries.
Thirdly, our scheme is applied to the Bay of Bengal region, which is non-convex, includes islands, and is, therefore, a complex domain.
We again observe the effect of the TBC for this realistic domain.
The artificial reflections are removed significantly from the transmission boundaries, which are set on open sea boundaries.
We also study the effect of a position of an open sea boundary with the TBC and reveal that it is sufficiently small to neglect.
In~\cite{murshed2021}, energy estimates for the SWEs were given, where the $L^2$-norm of the water level from the reference height was an important value related to the potential energy.
Focusing on the energy and the mass of the water level function, we observe the $L^2$-norm and the mass of the water level function, which show the effectiveness of the TBC.
\par
From the computations, we show that our new scheme has the following properties; (i)~the same advantages of Lagrange--Galerkin methods; (ii)~second-order accuracy in time; (iii)~mass preservation of the function of the water level from the reference height (until the contact with the transmission boundaries of the wave); and (iv) no significant artificial reflection from the transmission boundaries.
All of the numerical results in this paper, Part~I, are for the flat bottom topography, and the non-homogeneous bottom topography will be studied in our forthcoming paper, Part~II.
\par
The outline of this paper is as follows. 
Section~\ref{sec:scheme} presents a two-step Lagrange--Galerkin scheme for the SWEs together with a TBC, which is of second order in time.
In Section~\ref{sec:numer_square_domain}, numerical results for simple square domains are shown to observe the second-order accuracy in time and the effect of TBC.
In Section~\ref{sec:application}, our scheme is applied to the Bay of Bengal region, where the domain is non-convex and complex.
In Section~\ref{sec:conclusions}, conclusions are given.
The data for choosing the constant  $c_0$ required in the TBC is given in the Appendix.
\section{A two-step Lagrange--Galerkin scheme}
\label{sec:scheme}
We introduce some notations to be used in this paper.
$\Omega$ is a bounded spatial domain in~$\bR^{2}$, $\Gamma \defeq \partial \Omega$ is the boundary of~$\Omega$, and $(0, T)$ is a temporal domain in~$\bR_+ \, (\defeq \{x \in \bR; x > 0 \})$ for a positive constant~$T$.
We use the Lebesgue space~$L^2(\Omega)$ and the Sobolev space~$H^1(\Omega)$.
Let $(\cdot,\cdot)$ be the inner product in $L^2(\Omega)$, i.e., $(f,g) \defeq \int_{\Omega}f(x)g(x)dx$ for $f, g\in L^2(\Omega)$.
We employ the same notation~$(\cdot, \cdot)$ to represent the~$L^2(\Omega)$ inner product for scalar-, vector-, and matrix-valued functions.
Let $A : B$ be the tensor product defined by~$A : B \defeq \sum_{i,j=1}^{2}A_{ij}B_{ij} = {\rm tr}(AB^\top)$ for $A, B\in\bR^{2\times 2}$.
\subsection{Statement of the problem}
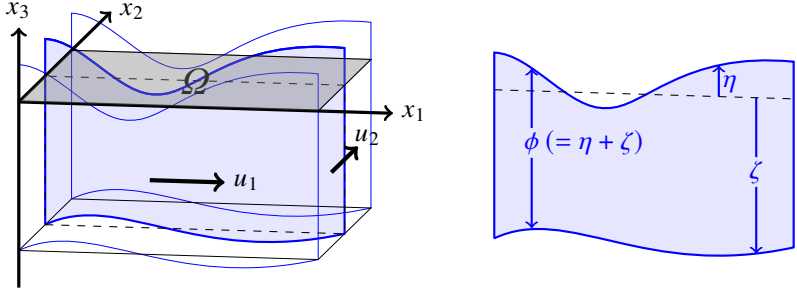
\begin{figure}[!htbp]
\centering
\tdplotsetmaincoords{80}{10}
\begin{tikzpicture}[tdplot_main_coords,scale=1]
\coordinate (O) at (0,0,0);
\coordinate (O1) at (0,0,-2.5);
\coordinate (O2) at (0,0,1);
\coordinate[label=right:$x_1$] (X) at (5,0,0);
\coordinate[label=right:$x_2$] (Y) at (0,7,0);
\coordinate[label=above:$x_3$] (Z) at (0,0,1);
\coordinate (A) at (0,0,-2);
\coordinate (B) at (4,0,-2);
\coordinate (C) at (4,4,-2);
\coordinate (D) at (0,4,-2);
\coordinate (E) at (0,0,0);
\coordinate (F) at (4,0,0);
\coordinate (G) at (4,4,0);
\coordinate (H) at (0,4,0);
\coordinate (P0) at (0,2,0.5);
\coordinate (P1) at (0.5,2,0.5);
\coordinate (P2) at (1,2,-0.25);
\coordinate (P3) at (1.5,2,-0.2);
\coordinate (P4) at (2,2,-0.2);
\coordinate (P5) at (2.5,2,0.6);
\coordinate (P6) at (4,2,0.5);
\draw[blue,thick] (P0) .. controls (P1) and (P2) .. (P3) .. controls (P4) and (P5) .. (P6);
\draw[blue,thick] (0,2,-2) .. controls (1,2,-1.5) and (2,2,-2.5) .. (4,2,-2) -- (4,2,0.5) .. controls (2.5,2,0.6) and (2,2,-0.2) .. (1.5,2,-0.2) .. controls (1,2,-0.25) and (0.5,2,0.5) .. (0,2,0.5) -- (0,2,-2);
\draw[blue,fill=blue,opacity=0.1] (0,2,-2) .. controls (1,2,-1.5) and (2,2,-2.5) .. (4,2,-2) -- (4,2,0.5) .. controls (2.5,2,0.6) and (2,2,-0.2) .. (1.5,2,-0.2) .. controls (1,2,-0.25) and (0.5,2,0.5) .. (0,2,0.5) -- (0,2,-2);
\draw[blue,very thin] (0,0,-2) .. controls (1,0,-1.5) and (2,0,-2.5) .. (4,0,-2) -- (4,0,0.5) .. controls (2.5,0,0.6) and (2,0,-0.2) .. (1.5,0,-0.2) .. controls (1,0,-0.25) and (0.5,0,0.5) .. (0,0,0.5) -- (0,0,-2);
\draw[blue,very thin] (0,4,-2) .. controls (1,4,-1.5) and (2,4,-2.5) .. (4,4,-2) -- (4,4,0.5) .. controls (2.5,4,0.6) and (2,4,-0.2) .. (1.5,4,-0.2) .. controls (1,4,-0.25) and (0.5,4,0.5) .. (0,4,0.5) -- (0,4,-2);
\coordinate[label=center:{\Large$\Omega$}] (Center) at (2,2,0);
\draw[blue,thick] (6,2,-2) .. controls (7,2,-1.5) and (8,2,-2.5) .. (10,2,-2) -- (10,2,0.5) .. controls (8.5,2,0.6) and (8,2,-0.2) .. (7.5,2,-0.2) .. controls (7,2,-0.25) and (6.5,2,0.5) .. (6,2,0.5) -- (6,2,-2);
\draw[blue,fill=blue,opacity=0.1] (6,2,-2) .. controls (7,2,-1.5) and (8,2,-2.5) .. (10,2,-2) -- (10,2,0.5) .. controls (8.5,2,0.6) and (8,2,-0.2) .. (7.5,2,-0.2) .. controls (7,2,-0.25) and (6.5,2,0.5) .. (6,2,0.5);
\draw[blue,<-,thick] (9.5,2,-2.1) -- (9.5,2,-1.2);
\draw[blue,-,thick] (9.5,2,-0.8) -- (9.5,2,0);
\draw[blue,<-,thick] (6.5,2,-1.85) -- (6.5,2,-0.9);
\draw[blue,->,thick] (6.5,2,-0.5) -- (6.5,2,0.3);
\draw[blue,->,thick] (9,2,0) -- (9,2,0.42);
\draw[dashed,-] (6,2,0) -- (10,2,0);
\coordinate[label=center:$\cB \zeta$] (Qmid) at (9.5,2,-1);
\coordinate[label=right:${\cB \phi \, (= \eta + \zeta)}$] (Rmid) at (6.27,2,-0.7);
\coordinate[label=center:$\cB \eta$] (Smid) at (9.15,2,0.18);
\draw[dashed] (0,2,-2) -- (4,2,-2) -- (4,2,0) -- (0,2,0) -- cycle;
\draw[->,very thick] (O) -- (X);
\draw[->,very thick] (O) -- (Y);
\draw[->,very thick] (O1) -- (O2);
\draw[very thin] (A) -- (B) -- (C) -- (D) -- cycle;
\draw[-] (E) -- (F) -- (G) -- (H) -- cycle;
\draw[fill=gray,opacity=0.4] (E) -- (F) -- (G) -- (H) -- cycle;
\draw[->,ultra thick] (1.75,0,-1) -- (2.75,0,-1);
\draw[->,ultra thick] (4,1.0,-1) -- (4,3,-1);
\coordinate[label=right:$u_1$] (u1) at (2.75,0,-1);
\coordinate[label=center:$u_2$] (u2) at (4,3.75,-1);
\end{tikzpicture}
\caption{Diagrams for the problem; left: the domain~$\Omega$ and the velocity~$u=(u_1, u_2)^\top$; right:~the total wave height~$\phi = \eta+\zeta$.}
\label{fig:diagrams}
\end{figure}

Our problem is to find $(\phi,u)\colon \Omega\times (0,T) \to \bR\times\bR^2$ such that
\begin{subequations}\label{eqn1}
	\begin{align}
		\frac{\partial \phi}{\partial t}+\nabla\cdot(u\phi) & = f  && \mbox{in}~\Omega\times(0,T), \label{eqn1_1} \\
		\rho\phi\Bigl[ \frac{\partial u}{\partial t}+(u\cdot\nabla)u\Bigr]-2\mu\nabla\cdot(\phi D(u)) & + \rho g\phi\nabla\eta = F && \mbox{in}~\Omega\times(0,T), \label{eqn1_2} \\
		\phi & = \eta+\zeta && \mbox{in}~\Omega\times(0,T), \label{eqn1_3} \\
	u&=0&& \mbox{on}~\Gamma_\rmD\times(0,T), \\
	u&=c_{0}\sqrt{g\zeta} \, \frac{\eta}{\phi}\, n&& \mbox{on}~\Gamma_\rmT\times(0,T), \\
	(\phi, u) & = (\phi^{0}, u^{0}) && \mbox{in}~\Omega,~\mbox{at}~t=0, 
\end{align}
\end{subequations}
where the total wave height and the velocity are denoted by~$\phi$ and~$u=(u_{1},u_{2})^\top$, respectively, 
the water level from the reference height and the depth of water level from the reference height, i.e., bottom topography, are represented by
$\eta\colon\Omega\times (0,T) \to\bR$ and $\zeta\colon \Omega\to\bR_{+}$, respectively,
a pair of external forces is given by~$(f, F)\colon\Omega\times (0,T)\to\bR\times\bR^{2}$,
a pair of initial values is given by~$(\phi^{0}, u^{0})\colon\Omega\to\bR\times\bR^{2}$,
density and viscosity constants of water are denoted by~$\rho>0$ and $\mu>0$,
the gravity constant is given by~$g>0$,
the strain-rate tensor~$D(u)$ is defined by
\[
D(u) \defeq \frac{1}{2} \left[ \nabla u+(\nabla u)^\top \right],
\]
and the outward unit normal vector is denoted by~$n\colon \Gamma\to\bR^2$, cf. Figure~\ref{fig:diagrams}.
We suppose that the boundary~$\Gamma$ is divided into two non-overlapping parts, $\Gamma_\rmD$ and~$\Gamma_\rmT$, i.e., $\overline{\Gamma} = \overline{\Gamma}_\rmD \cup \overline{\Gamma}_\rmT$ and $\Gamma_\rmD \cap \Gamma_\rmT = \emptyset$, where the subscripts ``$\rmD$'' and ``$\rmT$'' imply Dirichlet and transmission boundaries, respectively.
A positive constant~$c_{0}$ is chosen suitably to remove the artificial reflection, and, throughout this paper, we employ~$c_0 = 0.9$, which is determined based on numerical experiments given in Appendix.
We consider homogeneous flat bottom topography in this paper, Part~I, and non-homogeneous bottom topography in our forthcoming paper, Part~II.
\subsection{Presentation of the scheme}\label{subsection:scheme}
%\noindent Weak formulation of the problem ~\eqref{eqn1}:\\
\noindent Let $\Psi\coloneqq L^{2}(\Omega)$, $Y\coloneqq H^{1}(\Omega)^{2}$,
\[
V(G)\coloneqq \bigl\{ v\in Y;\ v=0~\mbox{on}~\Gamma_\rmD~ \mbox{and} ~v=G~\mbox{on}~\Gamma_\rmT \bigr\}
\]
for a function~$G\colon \Gamma_\rmT\to \bR^{2}$, and $V\defeq V(0)$.
%\par
We introduce a $\phi$-dependent function, $G(\phi) = G(\phi; \eta) \colon \Gamma_\rmT\to \bR^2$, defined by
\[
G(\phi) = G(\phi; \eta)\defeq c_{0}\sqrt{g\zeta} \, \frac{\eta}{\phi} n.
\]
Assume $\phi^0 \in \Psi$, $\eta^0 \defeq \phi^0-\zeta \in \Psi$ and $u^0\in V(G(\phi^0;\eta^0))$.
A weak formulation to problem~\eqref{eqn1} is to find $\{(\phi, u)(t)\in \Psi\times V(G(\phi(t); \eta(t)));~t\in(0,T)\}$ such that, for $t\in (0,T)$,
\begin{subequations}\label{weak}
	\begin{align}
	\Bigl( \frac{\partial \phi}{\partial t}+\nabla\cdot(u\phi),\psi \Bigr) & = (f,\psi)\qquad \forall \psi \in \Psi, \\
	\rho\biggl( \phi\Bigl[ \frac{\partial u}{\partial t}+(u\cdot\nabla)u \Bigr], v\biggr) + a(u,v;\phi) + b(\eta,v;\phi) & = (F,v) \qquad \forall v \in V, \\
	\phi & = \eta+\zeta,
\end{align}
\end{subequations}
with the initial condition~$(\phi, u) (0) = (\phi^{0}, u^{0}) \in \Psi \times V(G(\phi^0; \eta^0))$,
%
%\begin{align}
%	(\phi, u) (0) = (\phi^{0}, u^{0}) \in \Psi \times V(G(\phi^0; \eta^0)),
%\end{align}
%
where the bilinear forms~$a(\cdot,\cdot\,;\phi)\colon Y\times Y\to \bR$ and~$b(\cdot,\cdot\,;\phi)\colon \Psi\times Y\to \bR$ are defined by
\begin{align*}
a(u,v;\phi) & \defeq 2\mu \bigl( \phi D(u), D(v) \bigr), &
b(\eta,v;\phi) & \defeq \rho g \bigl( \phi \nabla\eta, v\bigr).
\end{align*}
Now, we present our scheme for solving problem~\eqref{eqn1}.
Let $\mathcal{T}_{h}=\{K\}$ be a partition of $\overline{\Omega}$ by triangular  elements, $h$ be the maximum diameter of~$K \in\mathcal{T}_h$, and $\Omega_h\defeq {\rm int}\( \bigcup_{K\in\mathcal{T}_h} K \)$ be an approximated domain.
Although it holds that~$\Omega \neq \Omega_h$ in general, we assume~$\Omega = \Omega_h$ throughout the paper to avoid the complexity of introducing many symbols.
We define finite element spaces, $\Psi_h$, $Y_{h}$ and $V_{h}(G)$, corresponding to~$\Psi$, $Y$ and~$V(G)$ by
\begin{align*}
\Psi_h & \defeq 
\{\psi_{h}\in C^{0}(\overline{\Omega}); \psi_{h \vert K}\in P_{1}(K)\ \forall K\in \mathcal{T}_{h}\}, \\
Y_{h} & \defeq \{v_{h}\in C^{0}(\overline{\Omega})^2;\ v_{h \vert K} \in P_{1}(K)^2\ \forall K\in \mathcal{T}_{h}\}, \\
V_h (G) & \defeq \{v_h\in Y_h;\ v_h = 0\ \mbox{on}\ \Gamma_\rmD\ \mbox{and}\ v_h = G\ \mbox{on}\ \Gamma_\rmT\},
\end{align*}
and set $V_h \defeq V_h (0)$, where the function~$G \colon \Gamma_\rmT\to\bR^2$ is assumed to be a piecewise linear function.
\par
Let $\Delta t$ be a time increment, $N_T\defeq \lfloor T/\Delta t\rfloor$ a total number of time steps, and $t^n\defeq n\Delta t$ a time at $n$-th time step.
For $v\colon \Omega\to \bR^2$, we define mappings $X_{1}[v], \tilde{X}_{1}[v]\colon \Omega\to \bR^2$ and $\gamma_{1}[v], \tilde{\gamma}_{1}[v]\colon \Omega\to \bR$ by
\begin{align*}
	X_{1}[v](x) & \defeq x-\Delta t\, v(x), & 
	\tilde{X}_{1}[v](x) & \defeq x- 2 \Delta t \, v(x),\\
	\gamma_{1}[v](x) & \defeq \det \biggr( \prz{X_{1}[v]}{x}(x) \biggr), &
	\tilde{\gamma}_{1}[v](x) & \defeq \det \biggl( \prz{\tilde{X}_{1}[v]}{x}(x) \biggr).
\end{align*}
For $\bigl\{\phi^{n}\bigr\}_{n=0}^{N_T}$ and $\bigl\{u^{n}\bigr\}_{n=0}^{N_T}$, we define an operator~$\mathcal{A}_{\Delta t}[u]\phi^n$ by, for $n=1,\ldots, N_T$,
\begin{align*}
\mathcal{A}_{\Delta t}[u]\phi^{n} & \coloneqq
\left\{
  \begin{aligned}
    &\mathcal{A}_{\Delta t}^{(1)}[u]\phi^n && (n=1),\\
	&\mathcal{A}_{\Delta t}^{(2)}[u]\phi^{n}&& (n\geq 2),
  \end{aligned}
\right.
\end{align*}
where
\begin{align*}
\mathcal{A}_{\Delta t}^{(1)}[u]\phi^{n} & \coloneqq \frac{\phi^{n}-\phi^{n-1}\circ X_{1}[u^{n-1}]\gamma_{1}[u^{n-1}]}{\Delta t}, \\
\mathcal{A}_{\Delta t}^{(2)}[u]\phi^{n} & \coloneqq \frac{3\phi^{n}-4\phi^{n-1}\circ X_{1}[u^{n\ast}]\gamma_{1}[u^{n\ast}]+\phi^{n-2}\circ\tilde{X}_{1}[u^{n\ast}]\tilde{\gamma}_{1}[u^{n\ast}]}{2\Delta t}.
\end{align*}
The composition of functions is represented by the symbol $\circ$, i.e.,
\[
\bigl( \psi\circ X_1[v] \bigr) (x) = \psi \bigl( X_1[v](x) \bigr),
\]
and the function~$u^{n\ast}\colon \Omega\to \bR^{2}$ is defined by
\[
u^{n\ast} \defeq 2u^{n-1}-u^{n-2},
\]
which is a second-order temporal approximation of $u^n$ if $u$ is sufficiently smooth.
We also define, for $\bigl\{w^{n}\bigr\}_{n=0}^{N_T}$,
\[
\mathcal{B}_{\Delta t}[w]u^{n}\coloneqq\left\{\begin{array}{ccc}
	\mathcal{B}_{\Delta t}^{(1)}[w]u^{n}&& (n=1),\\
	\mathcal{B}_{\Delta t}^{(2)}[w]u^{n}&& (n\geq 2),
\end{array}\right.
\]
where
\begin{align*}
\mathcal{B}_{\Delta t}^{(1)}[w]u^{n} & \coloneqq \frac{u^{n}-u^{n-1}\circ X_{1}[w^{n-1}]}{\Delta t}, \\
\mathcal{B}_{\Delta t}^{(2)}[w]u^{n} & \coloneqq \frac{3u^{n}-4u^{n-1}\circ X_{1}[w^{n\ast}]+u^{n-2}\circ\tilde{X}_{1}[w^{n\ast}]}{2\Delta t}.
\end{align*}
The two-step Lagrange--Galerkin scheme is to find $\{ (\phi_{h}^{n}, u_{h}^{n})\in \Psi_h \times V_h ( G(\phi_h^n;\eta_h^n) );$ $n=1, \ldots, N_T \}$ such that, for $n=1,2, \ldots, N_{T}$,
\begin{subequations}\label{scheme}
	\begin{align}
		\bigl(\mathcal{A}_{\Delta t}[u_{h}]\phi_{h}^{n},\psi_{h}\bigr) & = \bigl(f^n,\psi_{h}\bigr) & \forall \psi_{h}\in \Psi_{h}, 
		\label{eqn1LG} \\
		\rho\bigl(\phi_{h}^{n}\mathcal{B}_{\Delta t}[u_{h}]u_{h}^{n},v_{h}\bigr) +a\bigl(u_{h}^{n},v_{h};\phi_{h}^{n}\bigr) +b\bigl(\eta_{h}^{n}, & v_{h}; \phi_{h}^{n}\bigr) \notag\\
		& = \bigl(F^{n},v_{h}\bigr) &\forall v_{h}\in V_{h}, 
		\label{eqn2LG}\\
		\phi_{h}^{n}&=\eta_{h}^{n}+\Pi_{h}\zeta,
		\label{eqn3LG}
	\end{align}
with an initial condition
\begin{align}
	\bigl(\phi_h^0, u_h^0\bigr) = \bigl(\Pi_h\phi^{0}, \Pi_h u^{0}\bigr) \in \Psi_h \times Y_h,
\end{align}
\end{subequations}	
where the Lagrange interpolation operator is denoted by~$\Pi_{h}\colon C(\overline{\Omega})\to \Psi_{h}$, which is also used for the vector-valued function~$u^0$, i.e., $\Pi_h u^0 \in Y_h$.
\begin{Rmk}
(i)~At each time step, we obtain~$\phi_{h}^{n}\in\Psi_{h}$ from~\eqref{eqn1LG} and $u_{h}^{n}\in V_{h}(G(\phi_h^n;\eta_h^n))$ from~\eqref{eqn2LG} combined with~\eqref{eqn3LG}.
\smallskip\\
(ii)~We need $\mathcal{A}_{\Delta t}^{(1)}[u]$ and $\mathcal{B}_{\Delta t}^{(1)}[w]$ due to the lack of the functions~$\phi_h^{n-2}$ and~$u_h^{n-2}$ for~$n=1$, which are used for $\mathcal{A}_{\Delta t}^{(2)}[u_h]\phi_h^n$ and $\mathcal{B}_{\Delta t}^{(2)}[u_h]u_h^n$ for~$n \ge 2$.
\smallskip\\
(iii)~The two-step methods in conservative and non-conservative forms, $\mathcal{A}_{\Delta t}^{(2)}[u_h]\phi_h^n$ and $\mathcal{B}_{\Delta t}^{(2)}[u_h]u_h^n$, are developed and analyzed for convection-diffusion problems in~\cite{FutKolNotSuz-2022,EwiRus-1981}.
\smallskip\\
(iv)~It is discussed in~\cite{FutKolNotSuz-2022, notsutabata2norder2016} that the one-time use of first-order single-step methods, $\mathcal{A}_{\Delta t}^{(1)}[u_h]\phi_h^n$ and $\mathcal{B}_{\Delta t}^{(1)}[u_h]u_h^n$, has no loss of convergence orders in discrete versions of $L^\infty(0,T; L^2(\Omega))$- and $L^2(0,T; H^1(\Omega))$-norms for a convection-diffusion equation and the Navier--Stokes equations, respectively.
\end{Rmk}
\section{Numerical results in square domains}
\label{sec:numer_square_domain}
In this section, numerical results via FreeFem++~\cite{MR3043640} are presented to see the experimental order of convergence~(EOC) and the effect of the TBC in square domains.
We call scheme~\eqref{scheme} LG2, and also call scheme~\eqref{scheme} replacing~$\mathcal{A}_{\Delta t}$ and~$\mathcal{B}_{\Delta t}$ with~$\mathcal{A}_{\Delta t}^{(1)}$ and~$\mathcal{B}_{\Delta t}^{(1)}$, respectively, LG1 which is a (single-step) Lagrange--Galerkin scheme of first order in time.
\subsection{Experimental order of convergence}
\label{subsec:EOC}
We solve Examples~\ref{ex1} and~\ref{ex2} below by LG1 and LG2 and compare the experimental orders of convergence~(EOCs).
\begin{example}[$\Gamma = \Gamma_\rmD$]\label{ex1}
In problem~\eqref{eqn1}, we set $\Omega=(0,1)^{2}$, $\Gamma = \Gamma_\rmD$~$(\Gamma_\rmT = \emptyset)$, $T=1$, $g=\rho=\mu=\zeta=1$, and the function $\eta^{0}$, $u^{0}$, $f$ and $F$ are given so that the exact solution is
\begin{align*}
\phi(x,t) & = 1+\frac{\sin\pi x_{1}\sin\pi x_{2}(2+\sin\pi t)}{8}, &
u(x,t) & = \frac{\sin\pi x_{1}\sin\pi x_{2}(2+\sin\pi t)}{3}
\begin{bmatrix}
1\\
1
\end{bmatrix}.
\end{align*}
\vspace*{-1.5em}
\end{example}
\begin{example}[$\Gamma = \overline{\Gamma}_\rmD \cup \overline{\Gamma}_\rmT$]\label{ex2}
In Example~\ref{ex1}, we replace $\Gamma_\rmT$ and $\Gamma_\rmD$ with $\Gamma_\rmT = \{x\in\Gamma;\ x_2=0\}$ and $\Gamma_\rmD = \Gamma \setminus \overline{\Gamma}_\rmT$, respectively.
\end{example}
For a numerical solution~$z_h=\{z_h^n\}_{n=0}^{N_T}$ and its exact solution~$z=\{z^n\}_{n=0}^{N_T}$, we introduce notations of errors, $E_i(z)$, $i=0, 1$, defined by
\begin{align*}
E_0(z) & \defeq\fz{\|z_h-z\|_{\ell^\infty(L^2)}}{\|z\|_{\ell^\infty(L^2)}}, & 
E_1(z) & \defeq\fz{\|\nabla(z_h-z)\|_{\ell^\infty(L^2)}}{\|\nabla z\|_{\ell^\infty(L^2)}},
\end{align*}
where $\|\cdot\|_{\ell^\infty(L^2)}$ is a norm given by
\[
\|z\|_{\ell^\infty(L^2)} \defeq \max\{\|z^n\|_{L^2(\Omega)};\ n=0,\ldots,N_T\}.
\]
\par
\begin{figure}[!htbp]
\centering
\includegraphics[width=0.3\textwidth]{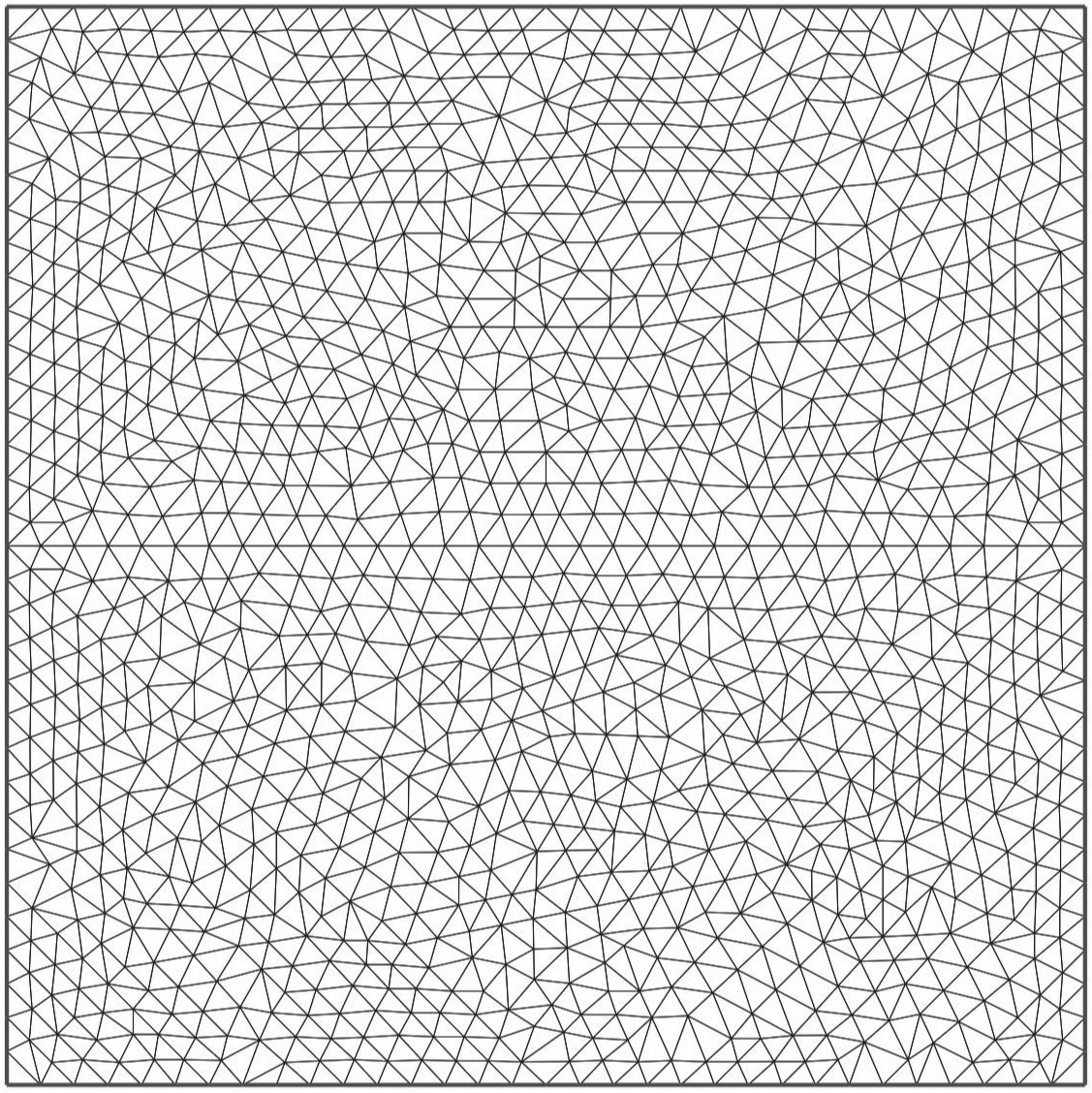}
\caption{A sample mesh with $N=32$ for Example~\ref{ex1}.}
\label{fig:sample_mesh}
\end{figure}
\begin{figure}[!htbp]
	\centering
		\begin{subfigure}[b]{0.21\textwidth}
			\centering
			\includegraphics[width=0.8\textwidth,clip]{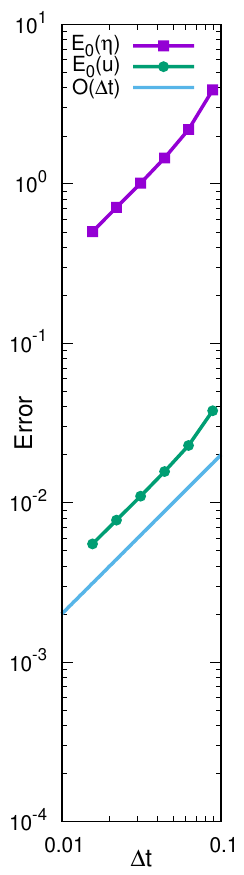}
			\caption{(i) LG1 for Ex.\ref{ex1}}
		\end{subfigure}
		\hspace{.1cm}   
		\begin{subfigure}[b]{0.21\textwidth}
			\centering
			\includegraphics[width=0.8\textwidth,clip]{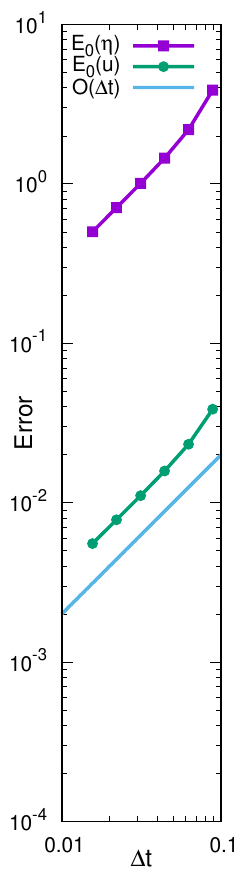}
			\caption{(ii) LG1 for Ex.\ref{ex2}}
		\end{subfigure}
		\qquad
		\begin{subfigure}[b]{0.21\textwidth}
			\centering
			\includegraphics[width=0.8\textwidth,clip]{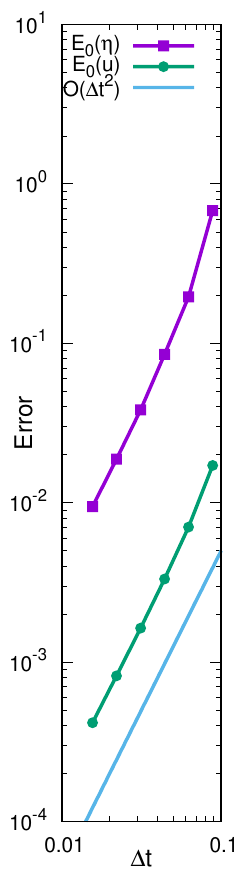}
			\caption{(iii) LG2 for Ex.\ref{ex1}}
		\end{subfigure} 
		\hspace{.1cm}
		\begin{subfigure}[b]{0.21\textwidth}
			\centering
			\includegraphics[width=0.8\textwidth,clip]{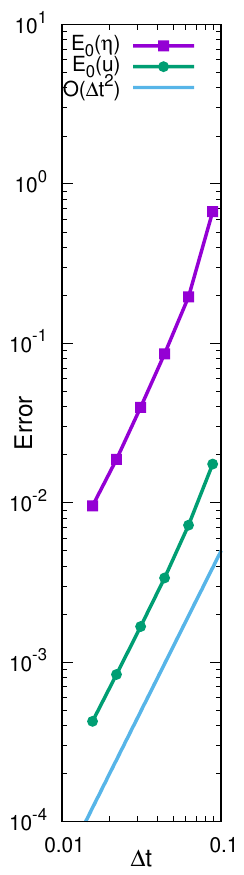}
			\caption{(iv) LG2 for Ex.\ref{ex2}}
		\end{subfigure}
	\caption{Graphs of errors~$E_0(\eta)$ and~$E_0(u)$ in logarithmic scale by LG1 for Example~\ref{ex1}~(i) and Example~\ref{ex2}~(ii), and by LG2 for Example~\ref{ex1}~(iii) and Example~\ref{ex2}~(iv).}
	\label{fig:error0}
\end{figure}
\begin{figure}[!htbp]
	\centering
	\begin{subfigure}[b]{0.21\textwidth}
		\centering
		\includegraphics[width=0.8\textwidth,clip]{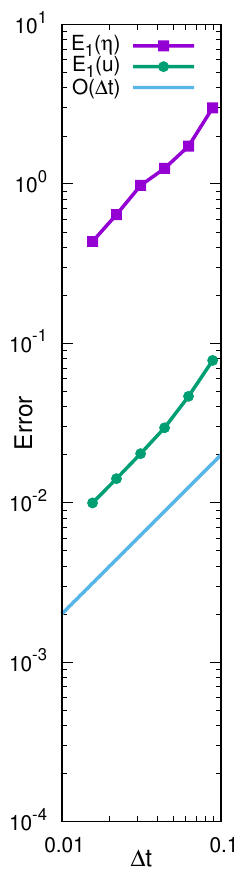}
		\caption{(i) LG1 for Ex.\ref{ex1}}
	\end{subfigure}
	\hspace{.1cm}   
	\begin{subfigure}[b]{0.21\textwidth}
		\centering
		\includegraphics[width=0.8\textwidth,clip]{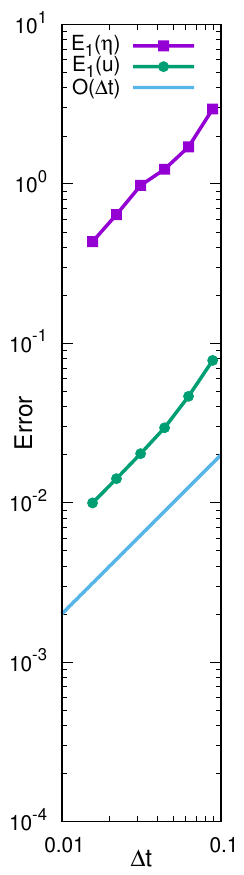}
		\caption{(ii) LG1 for Ex.\ref{ex2}}
	\end{subfigure}
	\qquad
	\begin{subfigure}[b]{0.21\textwidth}
		\centering
		\includegraphics[width=0.8\textwidth,clip]{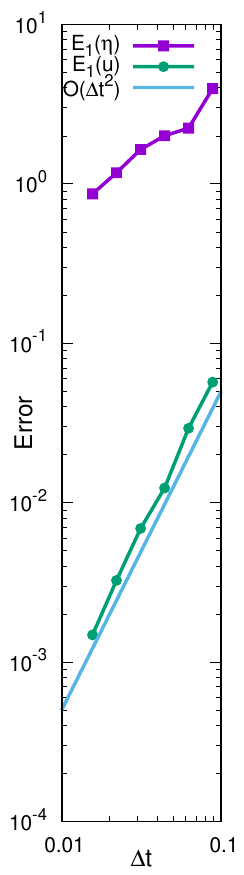}
		\caption{(iii) LG2 for Ex.\ref{ex1}}
	\end{subfigure} 
	\hspace{.1cm}   
	\begin{subfigure}[b]{0.21\textwidth}
		\centering
		\includegraphics[width=0.8\textwidth,clip]{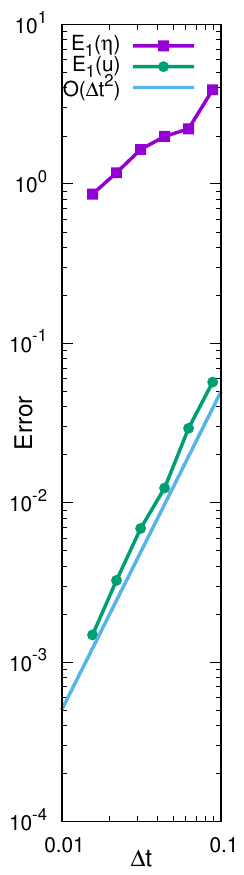}
		\caption{(iv) LG2 for Ex.\ref{ex2}}
	\end{subfigure}
	\caption{Graphs of errors~$E_1(\eta)$ and~$E_1(u)$ in logarithmic scale by LG1 for Example~\ref{ex1}~(i) and Example~\ref{ex2}~(ii), and by LG2 for Example~\ref{ex1}~(iii) and Example~\ref{ex2}~(iv).}
	\label{fig:error1}
\end{figure}

Let $N$ be a division number of each side of the unit square domain~$\Omega$ and $h\defeq 1/N$ a representative mesh size.
We prepare non-uniform triangulations of~$\Omega$, $\mathcal{T}_h$, for $N = 8, 16, 32, 64, 128$ and~$256$, cf.~Figure~\ref{fig:sample_mesh} for~$N=32$.
Choosing $\Delta t = 0.25\sqrt{h}$, we compute the errors, $E_i (\eta)$ and~$E_i (u)$, $i=0, 1$, by~LG1 and~LG2.
Figures~\ref{fig:error0} and~\ref{fig:error1} show graphs of the errors of $E_0(\cdot)$ and $E_1(\cdot)$, respectively, in logarithmic scale by LG1 for Example~\ref{ex1}~(i) and Example~\ref{ex2}~(ii), and by LG2 for Example~\ref{ex1}~(iii) and Example~\ref{ex2}~(iv), and the values of errors and their EOCs are given in Tables~\ref{table:ex1} and~\ref{table:ex2}.
We observe that LG2 is of second order in time numerically and that the order is higher than that of~LG1.
Although $E_1(\eta)$ is not of second order in time, it is natural as equation~\eqref{eqn1_1} for~$\phi \ (= \eta + \zeta)$ does not include any diffusion term.
%
%
%
%
%%%%%%%%%%%%%%%%%%%%%%%%%%%%%%%%%%%%%%%%%%%%%%%%%%%%%%%%%%%%%%%%%%%%%%%%%%%%%%%%%%%%%%%%%
\begin{table}[!htbp]
	\centering
	\caption{Values of $E_i(\eta)$ and $E_i(u)$, $i=0,1$, by schemes LG1 and LG2 for Example~\ref{ex1}~($\Gamma = \Gamma_\rmD$).}
	\label{table:ex1}
%\small
	\begin{tabular}{cccccc}
		\toprule
		& & \multicolumn{4}{c}{LG1} \\ \cline{3-6}
		$N$ & $\Delta t$&$E_0(\eta)$  & EOC & $E_0(u)$ & EOC  \\ \hline\hline
		8 & $8.84\times 10^{-2}$ &$3.89\times 10^{0} $  &- & $3.78\times 10^{-2} $ &-  \\
		16 & $6.25\times 10^{-2}$ & $2.20\times 10^{0} $ & 1.65 & $2.28\times 10^{-2} $ & 1.45 \\
		32 & $4.42\times 10^{-2}$ & $1.45\times 10^{0} $ & 1.19& $1.57\times 10^{-2} $  & 1.09 \\
		64 & $3.13\times 10^{-2}$ &$1.01\times 10^{0} $ &1.05& $1.10\times 10^{-2} $ & 1.03 \\
		128 & $2.21 \times 10^{-2}$ & $7.11\times 10^{-1} $ & 1.01& $7.77\times 10^{-3} $ & 1.00 \\
		256 & $1.56\times 10^{-2}$ &$5.02\times 10^{-1} $ & 1.00& $5.51\times 10^{-3} $ & 0.99 \\
		\midrule
		& & \multicolumn{4}{c}{LG1} \\ \cline{3-6}
	$N$ &$\Delta t$& $E_1(\eta)$  & EOC & $E_1(u)$ & EOC \\ \hline\hline
		8 &$8.84\times 10^{-2}$& $3.00\times 10^{0} $  &- & $7.78\times 10^{-2} $ &-  \\
		16 &$6.25\times 10^{-2}$ & $1.73\times 10^{0} $ & 1.59 & $4.63\times 10^{-2} $ & 1.49\\
		32 &$4.42\times 10^{-2}$ & $1.25\times 10^{0} $ & 0.93& $2.95\times 10^{-2} $  & 1.31\\
		64 & $3.13\times 10^{-2}$ &$9.78\times 10^{-1} $ &0.71& $2.04\times 10^{-2} $ & 1.06 \\
		128 &$2.21 \times 10^{-2}$ & $6.42\times 10^{-1} $ & 1.22& $1.42\times 10^{-2} $ & 1.04 \\
		256 &$1.56\times 10^{-2}$& $4.35\times 10^{-1} $ & 1.12& $1.00\times 10^{-2} $ & 1.01\\
		\toprule
		& & \multicolumn{4}{c}{LG2} \\ \cline{3-6}
		$N$ & $\Delta t$&$E_0(\eta)$  & EOC & $E_0(u)$ & EOC \\ \hline\hline
		8 & $8.84\times 10^{-2}$&$6.81\times 10^{-1}$  &- & $1.71\times 10^{-2}$ &-  \\
		16 & $6.25\times 10^{-2}$ &$1.96\times 10^{-1}$ & 3.60 & $7.03\times 10^{-3}$ & 2.57 \\
		32 & $4.42\times 10^{-2}$&$8.53\times 10^{-2}$ & 2.40& $3.32\times 10^{-3}$  & 2.16 \\
		64 &$3.13\times 10^{-2}$&$3.82\times 10^{-2}$ &2.32& $1.64\times 10^{-3}$ & 2.04 \\
		128 &$2.21 \times 10^{-2}$ &$1.87\times 10^{-2}$ & 2.05& $8.20\times 10^{-4}$ & 1.99  \\
		256 &$1.56\times 10^{-2}$ &$9.46\times 10^{-3}$ & 1.97& $4.17\times 10^{-4}$ & 1.95 \\
		\midrule
		& & \multicolumn{4}{c}{LG2} \\ \cline{3-6}
		$N$ & $\Delta t$& $E_1(\eta)$  & EOC & $E_1(u)$ & EOC \\ \hline\hline
		8 &$8.84\times 10^{-2}$&$3.97\times 10^{0} $  &- & $5.68\times 10^{-2}$ &-  \\
		16 & $6.25\times 10^{-2}$&$2.24\times 10^{0} $ & 1.65 & $2.90\times 10^{-2}$ & 1.94  \\
		32 & $4.42\times 10^{-2}$&$2.00\times 10^{0} $ & 0.33& $1.20\times 10^{-2}$ & 2.54\\
		64 &$3.13\times 10^{-2}$&$1.64\times 10^{0} $ &0.57& $6.72\times 10^{-3}$ & 1.67\\
		128 & $2.21 \times 10^{-2}$&$1.17\times 10^{0} $ & 0.97& $3.23\times 10^{-3}$ & 2.11 \\
		256 &$1.56\times 10^{-2}$&$8.64\times 10^{-1} $ & 0.88& $1.47\times 10^{-3}$ & 2.28\\
		\bottomrule
	\end{tabular}
\end{table}
%
%%%%%%%%%%%%%%%%%%%%%%%%%%%%%%%%%%%%%%%%%%%%%%%%%%%%%%%%%%%%%%%%%%%%%%%%%%%%%%%%%%%%%%%%%
\begin{table}[!htbp]
	\centering
	\caption{Values of $E_i(\eta)$ and $E_i(u)$, $i=0,1$, by schemes LG1 and LG2 for Example~\ref{ex2}~($\Gamma = \overline{\Gamma}_\rmD \cup \overline{\Gamma}_\rmT$).}
	\label{table:ex2}
%\footnotesize
	\begin{tabular}{cccccc}
		\toprule
		& & \multicolumn{4}{c}{LG1} \\ \cline{3-6}
		$N$ & $\Delta t$& $E_0(\eta)$  & EOC & $E_0(u)$ & EOC  \\
		\hline\hline
		8& $8.84\times 10^{-2}$& $3.88\times 10^{0} $  &- & $3.86\times 10^{-2} $ &-  \\
		16 &$6.25\times 10^{-2}$& $2.19\times 10^{0} $ & 1.65 & $2.33\times 10^{-2} $ & 1.46  \\
		32 &$4.42\times 10^{-2}$& $1.45\times 10^{0} $ & 1.19& $1.58\times 10^{-2} $  & 1.11 \\
		64 &$3.13\times 10^{-2}$& $101\times 10^{0} $ &1.05& $1.11\times 10^{-2} $ & 1.03 \\
		128 &$2.21 \times 10^{-2}$& $7.09\times 10^{-1} $ & 1.01& $7.82\times 10^{-3} $ & 1.01 \\
		256& $1.56\times 10^{-2}$& $5.01\times 10^{-1} $ & 1.00& $5.53\times 10^{-3} $ & 1.00 \\
		\midrule
		& & \multicolumn{4}{c}{LG1} \\ \cline{3-6}
		$N$ & $\Delta t$& $E_1(\eta)$  & EOC & $E_1(u)$ & EOC  \\
		\hline\hline
		8 &$8.84\times 10^{-2}$& $2.95\times 10^{0} $  &- & $7.80\times 10^{-2} $ &-  \\
		16& $6.25\times 10^{-2}$& $1.71\times 10^{0} $ & 1.57 & $4.64\times 10^{-2} $ & 1.50 \\
		32 &$4.42\times 10^{-2}$& $1.24\times 10^{0} $ & 0.94& $2.95\times 10^{-2} $  & 1.31\\
		64 &$3.13\times 10^{-2}$& $9.78\times 10^{-1} $ &0.67& $2.03\times 10^{-2} $ & 1.07 \\
		128 &$2.21 \times 10^{-2}$& $6.42\times 10^{-1} $ & 1.21& $1.41\times 10^{-2} $ & 1.04 \\
		256 &$1.56\times 10^{-2}$& $4.34\times 10^{-1} $ & 1.13& $9.96\times 10^{-3} $ & 1.01\\
		\toprule
		& & \multicolumn{4}{c}{LG2} \\ \cline{3-6}
	$N$ & $\Delta t$& $E_0(\eta)$  & EOC & $E_0(u)$ & EOC  \\
		\hline\hline
		8 &$8.84\times 10^{-2}$ &$6.70\times 10^{-1} $  &- & $1.75\times 10^{-2} $ &-  \\
		16 & $6.25\times 10^{-2}$&$1.95\times 10^{-1} $ & 3.56 & $7.23\times 10^{-3} $ & 2.55 \\
		32 & $4.42\times 10^{-2}$&$8.58\times 10^{-2} $ & 2.37& $3.37\times 10^{-3} $  & 2.20 \\
		64 & $3.13\times 10^{-2}$&$3.97\times 10^{-2} $ &2.22& $1.67\times 10^{-3} $ & 2.03 \\
		128 &$2.21 \times 10^{-2}$& $1.87\times 10^{-2} $ & 2.17& $8.37\times 10^{-4} $ & 2.00 \\
		256 &$1.56\times 10^{-2}$& $9.54\times 10^{-3} $ & 1.94& $4.25\times 10^{-4} $ & 1.96 \\
		\midrule
		& & \multicolumn{4}{c}{LG2} \\ \cline{3-6}
		$N$ &$\Delta t$& $E_1(\eta)$  & EOC & $E_1(u)$ & EOC \\
		\hline\hline
		8 & $8.84\times 10^{-2}$& $3.89\times 10^{0} $  &- & $5.70\times 10^{-2} $ &-  \\
		16 &$6.25\times 10^{-2}$&$2.21\times 10^{0} $ & 1.63 & $2.93\times 10^{-2} $ & 1.92 \\
		32 &$4.42\times 10^{-2}$ &$1.98\times 10^{0} $ & 0.32& $1.24\times 10^{-2} $  & 2.49\\
		64 & $3.13\times 10^{-2}$&$1.65\times 10^{0} $ &0.54& $6.90\times 10^{-3} $ & 1.69 \\
		128 &$2.21 \times 10^{-2}$& $1.17\times 10^{0} $ & 0.97& $3.26\times 10^{-3} $ & 2.16 \\
		256 &$1.56\times 10^{-2}$& $8.62\times 10^{-1} $ & 0.89& $1.48\times 10^{-3} $ & 2.27\\
		\bottomrule
	\end{tabular}
\end{table}
%
%
%
%
%
%%%%%%%%%%%%%%%%%%%%%%%%%%%%%%%%%%%%%%%%%%%%%%%%%%%%%%%%%%%%%%%%%%%%%%%%%%%%%%%%%%%%%%%%%%%%%%%%%%%%%%%%%%%
%
%
%
%
%
%
%
%
%
%
%
%
\subsection{Effect of the TBC}
\label{sub:problem}
\begin{figure}[!htbp]
	\begin{center}
		\begin{tabular}{ccccc}
			\hspace{0.8mm}
			\subfloat[(i) t=0\label{D0}]{%
				\begin{overpic}[width=0.175\linewidth]{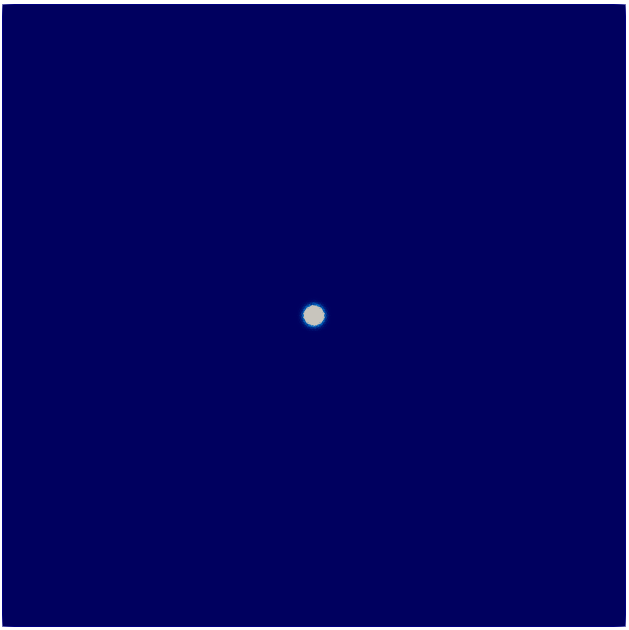}
					\put(-26.3,50){\large\textbf{(a)}}
					\put(45,87){$\color{white}\Gamma_\rmD$}
					\put(45,5){$\color{white}\Gamma_\rmD$}
					\put(80,45){$\color{white}\Gamma_\rmD$}
					\put(3,45){$\color{white}\Gamma_\rmD$}
				\end{overpic}
			}&
			\hspace{-5.3mm}
			\subfloat[(ii) t=25\label{D25}]{%
				\begin{overpic}[width=0.175\linewidth]{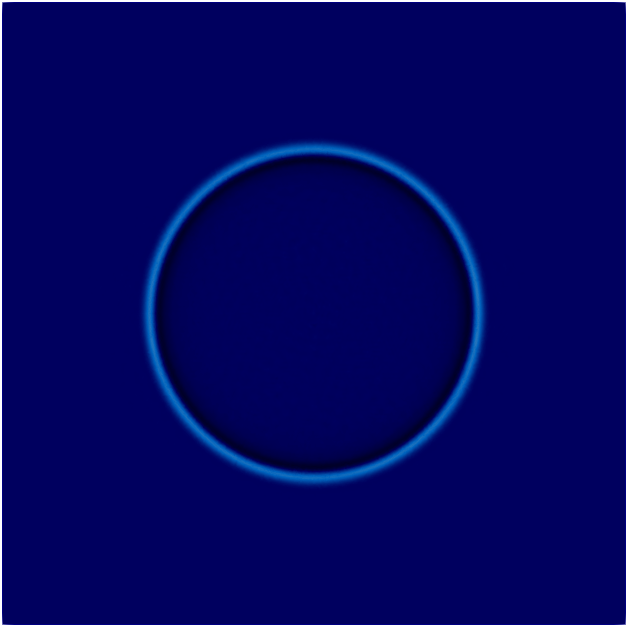}
					%	\put(40,40){some text again!}
					\put(45,87){$\color{white}\Gamma_\rmD$}
					\put(45,5){$\color{white}\Gamma_\rmD$}
					\put(80,45){$\color{white}\Gamma_\rmD$}
					\put(3,45){$\color{white}\Gamma_\rmD$}
				\end{overpic}
			}&
			\hspace{-5.3mm}
			\subfloat[(iii) t=50\label{D50}]{%
				\begin{overpic}[width=0.175\linewidth]{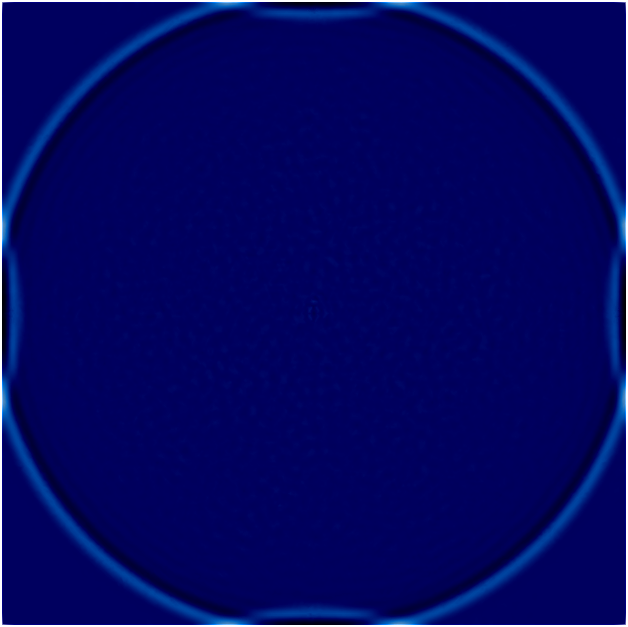}
					%	\put(-4.3,50){(a)}
					\put(45,87){$\color{white}\Gamma_\rmD$}
					\put(45,5){$\color{white}\Gamma_\rmD$}
					\put(80,45){$\color{white}\Gamma_\rmD$}
					\put(3,45){$\color{white}\Gamma_\rmD$}
				\end{overpic}
			}&
			\hspace{-5.3mm}
			\subfloat[(iv) t=75\label{D75}]{%
				\begin{overpic}[width=0.175\linewidth]{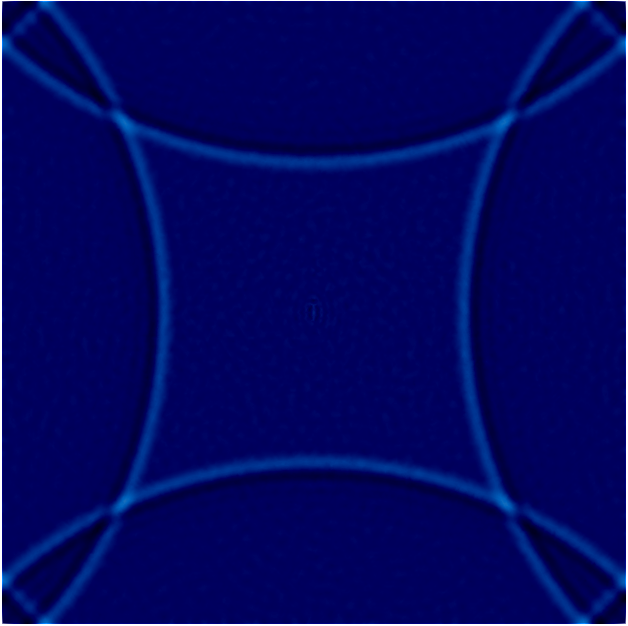}
					%	\put(40,40){some text again!}
					\put(45,87){$\color{white}\Gamma_\rmD$}
					\put(45,5){$\color{white}\Gamma_\rmD$}
					\put(80,45){$\color{white}\Gamma_\rmD$}
					\put(3,45){$\color{white}\Gamma_\rmD$}
				\end{overpic}
			}&
			\hspace{-5.3mm}
			\subfloat[(v) t=100\label{D100}]{%
				\begin{overpic}[width=2.8cm,height=2.1cm]{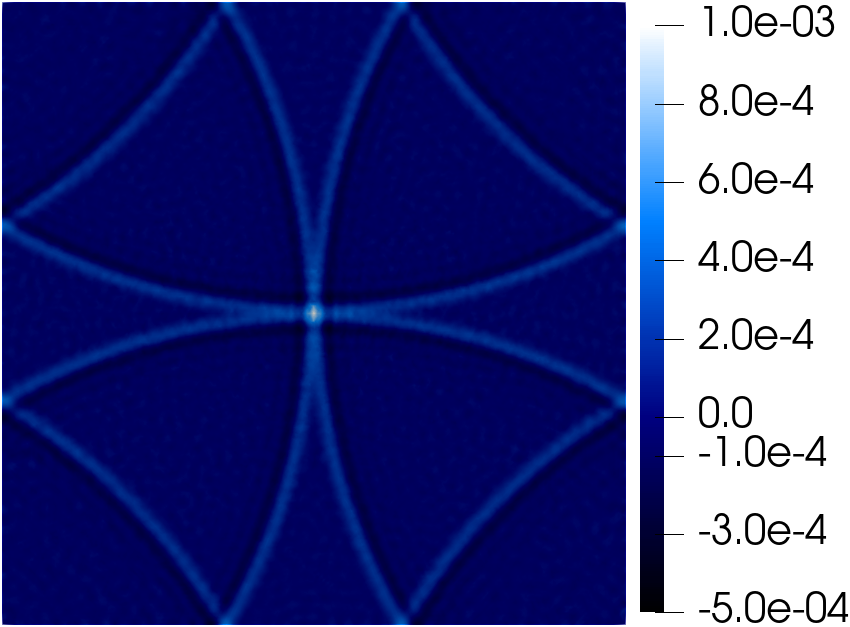}
					%\put(40,40){some text again!}
					\put(31,65){$\color{white}\Gamma_\rmD$}
					\put(31,2){$\color{white}\Gamma_\rmD$}
					\put(61,33){$\color{white}\Gamma_\rmD$}
					\put(1,33){$\color{white}\Gamma_\rmD$}
				\end{overpic}
			}	\\		
			\hspace{0.8mm}
			\subfloat[(i) t=0\label{1t0}]{%
				\begin{overpic}[width=0.175\linewidth]{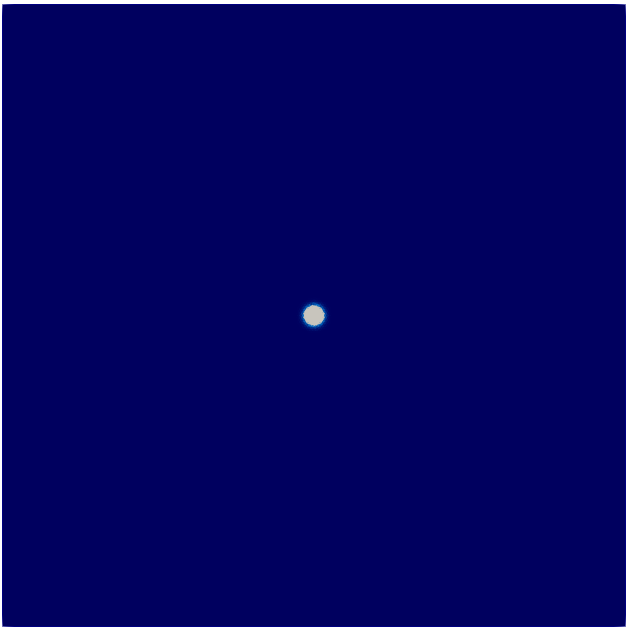}
					\put(-28.3,50){\large\textbf{(b)}}
					\put(45,87){$\color{white}\Gamma_\rmD$}
					\put(45,5){$\color{white}\Gamma_\rmT$}
					\put(80,45){$\color{white}\Gamma_\rmD$}
					\put(3,45){$\color{white}\Gamma_\rmD$}
				\end{overpic}
			}&
			\hspace{-5.3mm}
			\subfloat[(ii) t=25\label{1t25}]{%
				\begin{overpic}[width=0.175\linewidth]{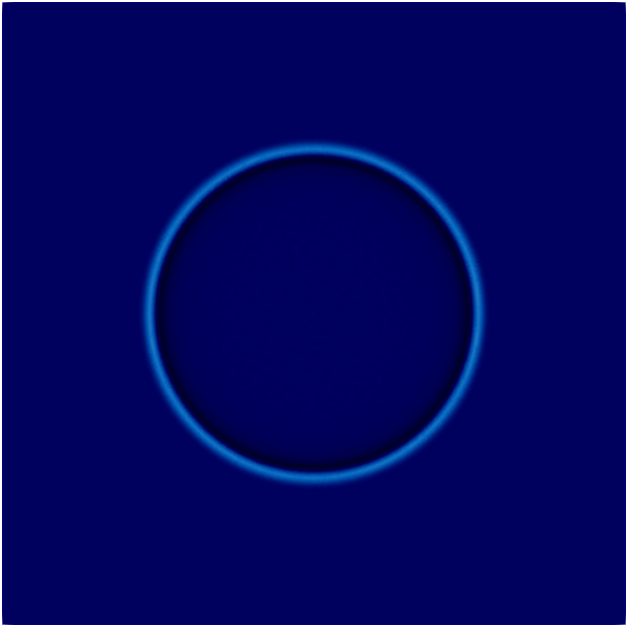}
					%	\put(40,40){some text again!}
					\put(45,87){$\color{white}\Gamma_\rmD$}
					\put(45,5){$\color{white}\Gamma_\rmT$}
					\put(80,45){$\color{white}\Gamma_\rmD$}
					\put(3,45){$\color{white}\Gamma_\rmD$}
				\end{overpic}
			}&
			\hspace{-5.3mm}
			\subfloat[(iii) t=50\label{1t50}]{%
				\begin{overpic}[width=0.175\linewidth]{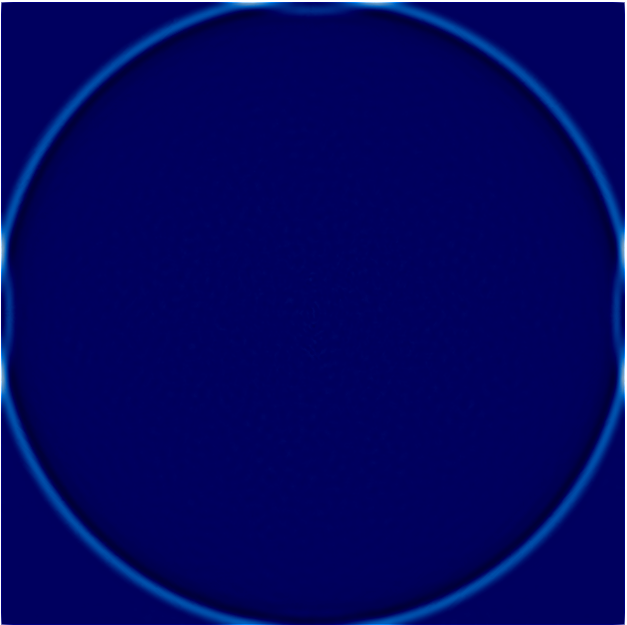}
					%	\put(-4.3,50){(a)}
					\put(45,87){$\color{white}\Gamma_\rmD$}
					\put(45,5){$\color{white}\Gamma_\rmT$}
					\put(80,45){$\color{white}\Gamma_\rmD$}
					\put(3,45){$\color{white}\Gamma_\rmD$}
				\end{overpic}
			}&
			\hspace{-5.3mm}
			\subfloat[(iv) t=75\label{1t75}]{%
				\begin{overpic}[width=0.175\linewidth]{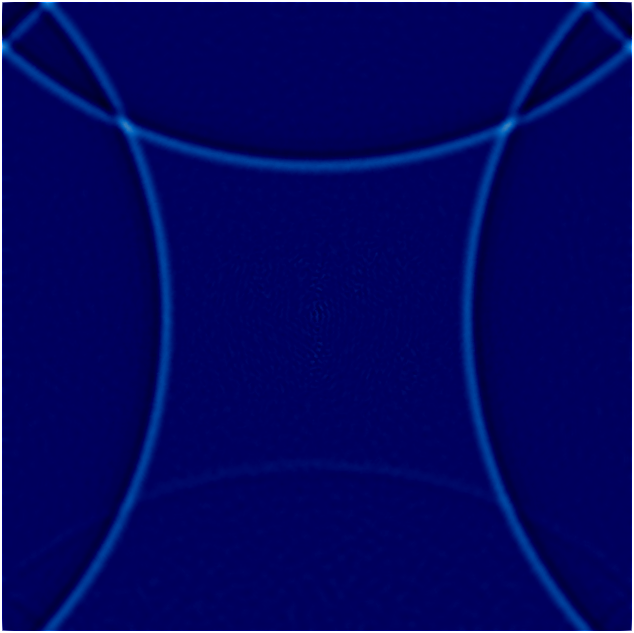}
					%	\put(40,40){some text again!}
					\put(45,87){$\color{white}\Gamma_\rmD$}
					\put(45,5){$\color{white}\Gamma_\rmT$}
					\put(80,45){$\color{white}\Gamma_\rmD$}
					\put(3,45){$\color{white}\Gamma_\rmD$}
				\end{overpic}
			}&
			\hspace{-5.3mm}
			\subfloat[(v) t=100\label{1t100}]{%
				\begin{overpic}[width=2.8cm,height=2.1cm]{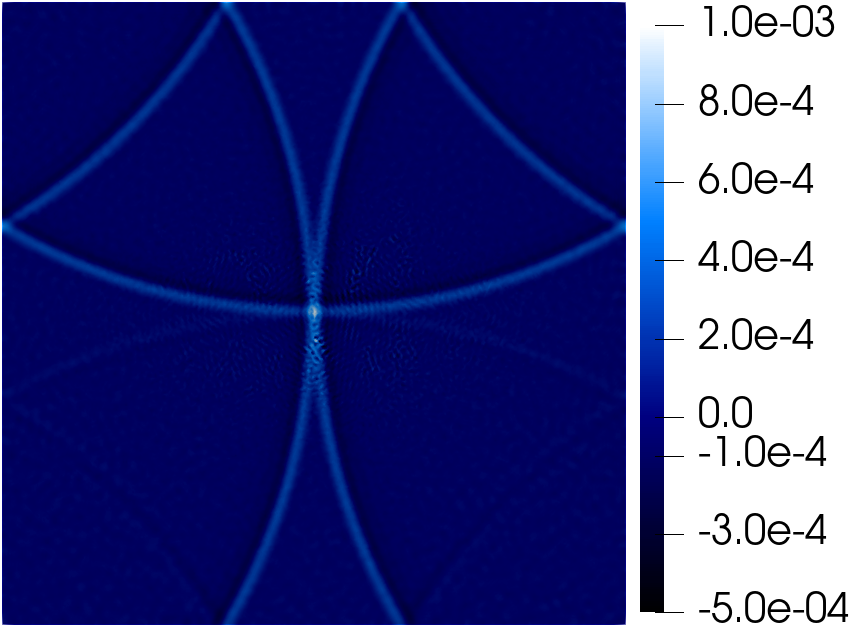}
					%\put(40,40){some text again!}
					\put(31,65){$\color{white}\Gamma_\rmD$}
					\put(31,2){$\color{white}\Gamma_\rmT$}
					\put(61,33){$\color{white}\Gamma_\rmD$}
					\put(1,33){$\color{white}\Gamma_\rmD$}
				\end{overpic}
			}	\\	
			\hspace{0.8mm}
			\subfloat[(i) t=0\label{12t0}]{%
				\begin{overpic}[width=0.175\linewidth]{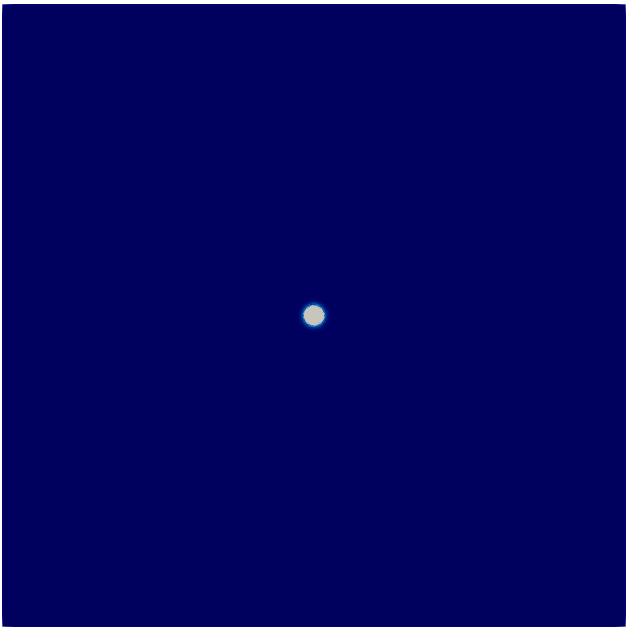}
					\put(-26.3,50){\large\textbf{(c)}}
					\put(45,87){$\color{white}\Gamma_\rmD$}
					\put(45,5){$\color{white}\Gamma_\rmT$}
					\put(80,45){$\color{white}\Gamma_\rmT$}
					\put(3,45){$\color{white}\Gamma_\rmD$}
				\end{overpic}
			}	&
			\hspace{-5.3mm}
			\subfloat[(ii) t=25\label{12t25}]{%
				\begin{overpic}[width=0.175\linewidth]{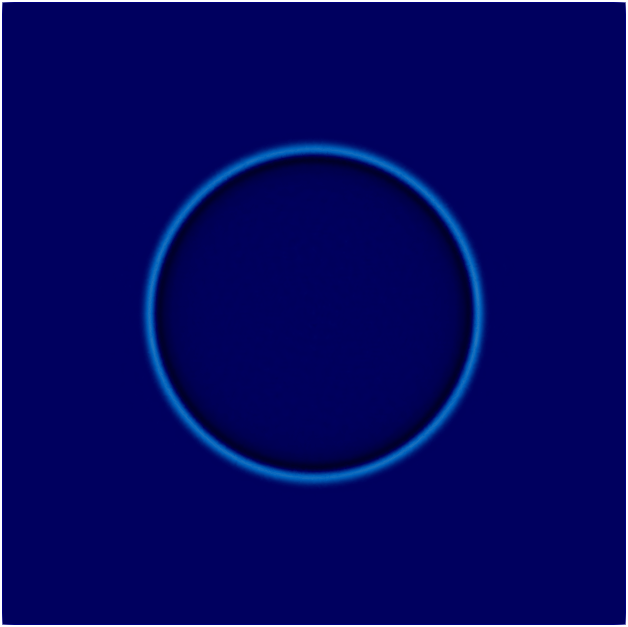}
					%	\put(40,40){some text again!}
					\put(45,87){$\color{white}\Gamma_\rmD$}
					\put(45,5){$\color{white}\Gamma_\rmT$}
					\put(80,45){$\color{white}\Gamma_\rmT$}
					\put(3,45){$\color{white}\Gamma_\rmD$}
				\end{overpic}
			}	&
			\hspace{-5.3mm}
			\subfloat[(iii) t=50\label{12t50}]{%
				\begin{overpic}[width=0.175\linewidth]{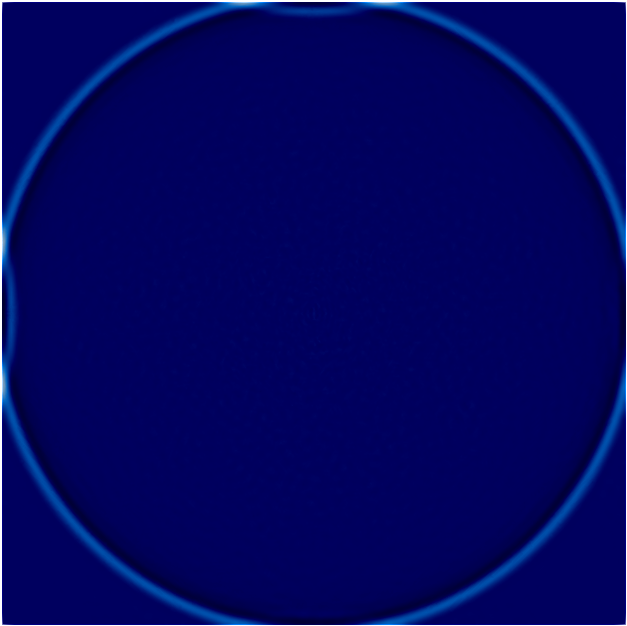}
					%	\put(-4.3,50){(a)}
					\put(45,87){$\color{white}\Gamma_\rmD$}
					\put(45,5){$\color{white}\Gamma_\rmT$}
					\put(80,45){$\color{white}\Gamma_\rmT$}
					\put(3,45){$\color{white}\Gamma_\rmD$}
				\end{overpic}
			}	&
			\hspace{-5.3mm}
			\subfloat[(iv) t=75\label{12t75}]{%
				\begin{overpic}[width=0.175\linewidth]{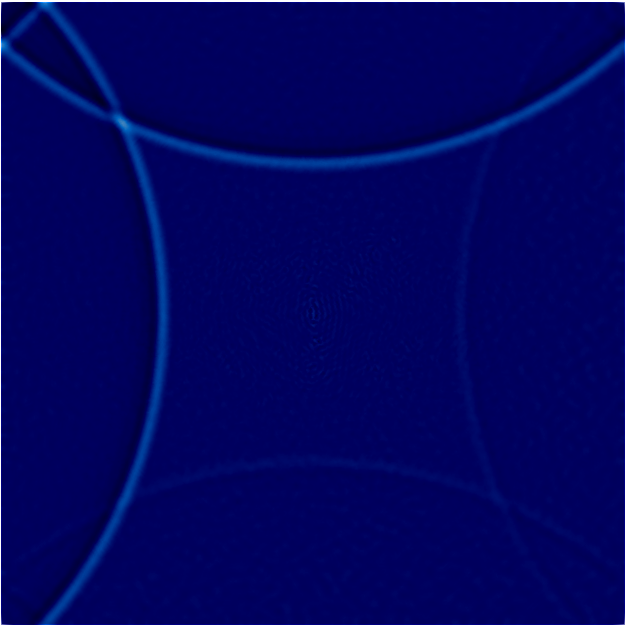}
					%	\put(40,40){some text again!}
					\put(45,87){$\color{white}\Gamma_\rmD$}
					\put(45,5){$\color{white}\Gamma_\rmT$}
					\put(80,45){$\color{white}\Gamma_\rmT$}
					\put(3,45){$\color{white}\Gamma_\rmD$}
				\end{overpic}
			}	&
			\hspace{-5.3mm}
			\subfloat[(v) t=100\label{12t100}]{%
				\begin{overpic}[width=2.8cm,height=2.1cm]{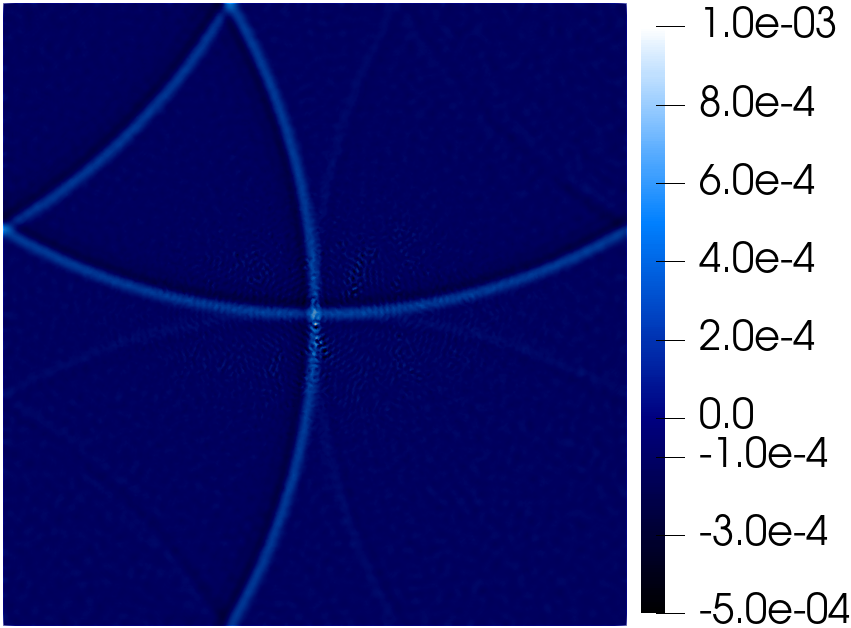}
					%\put(40,40){some text again!}
					\put(31,65){$\color{white}\Gamma_\rmD$}
					\put(31,2){$\color{white}\Gamma_\rmT$}
					\put(61,33){$\color{white}\Gamma_\rmT$}
					\put(1,33){$\color{white}\Gamma_\rmD$}
				\end{overpic}
			}	\\
			\hspace{0.8mm}
			\subfloat[(i) t=0\label{123t0}]{%
				\begin{overpic}[width=0.175\linewidth]{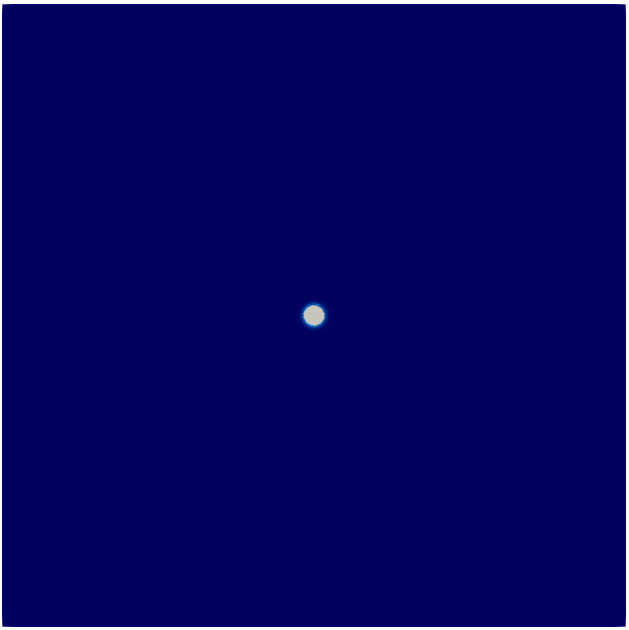}
					\put(-26.3,50){\large\textbf{(d)}}
					\put(45,87){$\color{white}\Gamma_\rmT$}
					\put(45,5){$\color{white}\Gamma_\rmT$}
					\put(80,45){$\color{white}\Gamma_\rmT$}
					\put(3,45){$\color{white}\Gamma_\rmD$}
				\end{overpic}
			}
			&
			\hspace{-5.3mm}
			\subfloat[(ii) t=25\label{123t25}]{%
				\begin{overpic}[width=0.175\linewidth]{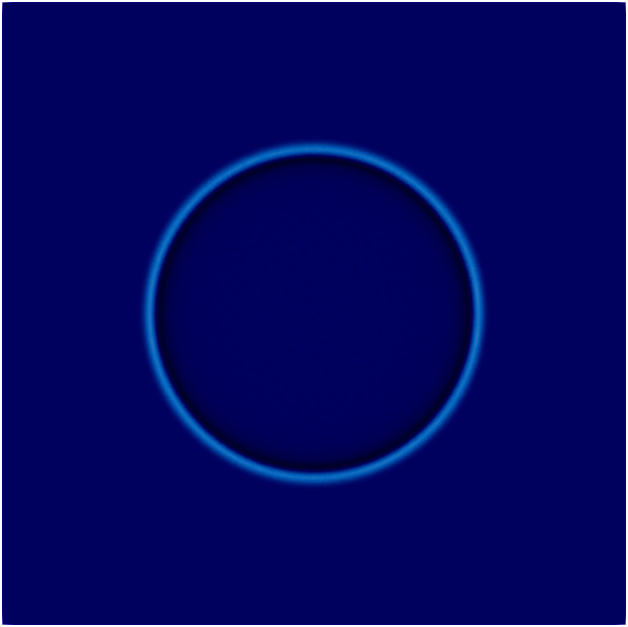}
					%	\put(40,40){some text again!}
					\put(45,87){$\color{white}\Gamma_\rmT$}
					\put(45,5){$\color{white}\Gamma_\rmT$}
					\put(80,45){$\color{white}\Gamma_\rmT$}
					\put(3,45){$\color{white}\Gamma_\rmD$}
				\end{overpic}
			}&
			\hspace{-5.3mm}
			\subfloat[(iii) t=50\label{123t50}]{%
				\begin{overpic}[width=0.175\linewidth]{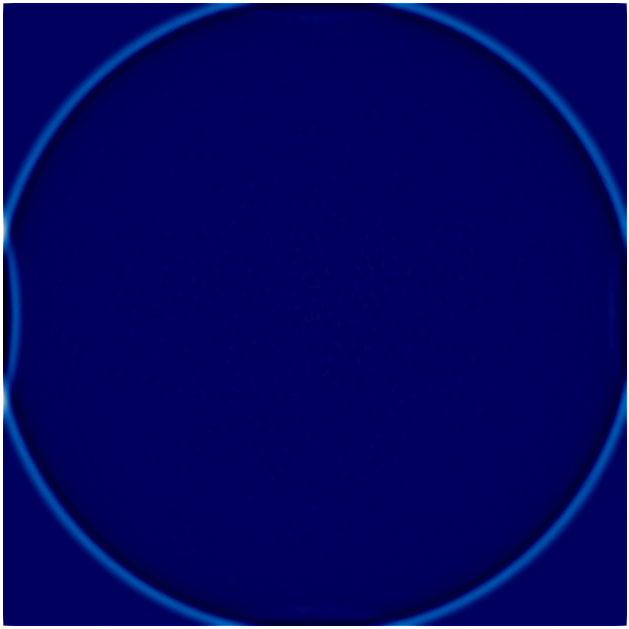}
					%	\put(-4.3,50){(a)}
					\put(45,87){$\color{white}\Gamma_\rmT$}
					\put(45,5){$\color{white}\Gamma_\rmT$}
					\put(80,45){$\color{white}\Gamma_\rmT$}
					\put(3,45){$\color{white}\Gamma_\rmD$}
				\end{overpic}
			}&
			\hspace{-5.3mm}
			\subfloat[(iv) t=75\label{123t75}]{%
				\begin{overpic}[width=0.175\linewidth]{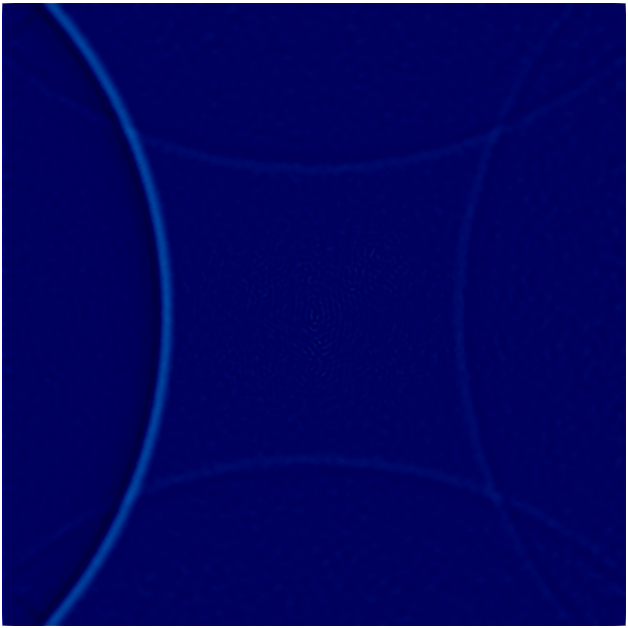}
					%	\put(40,40){some text again!}
					\put(45,87){$\color{white}\Gamma_\rmT$}
					\put(45,5){$\color{white}\Gamma_\rmT$}
					\put(80,45){$\color{white}\Gamma_\rmT$}
					\put(3,45){$\color{white}\Gamma_\rmD$}
				\end{overpic}
			}&
			\hspace{-5.3mm}
			\subfloat[(v) t=100\label{123t100}]{%
				\begin{overpic}[width=2.8cm,height=2.1cm]{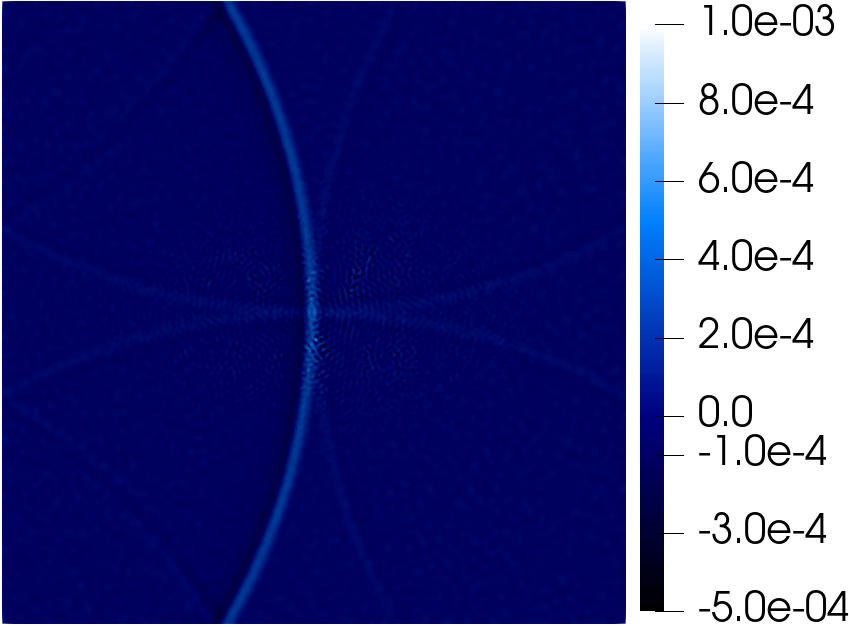}
					%\put(40,40){some text again!}
					\put(31,65){$\color{white}\Gamma_\rmT$}
					\put(31,2){$\color{white}\Gamma_\rmT$}
					\put(61,33){$\color{white}\Gamma_\rmT$}
					\put(1,33){$\color{white}\Gamma_\rmD$}
				\end{overpic}
			}	\\
			\hspace{0.8mm}
			\subfloat[(i) t=0\label{1234t0}]{%
				\begin{overpic}[width=0.175\linewidth]{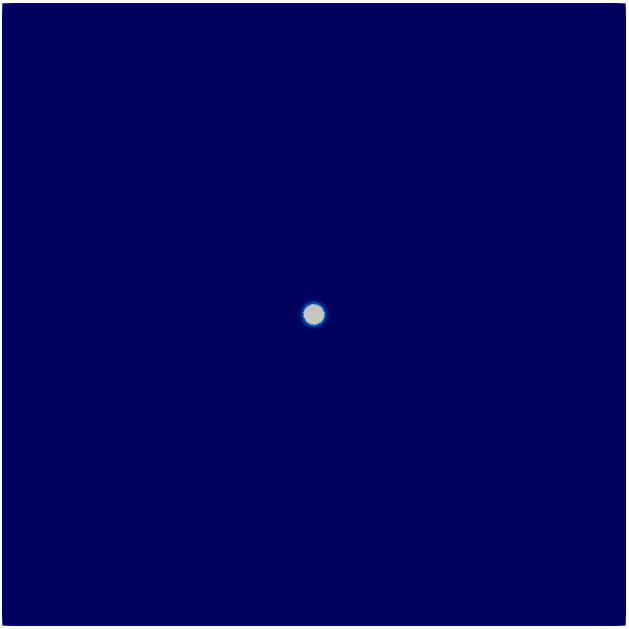}
					\put(-26.3,50){\large\textbf{(e)}}
					\put(45,87){$\color{white}\Gamma_\rmT$}
					\put(45,5){$\color{white}\Gamma_\rmT$}
					\put(80,45){$\color{white}\Gamma_\rmT$}
					\put(3,45){$\color{white}\Gamma_\rmT$}
				\end{overpic}
			}&
			\hspace{-5.3mm}
			\subfloat[(ii) t=25\label{1234t25}]{%
				\begin{overpic}[width=0.175\linewidth]{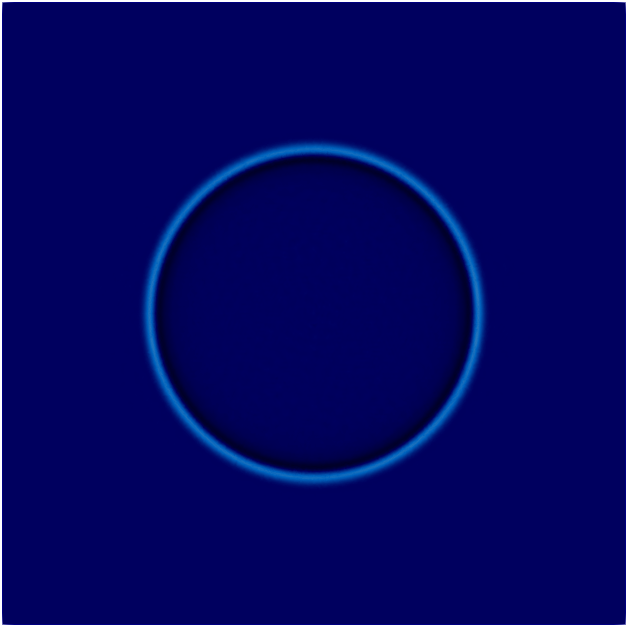}
					%	\put(40,40){some text again!}
					\put(45,87){$\color{white}\Gamma_\rmT$}
					\put(45,5){$\color{white}\Gamma_\rmT$}
					\put(80,45){$\color{white}\Gamma_\rmT$}
					\put(3,45){$\color{white}\Gamma_\rmT$}
				\end{overpic}
			}&
			\hspace{-5.3mm}
			\subfloat[(ii) t=50\label{1234t50}]{%
				\begin{overpic}[width=0.175\linewidth]{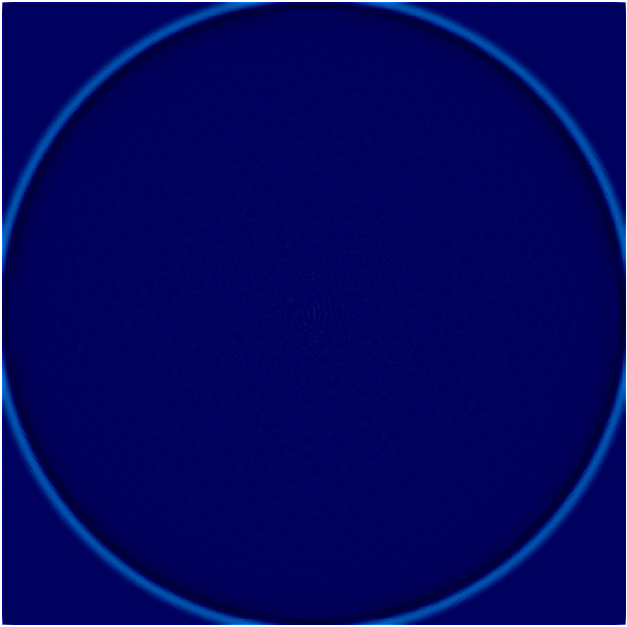}
					%	\put(-4.3,50){(a)}
					\put(45,87){$\color{white}\Gamma_\rmT$}
					\put(45,5){$\color{white}\Gamma_\rmT$}
					\put(80,45){$\color{white}\Gamma_\rmT$}
					\put(3,45){$\color{white}\Gamma_\rmT$}
				\end{overpic}
			}&
			\hspace{-5.3mm}
			\subfloat[(iv) t=75\label{1234t75}]{%
				\begin{overpic}[width=0.175\linewidth]{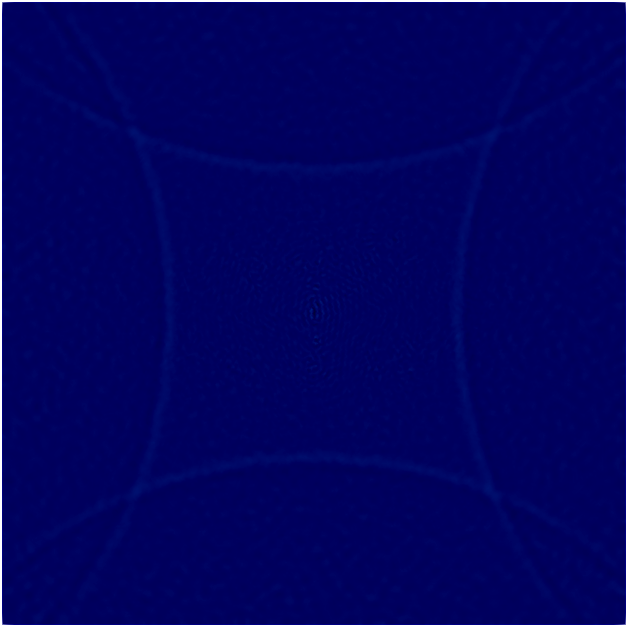}
					%	\put(40,40){some text again!}
					\put(45,87){$\color{white}\Gamma_\rmT$}
					\put(45,5){$\color{white}\Gamma_\rmT$}
					\put(80,45){$\color{white}\Gamma_\rmT$}
					\put(3,45){$\color{white}\Gamma_\rmT$}
				\end{overpic}
				
			}&
			\hspace{-5.3mm}
			\subfloat[(v) t=100\label{1234t100}]{%
				\begin{overpic}[width=2.8cm,height=2.1cm]{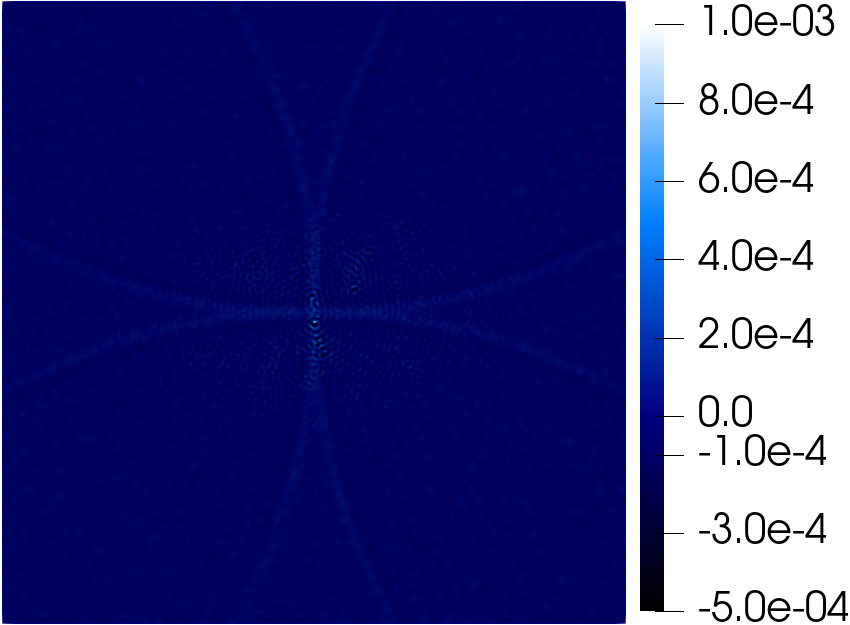}
					%\put(40,40){some text again!}
					\put(31,65){$\color{white}\Gamma_\rmT$}
					\put(31,2){$\color{white}\Gamma_\rmT$}
					\put(61,33){$\color{white}\Gamma_\rmT$}
					\put(1,33){$\color{white}\Gamma_\rmT$}
				\end{overpic}
			}\\
		\end{tabular}
		\caption{Color contours of $\eta_h^n$ by LG2 with and without the TBC for the five cases, (a)-(e), in Example~\ref{ex3}.}
		\label{fig:tm}
	\end{center}
\end{figure}
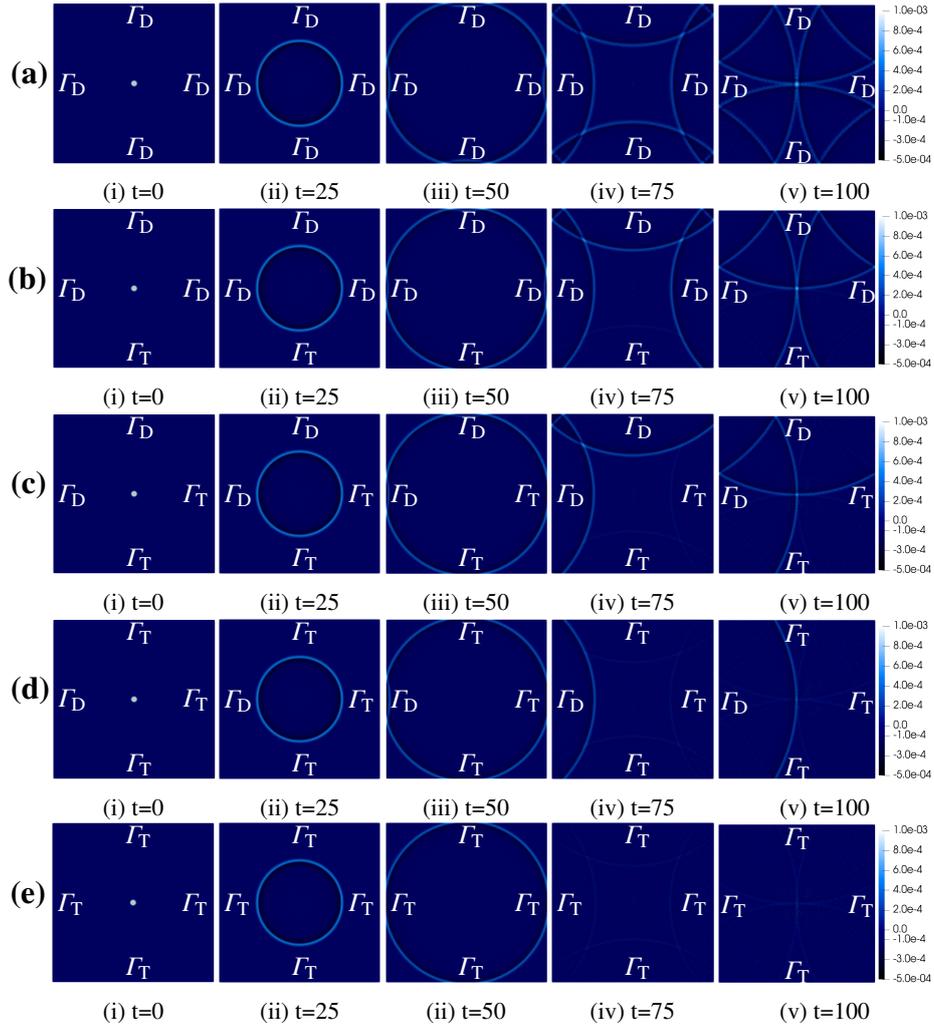

We consider the following example to see the effect of the TBC.
\begin{example}\label{ex3}
In problem~\eqref{eqn1}, we set $\Omega=(0,10)^{2}$, $T=100$, $g=\rho=\mu=\zeta=1$, $(f, F)=(0, 0)$, $\eta^{0} = c \exp( -100\,\vert x-p \vert^2 )$, $c=10^{-3}$, $p=(5, 5)^\top$, and~$u^{0}=0$.
We consider five cases of $\Gamma_\rmT$, \\
\indent {\rm (a)}~$\Gamma_\rmT = \emptyset$, i.e., $\Gamma = \Gamma_\rmD$,  \\
\indent {\rm (b)}~$\Gamma_\rmT = \{x\in\Gamma;\ x_2=0\}$~(bottom), $\Gamma_\rmD = \Gamma \setminus \overline{\Gamma}_\rmT$, \\
\indent {\rm (c)}~$\Gamma_\rmT = \{x\in\Gamma;\ x_1=10, x_2=0\}$~(right and bottom), $\Gamma_\rmD = \Gamma \setminus \overline{\Gamma}_\rmT$,  \\
\indent {\rm (d)}~$\Gamma_\rmT = \{x\in\Gamma;\ x_1=10, x_2=0, 10 \}$~(right, bottom and top), $\Gamma_\rmD = \Gamma \setminus \overline{\Gamma}_\rmT$,  \\
\indent {\rm (e)}~$\Gamma_\rmT = \Gamma$.
\end{example}
We solve Example~\ref{ex3} by LG2.
Figure~\ref{fig:tm} shows the color contours of~$\eta_{h}^{n}$ for~$t = 25k$, $k=0,\ldots, 4$, cf. (i)-(v), for the five cases, (a)-(e).
We can see the effect of the boundary conditions; 
the artificial reflection is observed and removed significantly when the wave touches the Dirichlet~$(\Gamma_\rmD)$ and the transmission~$(\Gamma_\rmT)$ boundaries, respectively.
Thus, LG2 works well for the SWEs with and without the TBC in the simple square domain.
\section{Application to the Bay of Bengal}
\label{sec:application}
\begin{figure}[!htbp]
	\centering
	\includegraphics[width=.55\textwidth]{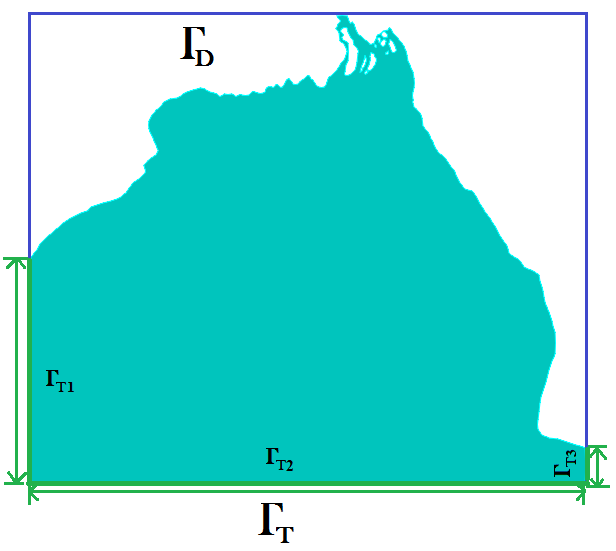}
	\caption{The domain for the Bay of Bengal region with the information of boundaries, $\Gamma_\rmD$ and~$\Gamma_\rmT~(=\Gamma_{\rmT1}\cup\Gamma_{\rmT2}\cup\Gamma_{\rmT3})$ used in Example~\ref{ex4}.}\label{Dmain}
\end{figure}

In this section, we apply LG2, i.e., scheme~\eqref{scheme} discussed in Subsection~\ref{subsection:scheme}, to a computational domain of the Bay of Bengal region, cf. Figure~\ref{Dmain}, which is an approximate domain of the original, cf. Figure~\ref{Bay}.
All the computations are performed via FreeFem++~\cite{MR3043640}.
\subsection{Numerical simulation with and without TBC}
\label{subsec:with_and_without_TBC}
We set the following example.
\begin{example}
\label{ex4}
Let $\Omega$ be the domain shown in Figure~\ref{Dmain}.
The domain is considered from $0$ to $1051.4~[\si{km}]$ in the horizontal direction and $0$ to $889.59~[\si{km}]$ in the vertical direction.
We employ two boundary conditions, the Dirichlet boundary condition on $\Gamma_\rmD$ and the TBC on $\Gamma_\rmT$, cf. Figure~\ref{Dmain}.
We set~$\Gamma_\rmD$ on the coastal and island boundaries and $\Gamma_\rmT$ on the artificial boundaries for the open sea.
As shown in Figure~\ref{Dmain}, there are three artificial boundaries on the open sea, i.e., $\Gamma_\rmT  = \Gamma_{\rmT 1}\cup\Gamma_{\rmT 2}\cup\Gamma_{\rmT 3}$.
In problem~\eqref{eqn1}, we set $T = 5{,}000~[\si{s}]$, $\zeta=2~[\si{km}]$, $\eta^{0}(x) = c_{1} \exp(-0.04\vert x-p\vert^2)~[\si{km}]$, $c_{1}=0.01~[\si{-}]$, $p = (559.56, 430.02)^\top$, $u^0 = 0$, $\mu=1~[\si{Pa.s}]$, $\rho=10^{12}~[\si{kg/km^3}]$, $g=9.8\times 10^{-3}~[\si{km/s^2}]$ and $(f, F) = (0, 0)$.
\end{example}
\par
\begin{figure}[H]
\centering
\includegraphics[width=.55\textwidth]{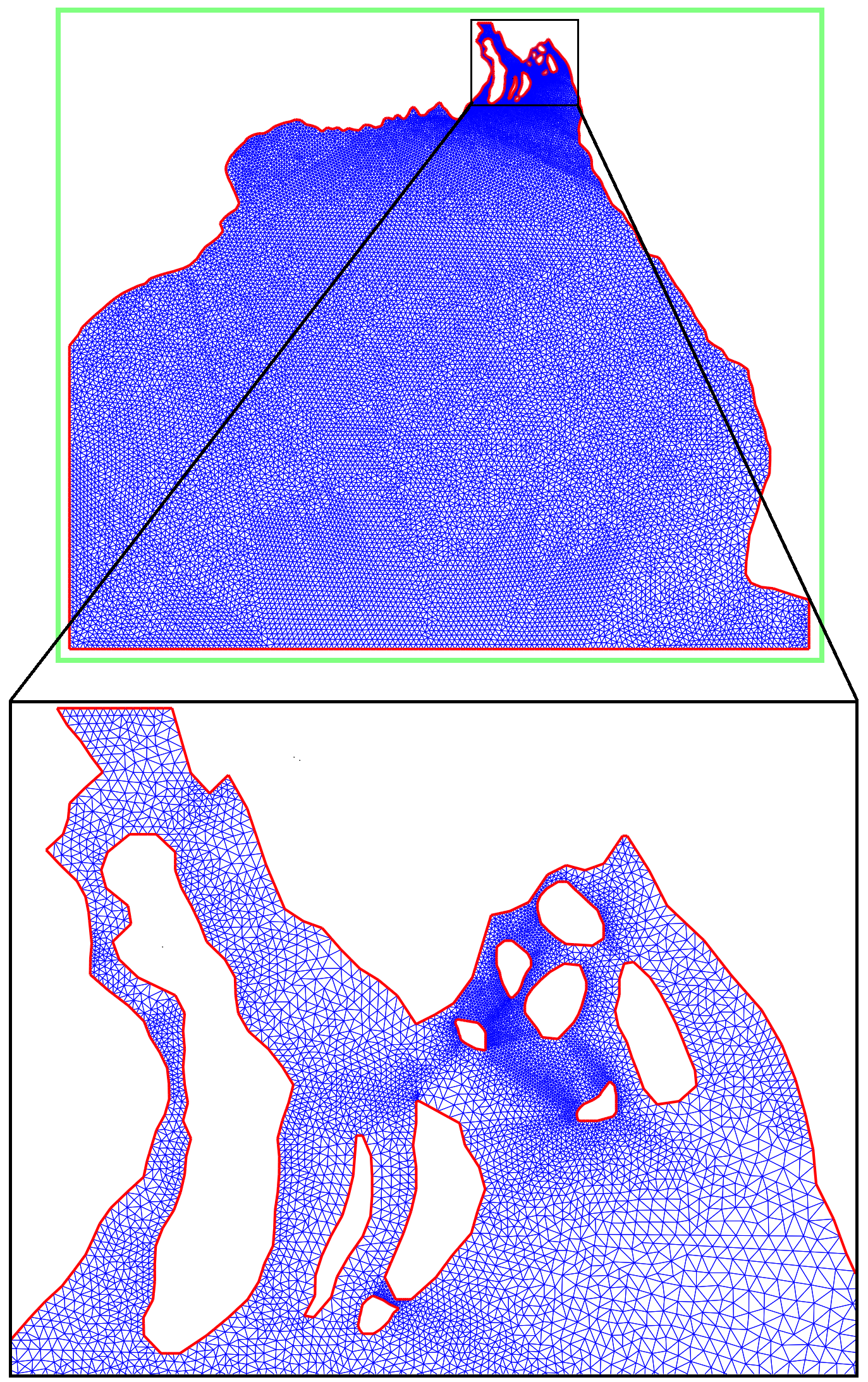}
\caption{The mesh for the Bay of Bengal region used for Example~\ref{ex4}.}\label{mesh}
\end{figure}
\begin{figure}[!htbp]
	\centering
	\begin{subfigure}[b]{0.45\textwidth}
		\centering
		\includegraphics[width=\textwidth]{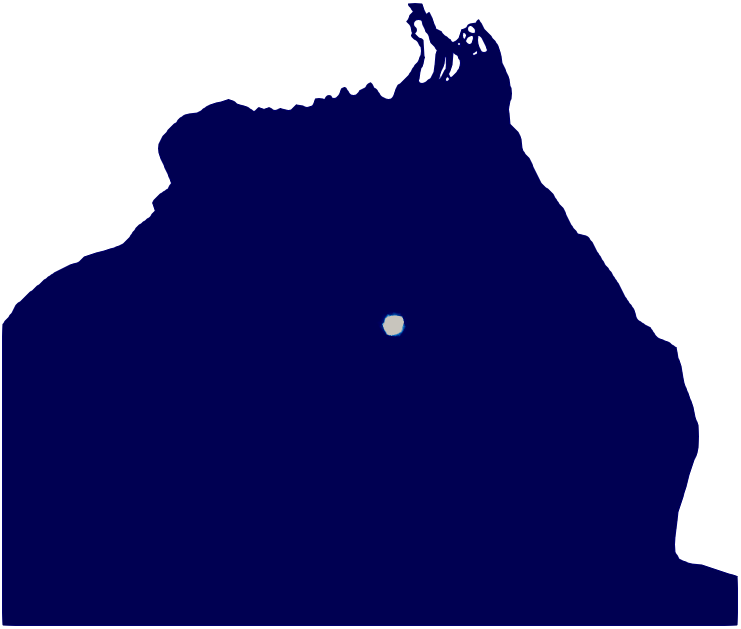}
		\caption{(a1) $ t=0$}\label{fig:b1}
	\end{subfigure}
	\begin{subfigure}[b]{0.45\textwidth}
		\centering
		\includegraphics[width=\textwidth]{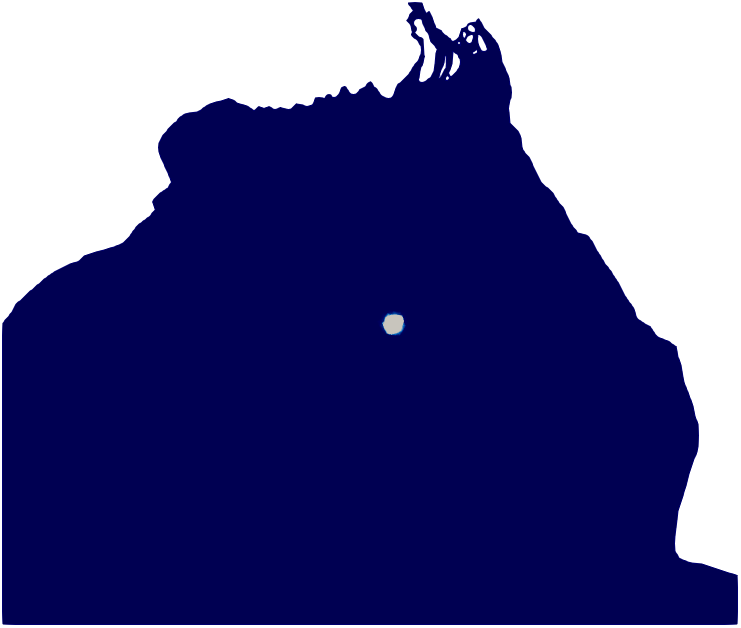}
		\caption{(a2) $t=0$}\label{fig:b2}
	\end{subfigure}
	\hspace*{0.08\textwidth} \\
	\begin{subfigure}[b]{0.45\textwidth}
		\centering
		\includegraphics[width=\textwidth]{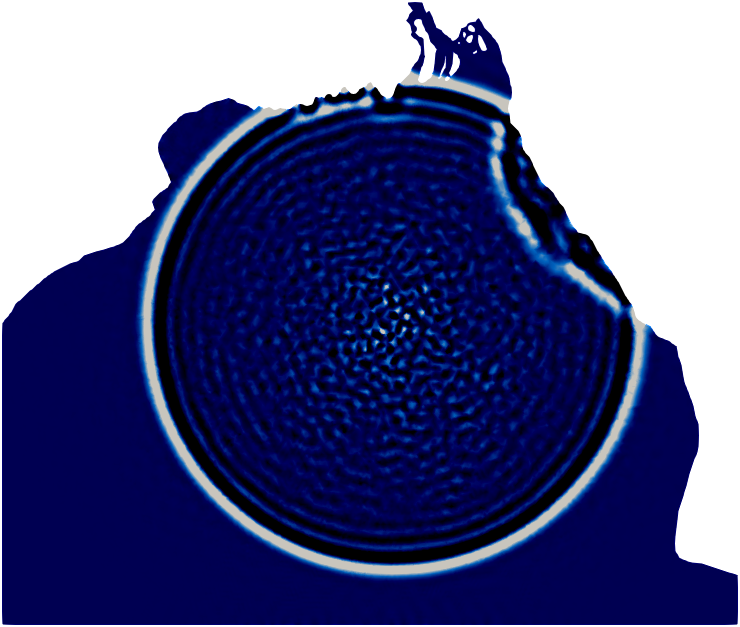}
		\caption{(b1) $t=2{,}500$}\label{fig:b3}
	\end{subfigure}
	\begin{subfigure}[b]{0.45\textwidth}
		\centering
		\includegraphics[width=\textwidth]{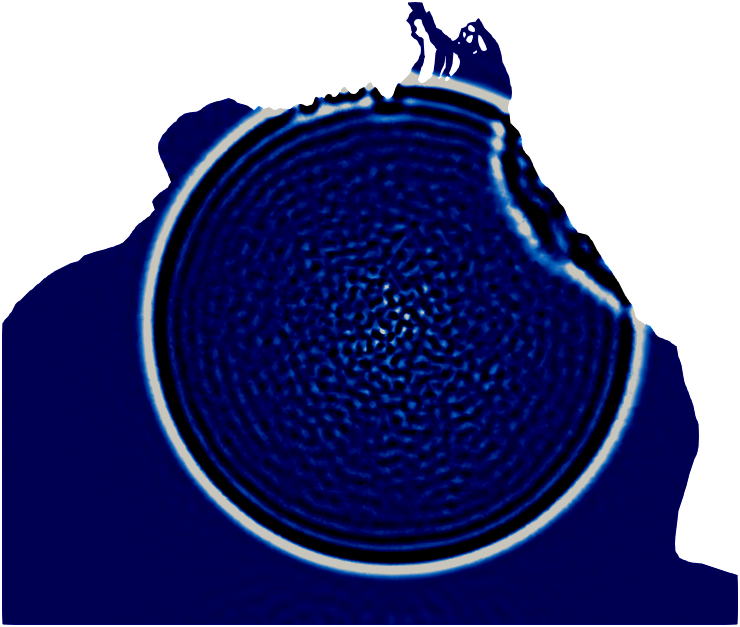}
		\caption{(b2) $t=2{,}500$}\label{fig:b4}
	\end{subfigure} 
	\hspace*{0.08\textwidth} \\
	\begin{subfigure}[b]{0.45\textwidth}
		\centering
		\includegraphics[width=\textwidth]{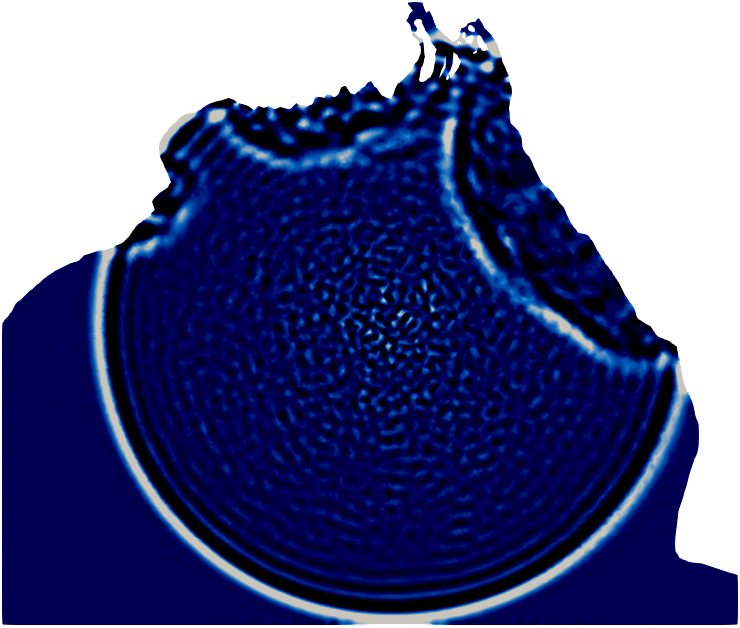}
		\caption{(c1) $t=3{,}000$}\label{fig:b5}
	\end{subfigure}
	\begin{subfigure}[b]{0.45\textwidth}
		\centering
		\includegraphics[width=\textwidth]{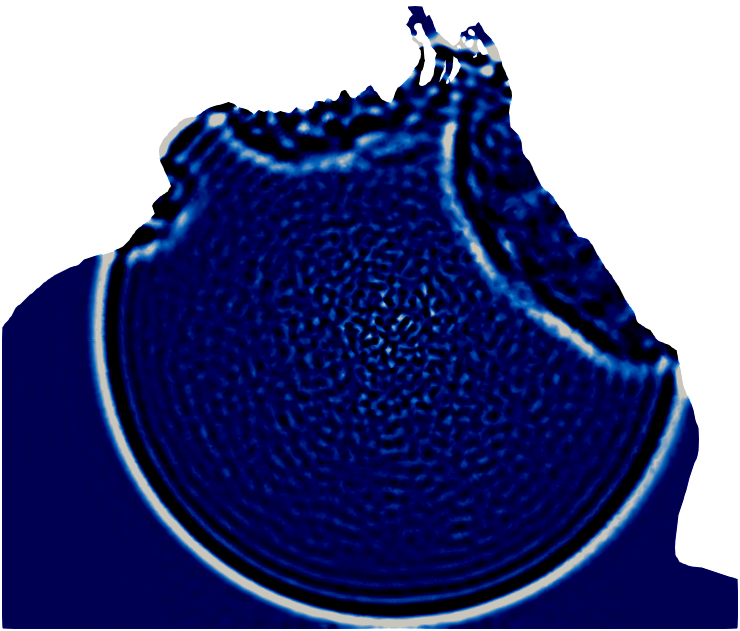}
		\caption{(c2) $t=3{,}000$}\label{figb6}
	\end{subfigure}
	\begin{subfigure}[b]{0.08\textwidth}
		\centering
		\includegraphics[width=\textwidth]{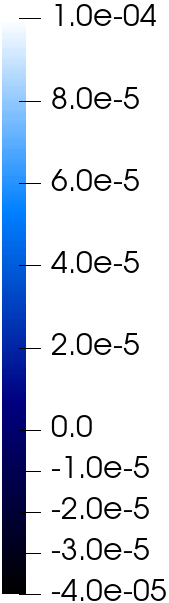}
		\caption{}
	\end{subfigure}
	\caption{Contour plot of $\eta_h^n$ by LG2 with $\Gamma = \Gamma_\rmD$~(left) and $\Gamma = \Gamma_\rmD \cup \Gamma_\rmT$~(right) on the Bay of Bengal for~$t=0, 2{,}500$ and $3{,}000$.}\label{fig:bob_part1}
\end{figure}
\begin{figure}[!htbp]
	\centering
	\begin{subfigure}[b]{0.45\textwidth}
		\centering
		\includegraphics[width=\textwidth]{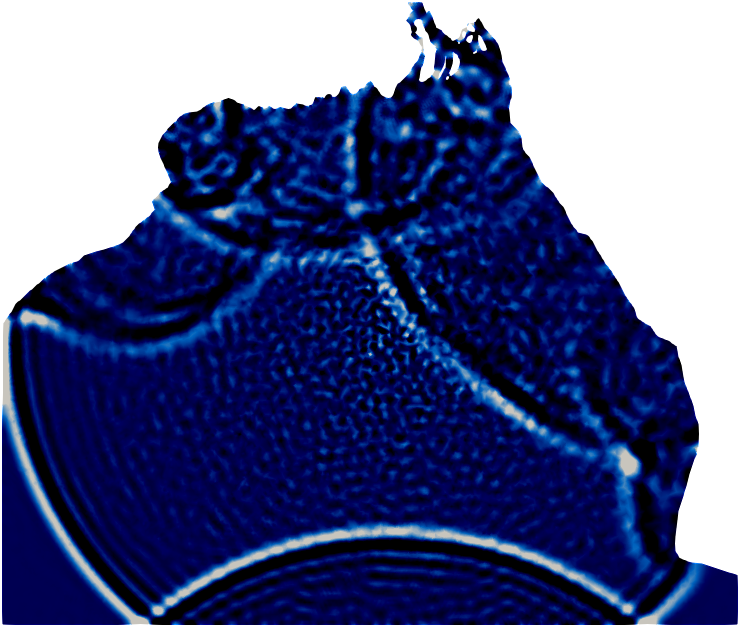}
		\caption{(d1) $t=4{,}000$}\label{fig:i1}
	\end{subfigure}
	\begin{subfigure}[b]{0.45\textwidth}
		\centering
		\includegraphics[width=\textwidth]{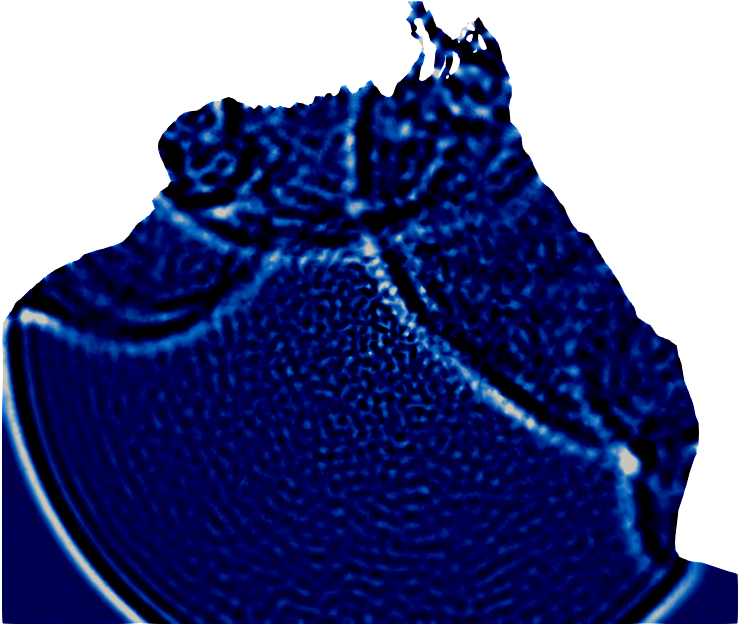}
		\caption{(d2) $t=4{,}000$}\label{fig:i2}
	\end{subfigure}
	\hspace*{0.08\textwidth} \\
	\begin{subfigure}[b]{0.45\textwidth}
		\centering
		\includegraphics[width=\textwidth]{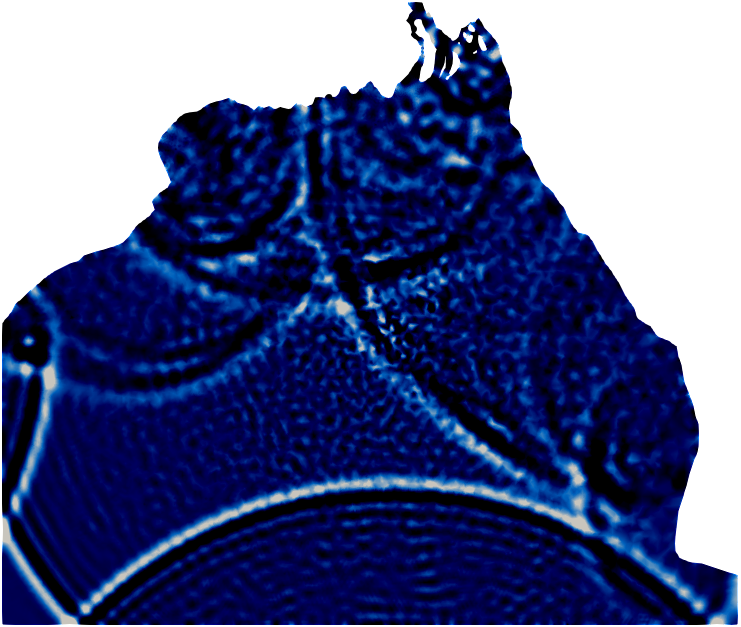}
		\caption{(e1) $t=4{,}500$}\label{fig:i3}
	\end{subfigure}
	\begin{subfigure}[b]{0.45\textwidth}
		\centering
		\includegraphics[width=\textwidth]{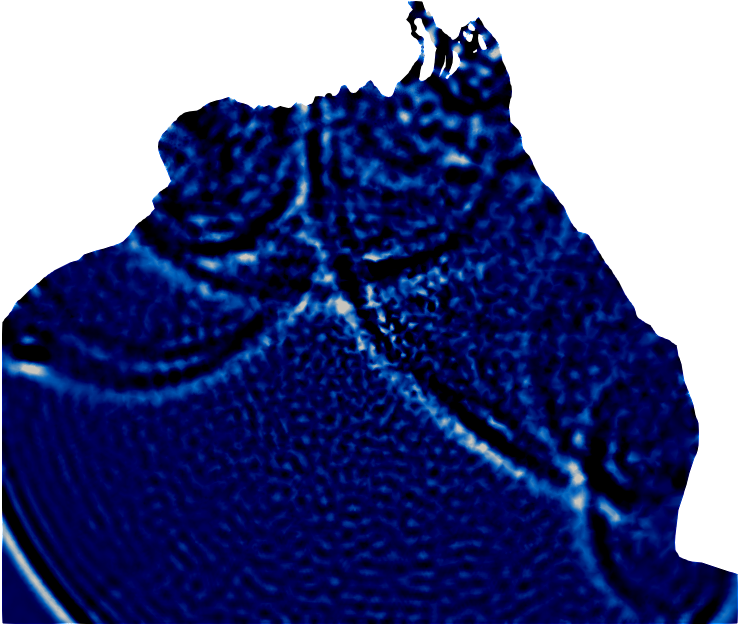}
		\caption{(e2) $t=4{,}500$}\label{fig:i4}
	\end{subfigure}
	\hspace*{0.08\textwidth} \\
	\begin{subfigure}[b]{0.45\textwidth}
		\centering
		\includegraphics[width=\textwidth]{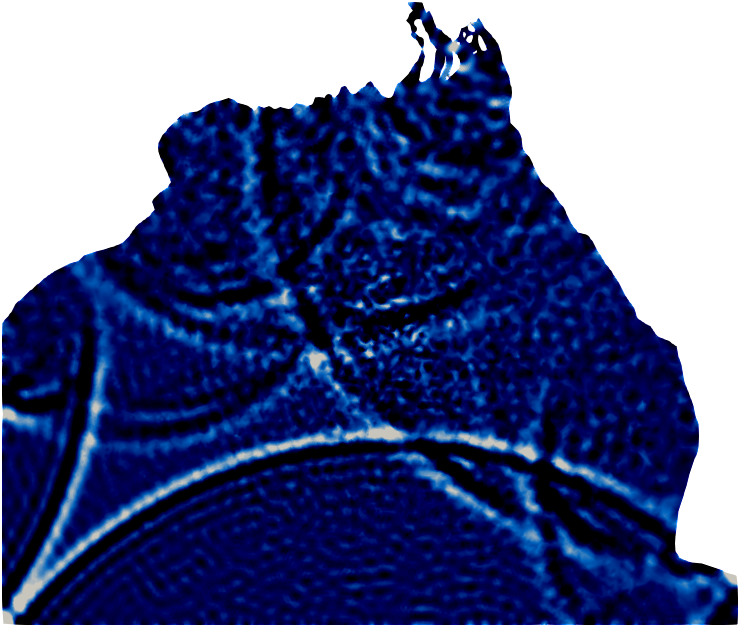}
		\caption{(f1) $t=5{,}000$}\label{fig:i5}
	\end{subfigure}    
	\begin{subfigure}[b]{0.45\textwidth}\hspace{6em}
		\centering
		\includegraphics[width=\textwidth]{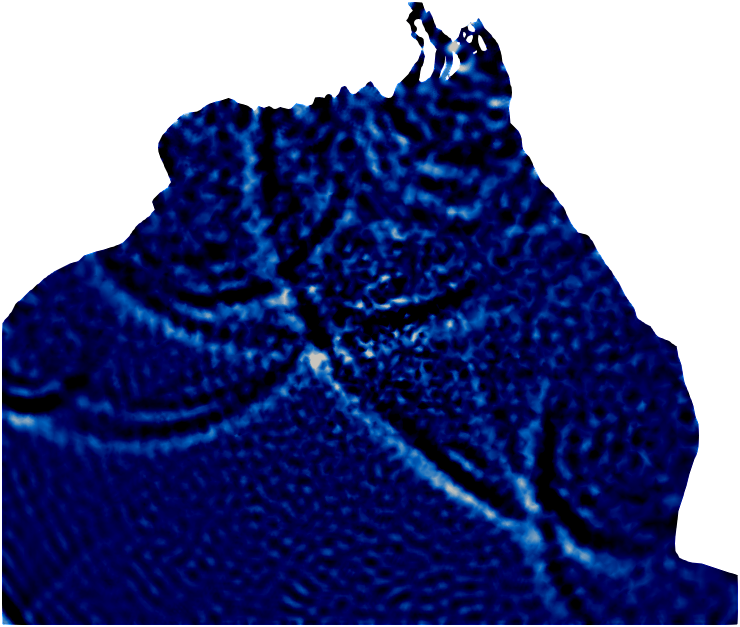}
		\caption{(f2) $t=5{,}000$}\label{fig:i6}
	\end{subfigure}
	\begin{subfigure}[b]{0.08\textwidth}
		\centering
		\includegraphics[width=\textwidth]{cb.png}
		\caption{}
	\end{subfigure}
	\caption{Contour plot of $\eta_h^n$ by LG2 with $\Gamma = \Gamma_\rmD$~(left) and $\Gamma = \bar{\Gamma}_\rmD \cup \bar{\Gamma}_\rmT$~(right) on the Bay of Bengal for~$t=4{,}000, 4{,}500$ and $5{,}000$.}\label{fig:bob_part2}
\end{figure}

We prepare a triangular mesh of the domain as shown in Figure~\ref{mesh}, where the numbers of elements and nodal points are $60{,}619$ and $31{,}120$, respectively.
Then, a numerical simulation is done by LG2 with $\Delta t = 0.2$~[\si{s}].
The results at $t = 0, 2{,}500, 3{,}000, 4{,}000, 4{,}500$ and $5{,}000$~[\si{s}] are presented in Figures~\ref{fig:bob_part1} and~\ref{fig:bob_part2}.
In the figures, for comparison to see the effect of the TBC, we compute Example~\ref{ex4} by replacing $\Gamma_\rmT$ with $\Gamma_\rmD$ and put it on the left.
From Figure~\ref{fig:bob_part1}, we can see that a circular wave is created at around the point~$p$, that it propagates towards the boundary over time, that reflections are found when the wave touches~$\Gamma_\rmD$, and that the results with $\Gamma = \Gamma_\rmD$~(left) and $\Gamma = \Gamma_\rmD \cup \Gamma_\rmT$~(right) are similar.
From Figure~\ref{fig:bob_part2}, we can observe that artificial reflections on the open sea boundaries are significantly removed when the wave touches~$\Gamma_\rmT$, cf. the right figures.
Thus, LG2 works well for a simple (square) domain and this complex domain, the Bay of Bengal region, which is non-convex and includes islands.
\par
\begin{figure}[htbp]
	\centering
	\includegraphics[width=0.7\textwidth]{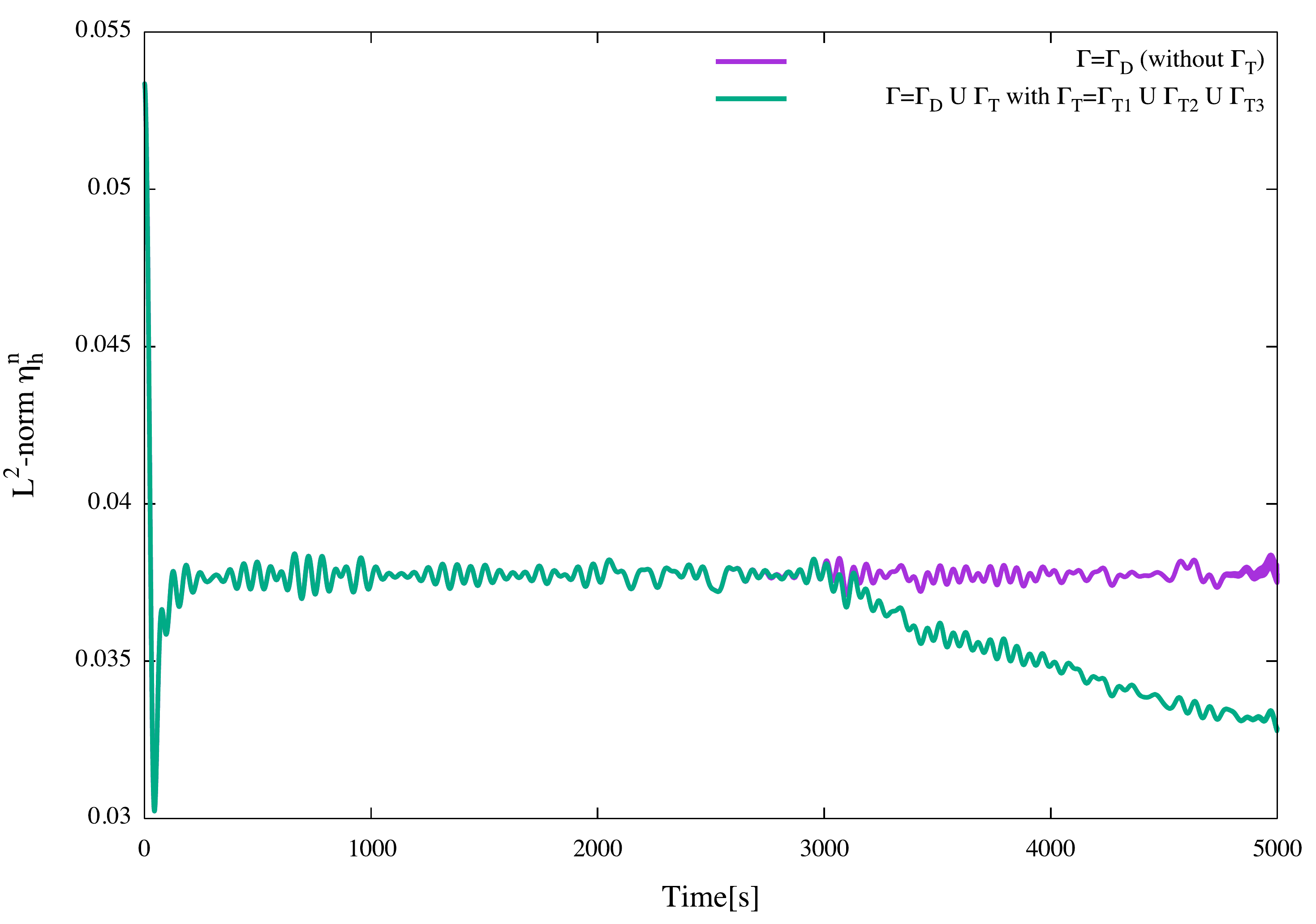}
	\caption{Graphs of~$\|\eta_h^n\|_{L^{2}(\Omega)}$ with respect to time~$(t=t^n)$ for Example~\ref{ex4} with~$\Gamma_\rmT$~$(\Gamma=\Gamma_\rmD\cup\Gamma_\rmT)$ and without~$\Gamma_\rmT$~$(\Gamma=\Gamma_\rmD)$.}
	\label{fig:L2_eta}
\end{figure}
\begin{figure}[htbp]
	\centering
	\includegraphics[width=0.7\textwidth]{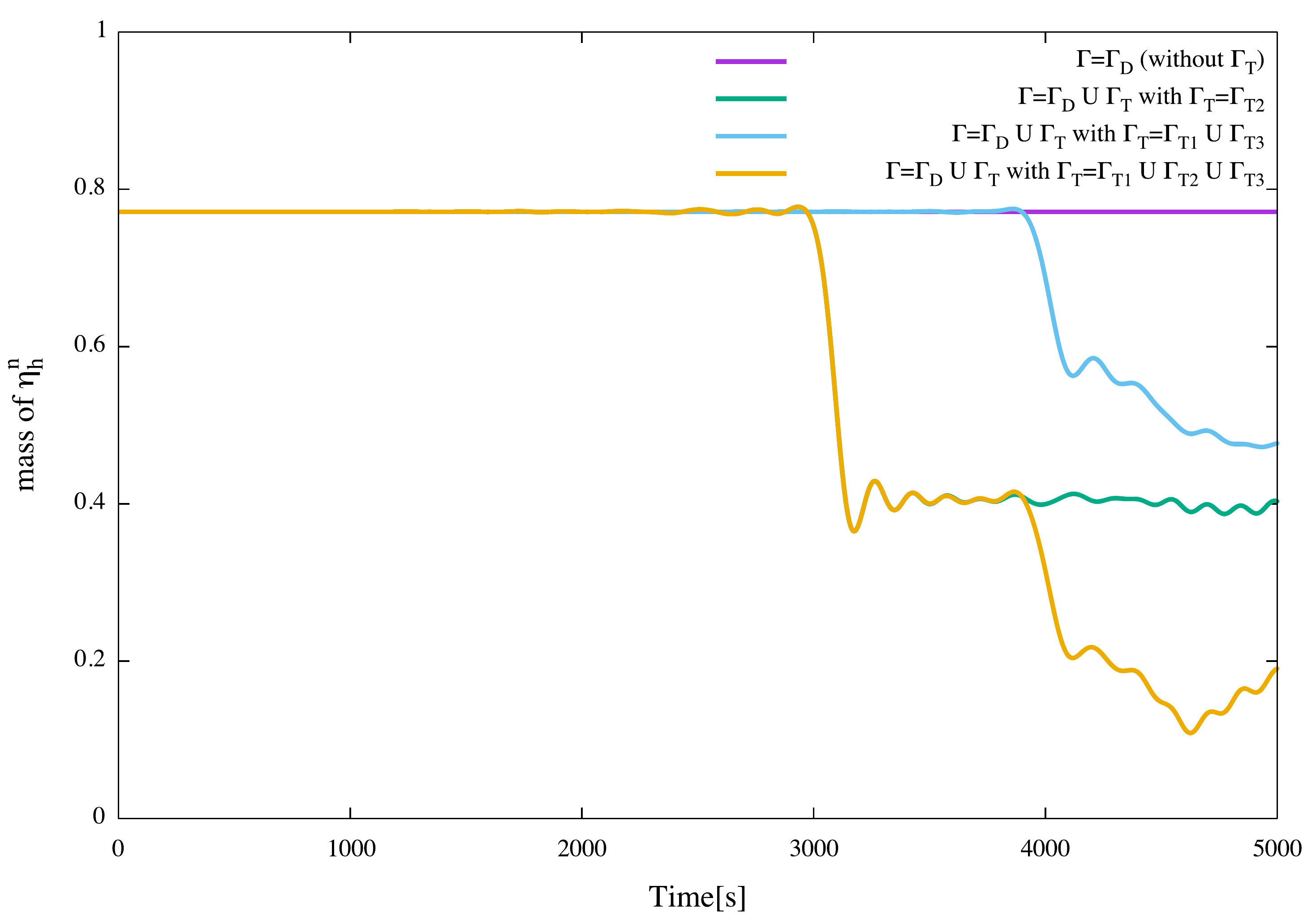}
	\caption{Graphs of the mass of~$\eta_h^n$ with respect to time~$(t=t^n)$ for Example~\ref{ex4} with the following four settings, (i) no transmission boundary, i.e., $\Gamma_\rmT=\emptyset$~(purple), (ii) one transmission boundary, i.e., $\Gamma_\rmT=\Gamma_{\rmT 2}$~(green) (iii) two transmission boundaries~(blue), i.e., $\Gamma_\rmT=\Gamma_{\rmT 1}\cup\Gamma_{\rmT 3}$~(blue), and (iv) three transmission boundaries, i.e., $\Gamma_\rmT=\Gamma_{\rmT 1}\cup\Gamma_{\rmT 2}\cup\Gamma_{\rmT 3}$~(yellow)}
	\label{fig:mass_eta}
\end{figure}

For any (smooth) solution to problem~\eqref{eqn1}, we define the total energy~$\mathcal{E} (t)$ by
\begin{align}
\mathcal{E}(t) \defeq \mathcal{E}_1(t) + \mathcal{E}_2(t)
\defeq \int_{\Omega}\frac{\rho}{2}\phi\vert u\vert^2 dx + \int_{\Omega}\frac{\rho g \vert\eta\vert^2}{2} dx,
\label{eq:energy}
\end{align}
where $\mathcal{E}_1(t)$ is the kinetic energy, and $\mathcal{E}_2(t)$ is the potential energy.
Then, it is worthy to note that the following energy estimate holds, cf.~\cite[Corollary 3.3-(i)]{murshed2021},
\begin{align*}
\fz{d}{dt} \mathcal{E} (t) & = - \frac{\rho}{2} \int_{\Gamma_\rmT} \phi |u|^2 (u \cdot n) \, ds 
- \rho g \int_{\Gamma_\rmT} \phi \, \eta (u \cdot n) \, ds \\
& \quad + 2\mu \int_{\Gamma_\rmT} \phi \bigl( [ D(u) n ] \cdot  u \bigr) \, ds
-2\mu \int_{\Omega} \phi \, \vert D(u) \vert^2 \, dx.
\end{align*}
Here, focusing on $\mathcal{E}_2(t)$~$(=\fz{1}{2}\int_{\Omega}\rho g \vert\eta\vert^2 dx)$ and the mass of~$\eta$, i.e., $\int_{\Omega}\eta\,dx$, we present the values of the $L^2(\Omega)$-norm of~$\eta_h^n$, i.e, $\|\eta_h^n\|_{L^2(\Omega)}$, and the mass of~$\eta_h^n$, i.e., $\int_\Omega\eta_h^n\,dx$, in Figures~\ref{fig:L2_eta} and~\ref{fig:mass_eta}, respectively.
In principle, we can say that the TBC works well numerically if $\|\eta_h^n\|_{L^2(\Omega)}$ and $\int_\Omega\eta_h^n\,dx$ decrease around the time that the wave touches the transmission boundaries.
Figure~\ref{fig:L2_eta} shows graphs of~$\|\eta_h^n\|_{L^2(\Omega)}$ for the two cases, with and without the transmission boundaries, i.e., $\Gamma = \Gamma_\rmD\cup\Gamma_\rmT$ and~$\Gamma = \Gamma_\rmD~(\Gamma_\rmT = \emptyset)$, respectively.
Figure~\ref{fig:mass_eta} shows the graphs of~$\int_\Omega\eta_h^n\,dx$ for the four cases of (transmission) boundaries, (i)~no transmission boundary, i.e., $\Gamma_\rmT=\emptyset$, (ii)~one transmission boundary, i.e., $\Gamma_\rmT=\Gamma_{\rmT 2}$, (iii)~two transmission boundaries, i.e., $\Gamma_\rmT=\Gamma_{\rmT 1}\cup\Gamma_{\rmT 3}$, and (iv)~three transmission boundaries, i.e., $\Gamma_\rmT=\Gamma_{\rmT 1}\cup\Gamma_{\rmT 2}\cup\Gamma_{\rmT 3}$.
From Figures~\ref{fig:L2_eta} and~\ref{fig:mass_eta}, we can see that there are decreasing phenomena of the value of $L^2(\Omega)$-norm as well as the value of the mass when the TBC is imposed.
From Figure~\ref{fig:bob_part1}, we can see that the wave touches the transmission boundary~$\Gamma_{\rmT 2}$ at time around $t = 3{,}000$~[s];  that is why, the mass of~$\eta_h^n$ decreases drastically from around $3{,}000$~[s] to $3{,}200$~[s], cf. Figure~\ref{fig:mass_eta}~(yellow and green lines).
Again, the mass started to decrease between the period from around $4{,}000$~[s] to $4{,}500$~[s], cf. Figure~\ref{fig:mass_eta}, since the wave reached the transmission boundary~$\Gamma_{\rmT 1}$ and~$\Gamma_{\rmT 3}$, cf. Figure~\ref{fig:bob_part2}.
\subsection{Effect of position of a transmission boundary}
\label{subsec:position_of_TBC}
\begin{figure}[!htbp] 
	\centering
	%\hline
	\begin{minipage}[t]{0.45\textwidth}
		\centering
		\begin{overpic}[width=\textwidth]{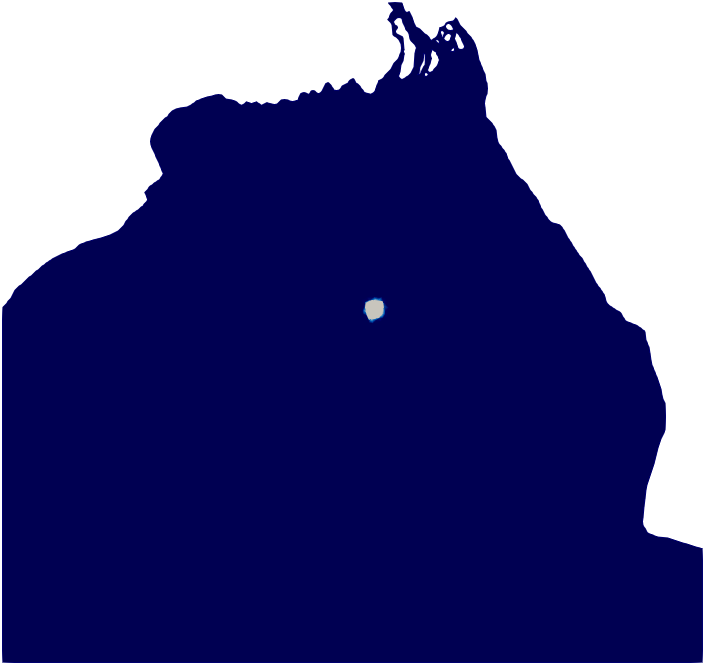}\put(1,10){\tikz \draw[dashed,green] (-0.0,0)--(5.4,0);}\end{overpic}
		\caption*{(a1) $ t=0$}\label{fig:b1}
	\end{minipage}
	\begin{minipage}[t]{0.45\textwidth}
		\centering
		\raisebox{14.7pt}{\includegraphics[width=\textwidth]{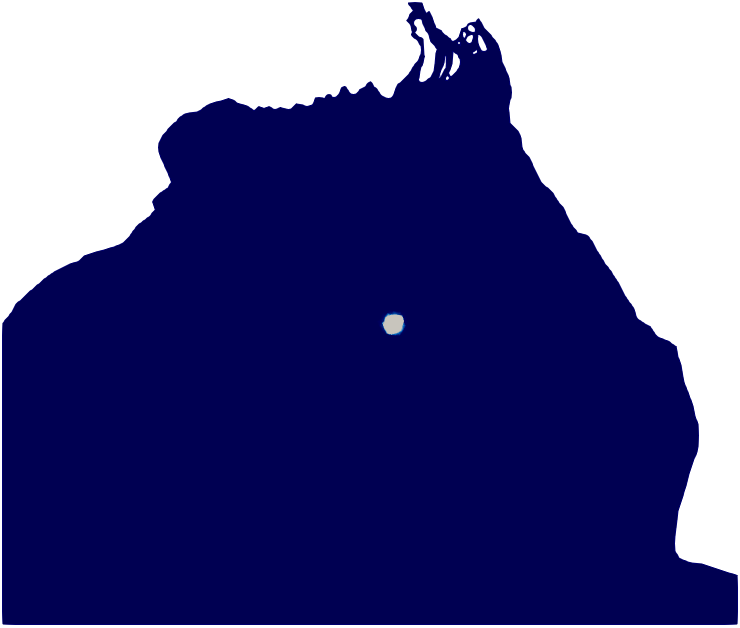}}
		\caption*{(a2) $t=0$}\label{fig:b2}
	\end{minipage}
	\hspace*{0.08\textwidth} \\
	\begin{minipage}[b]{0.45\textwidth}
		\centering
		\begin{overpic}[width=\textwidth]{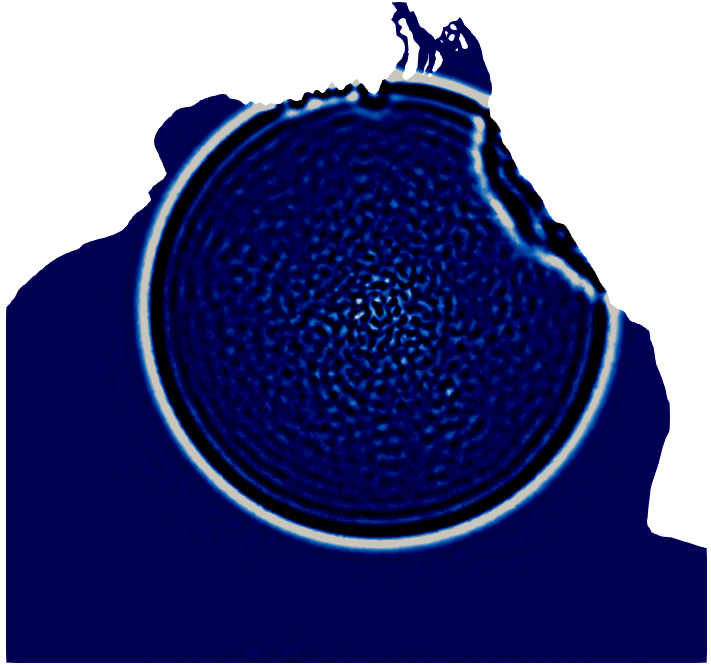}\put(1,10){\tikz \draw[dashed,green] (-0.0,0)--(5.4,0);}\end{overpic}
		\caption*{(b1) $t=2{,}500$}\label{fig:b3}
	\end{minipage}
	\begin{minipage}[b]{0.45\textwidth}
		\centering
		\raisebox{14.7pt}{\includegraphics[width=\textwidth]{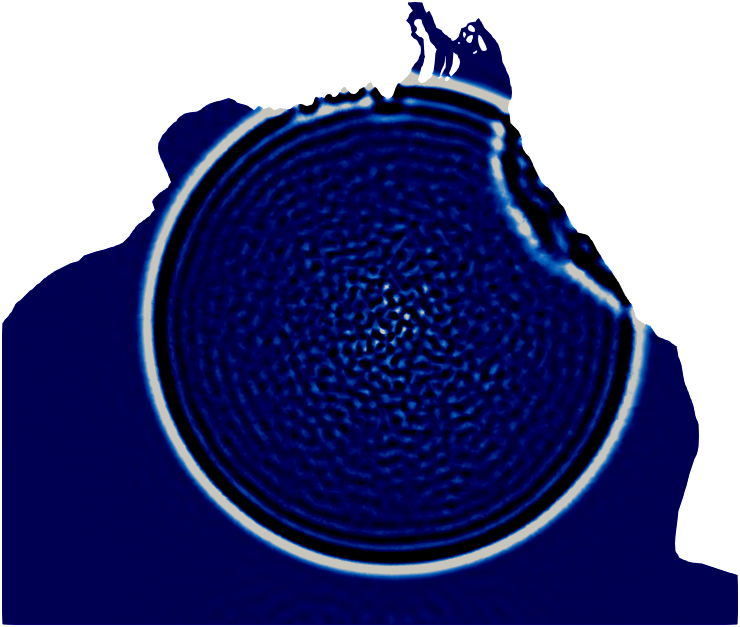}}
		\caption*{(b2) $t=2{,}500$}\label{fig:b4}
	\end{minipage} 
	\hspace*{0.08\textwidth} \\
	\begin{minipage}[b]{0.45\textwidth}
		\centering
		\begin{overpic}[width=\textwidth]{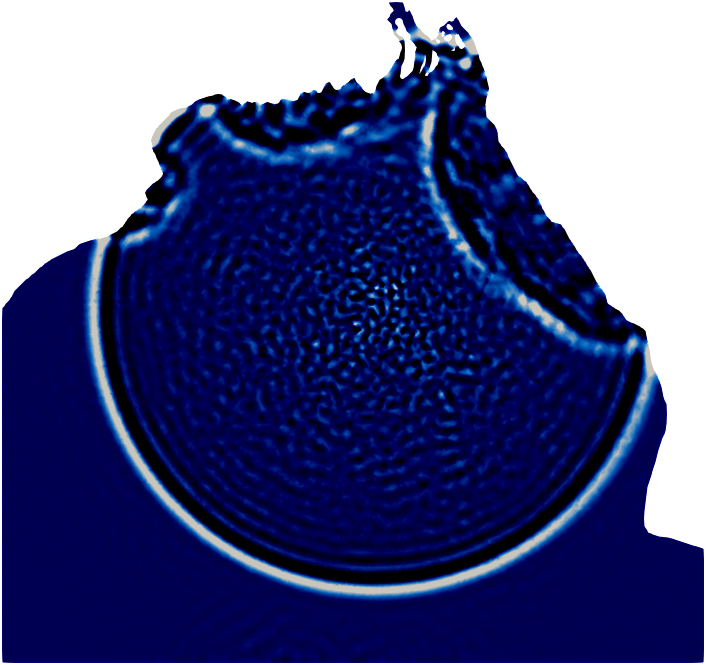}\put(1,10){\tikz \draw[dashed,green] (-0.0,0)--(5.4,0);}\end{overpic}
		\caption*{(c1) $t=3{,}000$}\label{fig:b5}
	\end{minipage}
	\begin{minipage}[b]{0.45\textwidth}
		\centering
		\raisebox{14.7pt}{\includegraphics[width=\textwidth]{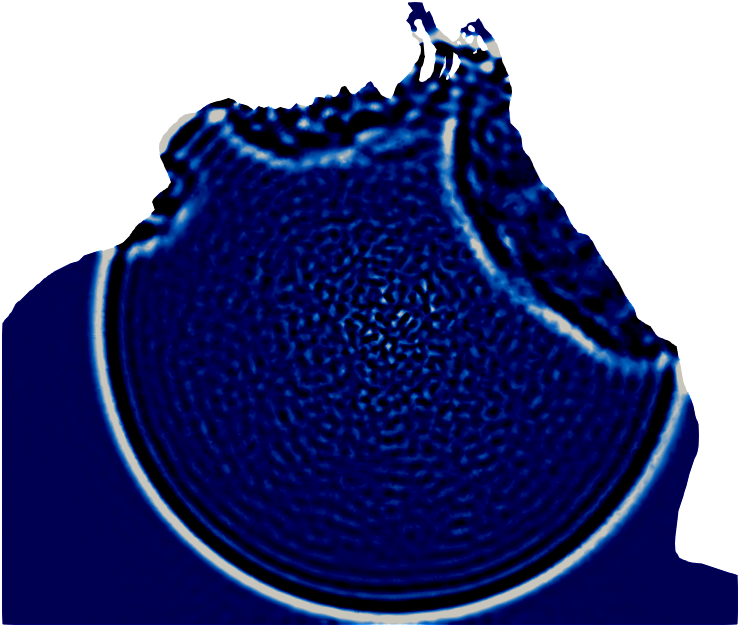}}
		\caption*{(c2) $t=3{,}000$}\label{figb6}
	\end{minipage}
	\begin{minipage}[b]{0.08\textwidth}
		\centering
		\includegraphics[width=\textwidth]{cb.png}
		\caption*{}
	\end{minipage}
	\caption{Contour plot of $\eta_h^n$ by LG2 with $\Gamma = \bar{\Gamma}_\rmD \cup \bar{\Gamma}_\rmT$~for the extended domain(left) and for the original domain (right) on the Bay of Bengal for~$t=0, 2{,}500$ and $3{,}000$.}\label{fig:exbob_part1}
\end{figure}

\begin{figure}[!htbp]
	\centering
	\begin{minipage}[b]{0.45\textwidth}
		\centering
		\begin{overpic}[width=\textwidth]{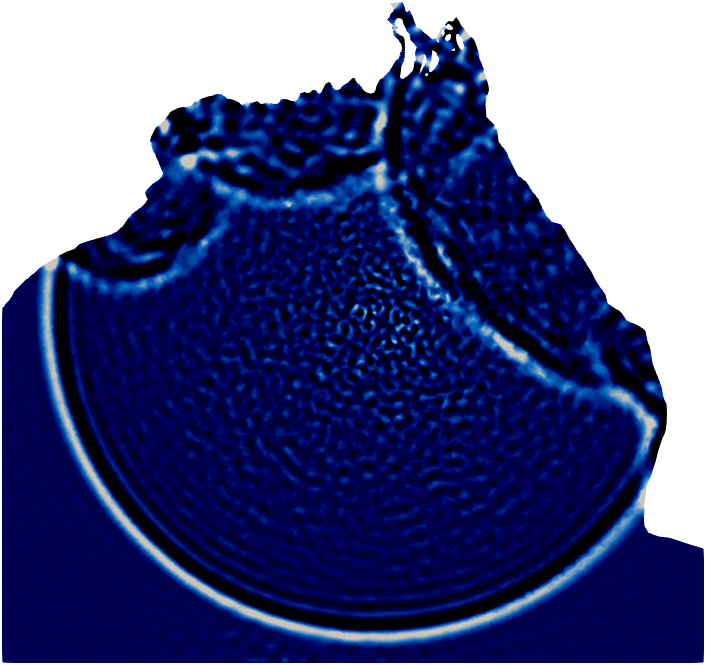}\put(1,10){\tikz \draw[dashed,green] (-0.0,0)--(5.4,0);}\end{overpic}
		\caption*{(d1) $t=3{,}500$}\label{fig:i1}
	\end{minipage}
	\begin{minipage}[b]{0.45\textwidth}
		\centering
		\raisebox{14.7pt}{\includegraphics[width=\textwidth]{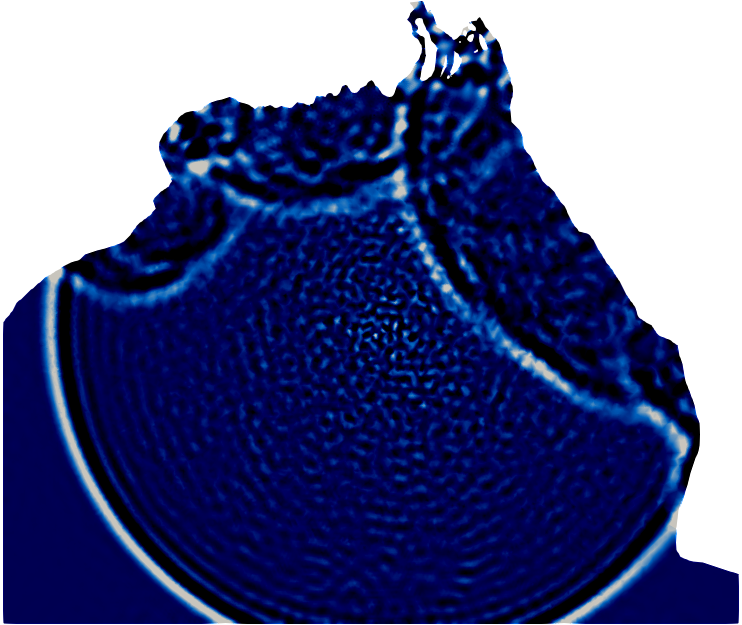}}
		\caption*{(d2) $t=3{,}500$}\label{fig:i2}
	\end{minipage}
	\hspace*{0.08\textwidth} \\
	\begin{minipage}[b]{0.45\textwidth}
		\centering
		\begin{overpic}[width=\textwidth]{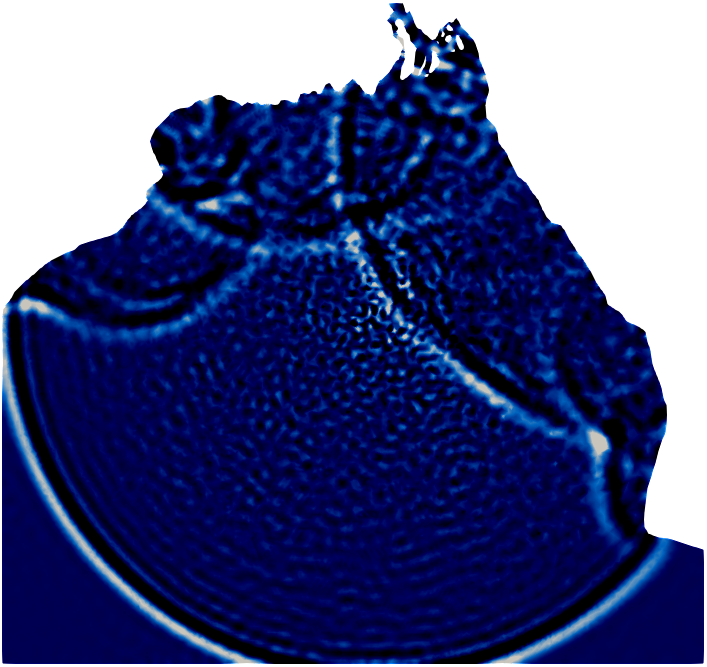}\put(1,10){\tikz \draw[dashed,green] (-0.0,0)--(5.4,0);}\end{overpic}
		\caption*{(e1) $t=4{,}000$}\label{fig:i3}
	\end{minipage}
	\begin{minipage}[b]{0.45\textwidth}
		\centering
		\raisebox{14.7pt}{\includegraphics[width=\textwidth]{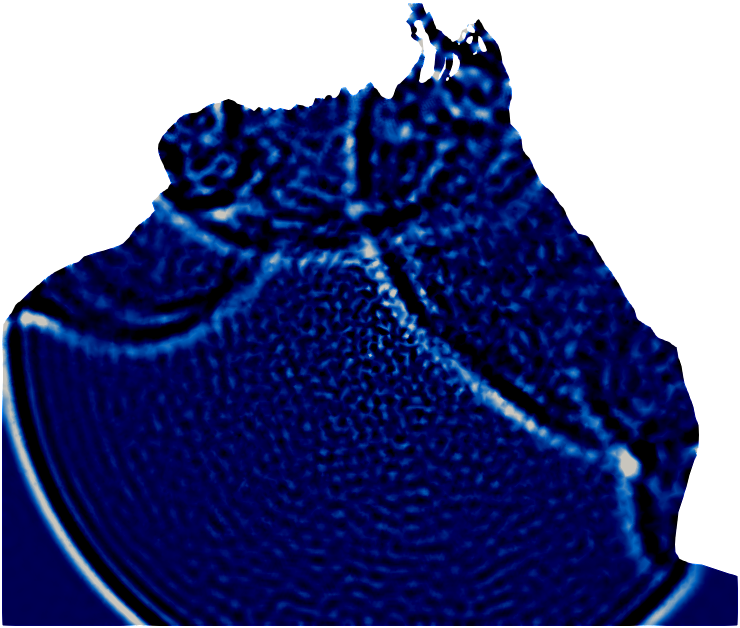}}
		\caption*{(e2) $t=4{,}000$}\label{fig:i4}
	\end{minipage}
	\hspace*{0.08\textwidth} \\
	\begin{minipage}[b]{0.45\textwidth}
		\centering
		\begin{overpic}[width=\textwidth]{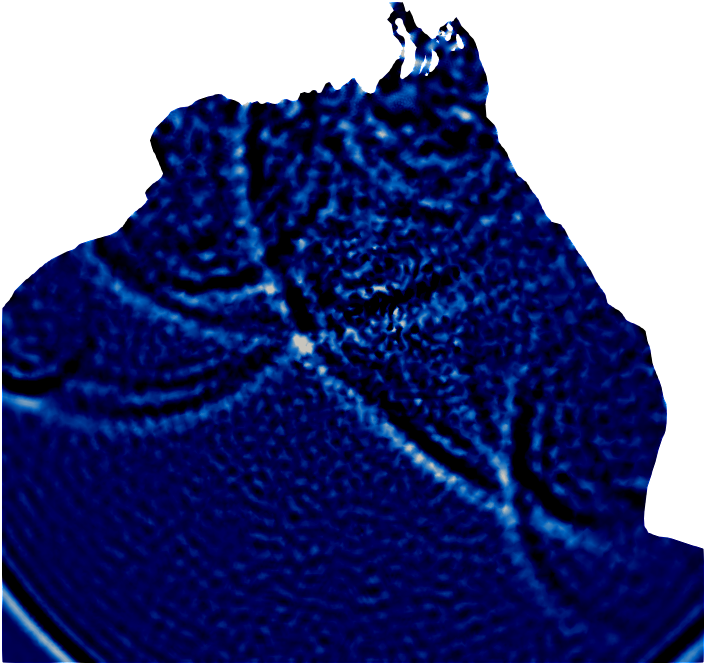}\put(1,10){\tikz \draw[dashed,green] (-0.0,0)--(5.4,0);}\end{overpic}
		\caption*{(f1) $t=5{,}000$}\label{fig:i5}
	\end{minipage}
	\begin{minipage}[b]{0.45\textwidth}\hspace{6em}
		\centering
		\raisebox{14.7pt}{\includegraphics[width=\textwidth]{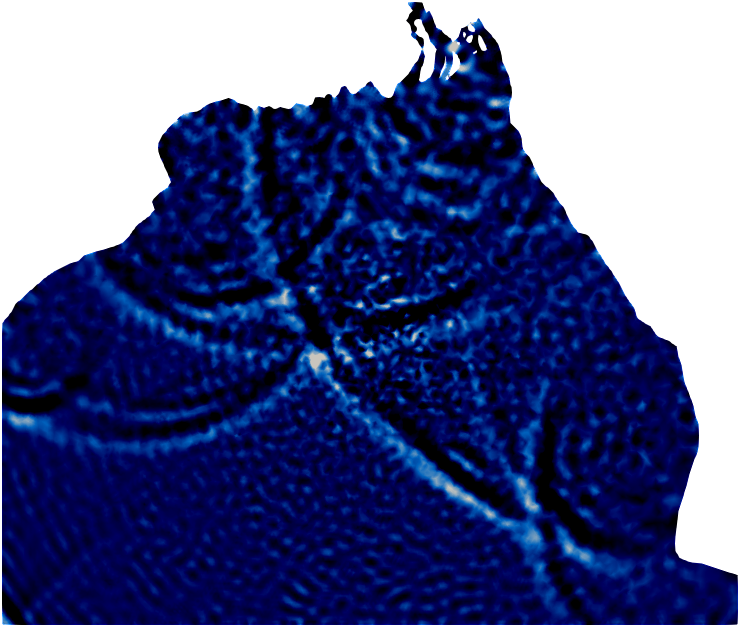}}
		\caption*{(f2) $t=5{,}000$}\label{fig:i6}
	\end{minipage}
	\begin{minipage}[b]{0.08\textwidth}
		\centering
		\includegraphics[width=\textwidth]{cb.png}
		\caption*{}
	\end{minipage}
	\caption{Contour plot of $\eta_h^n$ by LG2 with $\Gamma = \bar{\Gamma}_\rmD \cup \bar{\Gamma}_\rmT$~for the extended domain(left) and for the original domain (right) on the Bay of Bengal for~$t=3{,}500, 4{,}000$ and $5{,}000$.}\label{fig:exbob_part2}
\end{figure}

\begin{figure}[!htbp]
	\centering
	\includegraphics[width=0.8\textwidth]{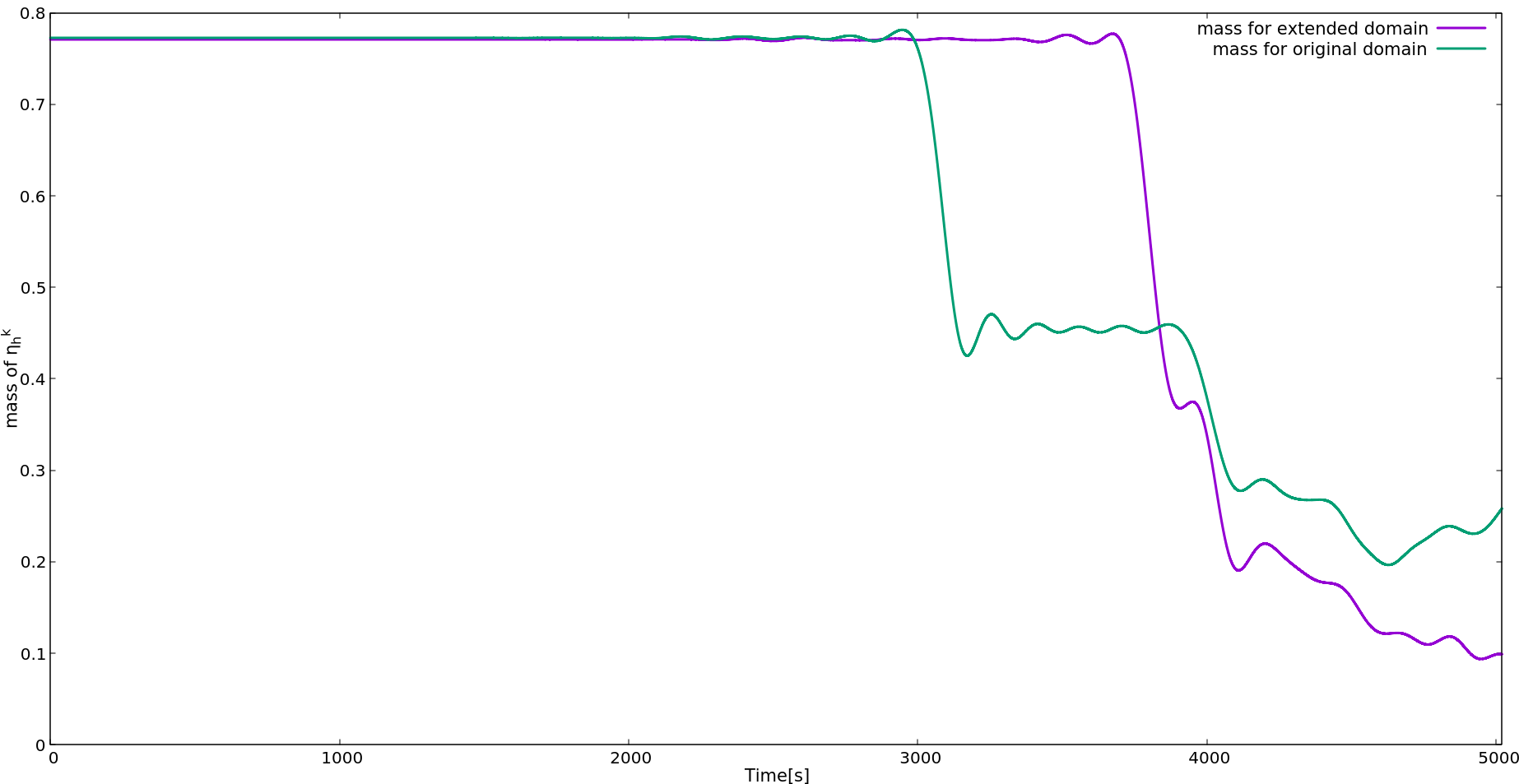}
	\caption{Graphs of~mass of $\eta$ for the extended and original domain with a TBC.}
	\label{fig:extmass}
\end{figure}

We consider Example~\ref{ex4} again to see the effect of the TBC with an extension of the domain ($\Omega$), where the size of the domain in the vertical direction is extended from $889.59~[{\rm km}]$ to $989.59~[{\rm km}]$, i.e., $100~[{\rm km}]$ extension. We employ the same boundary conditions on $\Gamma = \Gamma_{\rmD} \cup \Gamma_{\rmT}$ for both original and extended domains, where $\Gamma_\rmT=\Gamma_{\rmT 1}\cup\Gamma_{\rmT 2}\cup\Gamma_{\rmT 3}$. We compare the numerical results for the extended domain with the ones for the original domain, cf. Figures~\ref{fig:exbob_part1} and~\ref{fig:exbob_part2}, where the left and right figures show the results for the extended and original domains, respectively.
It is observed that there is no significant effect of the vertical position of the bottom transmission boundary~$\Gamma_{\rmT 2}$.
We also computed the mass of~$\eta$ for both domains, cf. Figure~\ref{fig:extmass}. From Figure~\ref{fig:extmass}, we can see that the mass of $\eta_{h}^{k}$ started to decrease at time $t=3{,}000$ for the original domain, cf. Figure \ref{fig:exbob_part1}-(c2), while the mass of~$\eta_{h}^{k}$ started to decrease at time $t = 4{,}000$ for the extended domain, cf. Figure~\ref{fig:exbob_part2}-(e1), because the wave touches the boundary~$\Gamma_{\rmT 2}$ at these times~($t = 3{,}000$ and $t = 4{,}000$) for the original and extended domains, respectively. 
%That means the wave is passing through the boundaries. 
A similar decreasing property of mass of~$\eta_{h}^{k}$ can be observed from Figure~\ref{fig:extmass} when the wave touches the transmission boundaries.
%~$\Gamma_{\rmT 1}$ and~$\Gamma_{\rmT 3}$ for both original and extended domain which 
The results confirm that the TBC works well numerically and that we can choose the vertical position of the bottom transmission boundary~$\Gamma_{\rmT 2}$ without significant effect.
\section{Conclusions}
\label{sec:conclusions}
We have presented a two-step Lagrange--Galerkin scheme for the shallow water equations with a TBC.
For the scheme, the EOCs have been computed (cf. Examples~\ref{ex1} and~\ref{ex2} in Subsection~\ref{subsec:EOC}) and the second-order accuracy in time has been confirmed.
From numerical experiments on a simple square domain (cf. Example~\ref{ex3} in Subsection~\ref{sub:problem}), it has been observed that the effect of the TBC works well.
Our scheme has been applied to a realistic domain, the Bay of Bengal, and numerical experiments have been performed for two different types of boundary conditions, i.e., with and without the TBC (cf. Subsection~\ref{subsec:with_and_without_TBC}).
There have been no significant reflections from~$\Gamma_\rmT$ and the wave has passed through $\Gamma_\rmT$ while reflections have been observed from~$\Gamma_\rmD$, and, in the graphs of~$\|\eta_h^n\|_{L^2(\Omega)}$ and the mass of $\eta_h^n$~(cf.
Figures~\ref{fig:L2_eta} and~\ref{fig:mass_eta}), natural decays of the values of~$\|\eta_h^n\|_{L^2(\Omega)}$ as well as the mass of~$\eta_h^n$ have been observed when the TBC is imposed.
In addition, for the domain extended by $100$~[\si{km}] in the vertical direction, it has been confirmed that there is no significant effect of changing the position of the transmission boundary (cf. Subsection~\ref{subsec:position_of_TBC}).
From these numerical experiments, we conclude that our two-step Lagrange--Galerkin scheme, cf.~\eqref{scheme}, works well numerically not only for a simple domain but also for a complex domain with the TBC if the bottom topography is flat.
In our forthcoming paper, Part~II, the scheme will be applied to rapidly varying bottom surfaces and a real bottom topography of the Bay of Bengal region to investigate the effect of non-homogeneity of the bottom topography.

\subsection*{Acknowledgements}
M.M.R. is supported by the MEXT scholarship. This work is partially supported by JSPS KAKENHI Grant Numbers JP20KK0058, JP21H00999, JP20H00117, JP20H01812, JP18H01135, JP21H04431, and~JP20H01823, and JST CREST Grant Number~JPMJCR2014.
%
%
%
%
%%%%%%%%%%%%%%%%%%%%%%%%%%%%%%%%%%%%%%%%%%%%%
\appendix
\renewcommand{\thesection}{A}
\setcounter{table}{0}
\renewcommand{\thetable}{\thesection.\arabic{table}}
\setcounter{equation}{0}
\makeatletter
    \renewcommand{\theequation}{%
    \thesection.\arabic{equation}}
    \@addtoreset{equation}{section}
  \makeatother
%
%
%%%%%%%%%%%%%%%%%%%%%%%%%%%%%%%%%%%%%%%%%%%%%
\section*{Appendix}
%%%%%%%%%%%%%%%%%%%%%%%%%%%%%%%%%%%%%%%%%%%%%
\subsection{Choice of $c_0$}
\label{subsec:choice_c0}
\label{subsec:choice_c0}

Based on~\cite{murshed2021}, focusing on the potential energy~$\mathcal{E}_2(t)$, cf.~\eqref{eq:energy}, we perform numerical experiments for the choice of $c_{0}$ for two cases with the following settings: \bigskip\\
\textbf{Case~I}~(the square domain).\ 
In problem~\eqref{eqn1}, we set $\Omega=(0,10)^2$, $T=100$, $g=9.8\times 10^{-3}$, $\rho=10^{12}$, $\mu = \zeta = 1$, $(f, F)=(0, 0)$, $c=10^{-3}$, $\eta^{0}=c \exp(-100\vert x-p \vert^2)$, $p=(5, 5)^\top$, $u^{0}=0$ and~$\Gamma = \Gamma_\rmT$~$(\Gamma_\rmD = \emptyset)$.
We employ discretization parameters, $N=200$~$(h=1/N)$, and~$\Delta t=0.25\sqrt{h}$.
\bigskip\\
\textbf{Case~II}~(the Bay of Bengal).\ 
The parameters are the same as Example~\ref{ex4} except the value of~$c_0$. We employ the same mesh and $\Delta t~(= 0.2)$ in Section~\ref{sec:application}.
\bigskip
\par
For $\eta_h = \{\eta_h^n\}_{n=1}^{N_T}$, let $\|\eta_h\|_{\ell^2(L^2)}$ be a norm of~$\eta_h$ defined by
\[
\|\eta_h\|_{\ell^2(L^2)} \defeq \sqrt{ \Delta t \sum_{n=1}^{N_T} \|\eta_h^n\|_{L^2(\Omega)}^2 }
\quad (\approx \|\eta\|_{L^2(0,T; L^2(\Omega))}).
\]
We compute the two cases for $c_0 = 0.5, 0.6, \ldots$, and~$1.2$.
The results are shown in Table~\ref{choice} and imply that, for both cases, we have minimum values of $\|\eta_h\|_{\ell^2(L^2)}$ for $c_{0}=0.9$.
\begin{table}[!htbp]
	\centering
	\caption{Values of $c_{0}$ and $\|\eta_h\|_{\ell^2(L^2)}$.}\label{choice}
	\begin{tabular}{ccc}
		\toprule
		& \multicolumn{2}{c}{$\|\eta_h\|_{\ell^2(L^2)}$} \\ \cline{2-3}
		Value of $c_{0}$ & Case~I~(the square domain)  & Case~II~(the Bay of Bengal) \\
		\hline\hline
		0.5 & $8.16\phantom{0}\times 10^{-2}$ &13.55\phantom{00} \\
		0.6 & $8.08\phantom{0} \times 10^{-2}$ &13.54\phantom{00} \\
		0.7& $8.03\phantom{0} \times 10^{-2}$ & 13.5342 \\
		0.8 & $8.002\times 10^{-2}$&13.5323 \\
		\color{blue}0.9 & \color{blue}$7.997\times 10^{-2}$ & \color{blue}13.5319 \\
		1.0 & $8.006 \times 10^{-2}$ & 13.5328 \\
		1.1 & $8.02\phantom{4} \times 10^{-2}$ & 13.5354\\
		1.2 & $8.05\phantom{0} \times 10^{-2}$ & 13.5375\\
		\hline
	\end{tabular}
\end{table}
%
%
%
%
%
%
%
%\bibliography{mybibfile}

\begin{thebibliography}{43}
\expandafter\ifx\csname natexlab\endcsname\relax\def\natexlab#1{#1}\fi
\providecommand{\url}[1]{\texttt{#1}}
\providecommand{\href}[2]{#2}
\providecommand{\path}[1]{#1}
\providecommand{\DOIprefix}{doi:}
\providecommand{\ArXivprefix}{arXiv:}
\providecommand{\URLprefix}{URL: }
\providecommand{\Pubmedprefix}{pmid:}
\providecommand{\doi}[1]{\href{http://dx.doi.org/#1}{\path{#1}}}
\providecommand{\Pubmed}[1]{\href{pmid:#1}{\path{#1}}}
\providecommand{\bibinfo}[2]{#2}
\ifx\xfnm\undefined \def\xfnm[#1]{\unskip,\space#1}\fi
%Type = Article
\bibitem{AchGue-2000}
\bibinfo{author}{Achdou\xfnm[ Y.]}, \bibinfo{author}{Guermond\xfnm[ J.L.]}.
\newblock \bibinfo{title}{Convergence analysis of a finite element
  projection/{L}agrange--{G}alerkin method for the incompressible
  {N}avier--{S}tokes equations}.
\newblock \bibinfo{journal}{SIAM Journal on Numerical Analysis}
  \bibinfo{year}{2000};\bibinfo{volume}{37}:\bibinfo{pages}{799--826}.
%Type = Article
\bibitem{BenBer-2011}
\bibinfo{author}{Ben\'{i}tez\xfnm[ M.]}, \bibinfo{author}{Berm\'{u}dez\xfnm[
  A.]}.
\newblock \bibinfo{title}{A second order characteristics finite element scheme
  for natural convection problems}.
\newblock \bibinfo{journal}{Journal of Computational and Applied Mathematics}
  \bibinfo{year}{2011};\bibinfo{volume}{235}:\bibinfo{pages}{3270--3284}.
%Type = Article
\bibitem{BenBer-2012_part1}
\bibinfo{author}{Ben\'{i}tez\xfnm[ M.]}, \bibinfo{author}{Berm\'{u}dez\xfnm[
  A.]}.
\newblock \bibinfo{title}{Numerical analysis of a second order pure
  {L}agrange--{G}alerkin method for convection-diffusion problems. {P}art~{I}:
  {T}ime discretization}.
\newblock \bibinfo{journal}{SIAM Journal on Numerical Analysis}
  \bibinfo{year}{2012}{\natexlab{a}};\bibinfo{volume}{50}:\bibinfo{pages}{858--882}.
%Type = Article
\bibitem{BenBer-2012_part2}
\bibinfo{author}{Ben\'{i}tez\xfnm[ M.]}, \bibinfo{author}{Berm\'{u}dez\xfnm[
  A.]}.
\newblock \bibinfo{title}{Numerical analysis of a second order pure
  {L}agrange--{G}alerkin method for convection-diffusion problems. {P}art~{II}:
  {F}ully discretized scheme and numerical results}.
\newblock \bibinfo{journal}{SIAM Journal on Numerical Analysis}
  \bibinfo{year}{2012}{\natexlab{b}};\bibinfo{volume}{50}:\bibinfo{pages}{2824--2844}.
%Type = Article
\bibitem{BerSaa-2012}
\bibinfo{author}{Bermejo\xfnm[ R.]}, \bibinfo{author}{Saavedra\xfnm[ L.]}.
\newblock \bibinfo{title}{Modified {L}agrange--{G}alerkin methods of first and
  second order in time for convection-diffusion problems}.
\newblock \bibinfo{journal}{Numerische {M}athematik}
  \bibinfo{year}{2012};\bibinfo{volume}{120}:\bibinfo{pages}{601--638}.
%Type = Article
\bibitem{BerGalSaa-2012}
\bibinfo{author}{Bermejo\xfnm[ R.]}, \bibinfo{author}{G\'alan~del Sastre\xfnm[
  P.]}, \bibinfo{author}{Saavedra\xfnm[ L.]}.
\newblock \bibinfo{title}{A second order in time modified
  {L}agrange--{G}alerkin finite element method for the incompressible
  {N}avier--{S}tokes equations}.
\newblock \bibinfo{journal}{SIAM {J}ournal on {N}umerical {A}nalysis}
  \bibinfo{year}{2012};\bibinfo{volume}{50}:\bibinfo{pages}{3084--3109}.
%Type = Article
\bibitem{BerNogVaz-2006_part1}
\bibinfo{author}{Berm\'{u}dez\xfnm[ A.]}, \bibinfo{author}{Nogueiras\xfnm[
  M.R.]}, \bibinfo{author}{V\'{a}zquez\xfnm[ C.]}.
\newblock \bibinfo{title}{Numerical analysis of convection‐diffusion‐reaction
  problems with higher order characteristics/finite elements. part i: Time
  discretization}.
\newblock \bibinfo{journal}{SIAM Journal on Numerical Analysis}
  \bibinfo{year}{2006}{\natexlab{a}};\bibinfo{volume}{44}(\bibinfo{number}{5}):\bibinfo{pages}{1829--1853}.
%Type = Article
\bibitem{BerNogVaz-2006_part2}
\bibinfo{author}{Berm\'{u}dez\xfnm[ A.]}, \bibinfo{author}{Nogueiras\xfnm[
  M.R.]}, \bibinfo{author}{V\'{a}zquez\xfnm[ C.]}.
\newblock \bibinfo{title}{Numerical analysis of convection‐diffusion‐reaction
  problems with higher order characteristics/finite elements. part ii: Fully
  discretized scheme and quadrature formulas}.
\newblock \bibinfo{journal}{SIAM Journal on Numerical Analysis}
  \bibinfo{year}{2006}{\natexlab{b}};\bibinfo{volume}{44}(\bibinfo{number}{5}):\bibinfo{pages}{1854--1876}.
%Type = Article
\bibitem{BMMR-1997}
\bibinfo{author}{Boukir\xfnm[ K.]}, \bibinfo{author}{Maday\xfnm[ Y.]},
  \bibinfo{author}{M\'etivet\xfnm[ B.]}, \bibinfo{author}{Razafindrakoto\xfnm[
  E.]}.
\newblock \bibinfo{title}{A high-order characteristics/finite element method
  for the incompressible {N}avier--{S}tokes equations}.
\newblock \bibinfo{journal}{International Journal for Numerical Methods in
  Fluids}
  \bibinfo{year}{1997};\bibinfo{volume}{25}:\bibinfo{pages}{1421--1454}.
%Type = Article
\bibitem{ChrWal-2008}
\bibinfo{author}{Chrysafinos\xfnm[ K.]}, \bibinfo{author}{Walkington\xfnm[
  N.J.]}.
\newblock \bibinfo{title}{Lagrangian and moving mesh methods for the convection
  diffusion equation}.
\newblock \bibinfo{journal}{ESAIM: Mathematical Modelling and Numerical
  Analysis} \bibinfo{year}{2008};\bibinfo{volume}{42}:\bibinfo{pages}{25--55}.
%Type = Article
\bibitem{ColCarBer-2020}
\bibinfo{author}{Colera\xfnm[ M.]}, \bibinfo{author}{Carpio\xfnm[ J.]},
  \bibinfo{author}{Bermejo\xfnm[ R.]}.
\newblock \bibinfo{title}{A nearly-conservative high-order
  {L}agrange--{G}alerkin method for the resolution of scalar
  convection-dominated equations in non-divergence-free velocity fields}.
\newblock \bibinfo{journal}{Computer Methods in Applied Mechanics and
  Engineering}
  \bibinfo{year}{2020};\bibinfo{volume}{372}:\bibinfo{pages}{113366}.
%Type = Article
\bibitem{ColCarBer-2021}
\bibinfo{author}{Colera\xfnm[ M.]}, \bibinfo{author}{Carpio\xfnm[ J.]},
  \bibinfo{author}{Bermejo\xfnm[ R.]}.
\newblock \bibinfo{title}{A nearly-conservative, high-order, forward
  {L}agrange--{G}alerkin method for the resolution of scalar hyperbolic
  conservation laws}.
\newblock \bibinfo{journal}{Computer {M}ethods in {A}pplied {M}echanics and
  {E}ngineering}
  \bibinfo{year}{2021};\bibinfo{volume}{376}:\bibinfo{pages}{113654}.
%Type = Article
\bibitem{das1972}
\bibinfo{author}{Das\xfnm[ P.K.]}.
\newblock \bibinfo{title}{Prediction model for storm surges in the {B}ay of
  {B}engal}.
\newblock \bibinfo{journal}{Nature}
  \bibinfo{year}{1972};\bibinfo{volume}{239}(\bibinfo{number}{5369}):\bibinfo{pages}{211--213}.
%Type = Article
\bibitem{debsarma2009}
\bibinfo{author}{Debsarma\xfnm[ S.K.]}.
\newblock \bibinfo{title}{Simulations of storm surges in the {B}ay of
  {B}engal}.
\newblock \bibinfo{journal}{Marine Geodesy}
  \bibinfo{year}{2009};\bibinfo{volume}{32}(\bibinfo{number}{2}):\bibinfo{pages}{178--198}.
%Type = Article
\bibitem{douglasrussell1982}
\bibinfo{author}{Douglas\xfnm[ J.J.]}, \bibinfo{author}{Russell\xfnm[ T.F.]}.
\newblock \bibinfo{title}{Numerical methods for convection-dominated diffusion
  problems based on combining the method of characteristics with finite element
  or finite difference procedures}.
\newblock \bibinfo{journal}{SIAM Journal on Numerical Analysis}
  \bibinfo{year}{1982};\bibinfo{volume}{19}(\bibinfo{number}{5}):\bibinfo{pages}{871--885}.
%Type = Inproceedings
\bibitem{EwiRus-1981}
\bibinfo{author}{Ewing\xfnm[ R.]}, \bibinfo{author}{Russell\xfnm[ T.]}.
\newblock \bibinfo{title}{Multistep {G}alerkin methods along characteristics
  for convection-diffusion problems}.
\newblock In: \bibinfo{editor}{Vichnevetsky\xfnm[ R.]},
  \bibinfo{editor}{Stepleman\xfnm[ R.]}, editors.
  \bibinfo{booktitle}{{A}dvances in {C}omputer {M}ethods for {P}artial
  {D}ifferential {E}quations {IV}}. \bibinfo{publisher}{IMACS};
  \bibinfo{year}{1981}. p. \bibinfo{pages}{28--36}.
%Type = Inproceedings
\bibitem{EwiRusWhe-1983}
\bibinfo{author}{Ewing\xfnm[ R.]}, \bibinfo{author}{Russell\xfnm[ T.]},
  \bibinfo{author}{Wheeler\xfnm[ M.]}.
\newblock \bibinfo{title}{Simulation of miscible displacement using mixed
  methods and a modified method of characteristics}.
\newblock In: \bibinfo{booktitle}{Proceedings of the Seventh Reservoir
  Simulation Symposium}. \bibinfo{publisher}{Society of Petroleum Engineers of
  AIME}; \bibinfo{year}{1983}. p. \bibinfo{pages}{71--81}.
%Type = Article
\bibitem{FutKolNotSuz-2022}
\bibinfo{author}{Futai\xfnm[ K.]}, \bibinfo{author}{Kolbe\xfnm[ N.]},
  \bibinfo{author}{Notsu\xfnm[ H.]}, \bibinfo{author}{Suzuki\xfnm[ T.]}.
\newblock \bibinfo{title}{A mass-preserving two-step {L}agrange--{G}alerkin
  scheme for convection-diffusion problems}.
\newblock \bibinfo{journal}{Journal of Scientific Computing}
  \bibinfo{year}{2022};\bibinfo{volume}{92}(\bibinfo{number}{2}):\bibinfo{pages}{37}.
%Type = Article
\bibitem{MR3043640}
\bibinfo{author}{Hecht\xfnm[ F.]}.
\newblock \bibinfo{title}{New development in {F}ree{F}em++}.
\newblock \bibinfo{journal}{Journal of Numerical Mathematics}
  \bibinfo{year}{2012};\bibinfo{volume}{20}(\bibinfo{number}{3-4}):\bibinfo{pages}{251--265}.
%Type = Article
\bibitem{johns1981}
\bibinfo{author}{Johns\xfnm[ B.]}.
\newblock \bibinfo{title}{Numerical simulation of storm surges in the {B}ay of
  {B}engal}.
\newblock \bibinfo{journal}{Monsoon Dynamics}
  \bibinfo{year}{1981};:\bibinfo{pages}{689--706}.
%Type = Article
\bibitem{kanayamadan2006}
\bibinfo{author}{Kanayama\xfnm[ H.]}, \bibinfo{author}{Dan\xfnm[ H.]}.
\newblock \bibinfo{title}{A finite element scheme for two-layer viscous
  shallow-water equations}.
\newblock \bibinfo{journal}{Japan Journal of Industrial and Applied
  Mathematics}
  \bibinfo{year}{2006};\bibinfo{volume}{23}(\bibinfo{number}{2}):\bibinfo{pages}{163--191}.
%Type = Article
\bibitem{kanayamadan2016}
\bibinfo{author}{Kanayama\xfnm[ H.]}, \bibinfo{author}{Dan\xfnm[ H.]}.
\newblock \bibinfo{title}{Tsunami propagation from the open sea to the coast}.
\newblock \bibinfo{journal}{Tsunami} \bibinfo{year}{2016};.
%Type = Article
\bibitem{LMNT-2017_Peterlin_Oseen_Part_I}
\bibinfo{author}{Luk\'{a}\v{c}ov\'{a}-Medvid'ov\'{a}\xfnm[ M.]},
  \bibinfo{author}{Mizerov\'{a}\xfnm[ H.]}, \bibinfo{author}{Notsu\xfnm[ H.]},
  \bibinfo{author}{Tabata\xfnm[ M.]}.
\newblock \bibinfo{title}{Numerical analysis of the {Oseen}-type {Peterlin}
  viscoelastic model by the stabilized {L}agrange--{G}alerkin method,
  {P}art~{I}: A linear scheme}.
\newblock \bibinfo{journal}{ESAIM:~M2AN}
  \bibinfo{year}{2017}{\natexlab{a}};\bibinfo{volume}{51}:\bibinfo{pages}{1637--1661}.
%Type = Article
\bibitem{LMNT-2017_Peterlin_Oseen_Part_II}
\bibinfo{author}{Luk\'{a}\v{c}ov\'{a}-Medvid'ov\'{a}\xfnm[ M.]},
  \bibinfo{author}{Mizerov\'{a}\xfnm[ H.]}, \bibinfo{author}{Notsu\xfnm[ H.]},
  \bibinfo{author}{Tabata\xfnm[ M.]}.
\newblock \bibinfo{title}{Numerical analysis of the {Oseen}-type {Peterlin}
  viscoelastic model by the stabilized {L}agrange--{G}alerkin method,
  {P}art~{II}: A nonlinear scheme}.
\newblock \bibinfo{journal}{ESAIM:~M2AN}
  \bibinfo{year}{2017}{\natexlab{b}};\bibinfo{volume}{51}:\bibinfo{pages}{1663--1689}.
%Type = Article
\bibitem{LNS-2015}
\bibinfo{author}{Luk\'{a}\v{c}ov\'{a}-Medvi\v{d}ov\'{a}\xfnm[ M.]},
  \bibinfo{author}{Notsu\xfnm[ H.]}, \bibinfo{author}{She\xfnm[ B.]}.
\newblock \bibinfo{title}{Energy dissipative characteristic schemes for the
  diffusive {O}ldroyd-{B} viscoelastic fluid}.
\newblock \bibinfo{journal}{International Journal for Numerical Methods in
  Fluids} \bibinfo{year}{2015};.
%Type = Phdthesis
\bibitem{murshed_thesis2019}
\bibinfo{author}{Murshed\xfnm[ M.M.]}.
\newblock \bibinfo{title}{Theoretical and Numerical Studies of the Shallow
  Water Equations with a Transmission Boundary Condition}.
\newblock Ph.D. thesis; Kanazawa University, Japan; \bibinfo{year}{2019}.
%Type = Article
\bibitem{murshed2021}
\bibinfo{author}{Murshed\xfnm[ M.M.]}, \bibinfo{author}{Futai\xfnm[ K.]},
  \bibinfo{author}{Kimura\xfnm[ M.]}, \bibinfo{author}{Notsu\xfnm[ H.]}.
\newblock \bibinfo{title}{Theoretical and numerical studies for energy
  estimates of the shallow water equations with a transmission boundary
  condition}.
\newblock \bibinfo{journal}{Discrete and Continuous Dynamical Systems - S}
  \bibinfo{year}{2021};\bibinfo{volume}{14}(\bibinfo{number}{3}):\bibinfo{pages}{1063--1078}.
%Type = Article
\bibitem{N-2008-JSCES}
\bibinfo{author}{Notsu\xfnm[ H.]}.
\newblock \bibinfo{title}{Numerical computations of cavity flow problems by a
  pressure stabilized characteristic-curve finite element scheme}.
\newblock \bibinfo{journal}{Transactions of Japan Society for Computational
  Engineering and Science}
  \bibinfo{year}{2008};\bibinfo{volume}{2008}:\bibinfo{pages}{20080032}.
%Type = Article
\bibitem{notsuruitabata2013}
\bibinfo{author}{Notsu\xfnm[ H.]}, \bibinfo{author}{Rui\xfnm[ H.]},
  \bibinfo{author}{Tabata\xfnm[ M.]}.
\newblock \bibinfo{title}{Development and {L}2-analysis of a single-step
  characteristics finite difference scheme of second order in time for
  convection-diffusion problems}.
\newblock \bibinfo{journal}{Journal of Algorithms \& Computational Technology}
  \bibinfo{year}{2013};\bibinfo{volume}{7}(\bibinfo{number}{3}):\bibinfo{pages}{343--380}.
%Type = Article
\bibitem{notsutabata2015}
\bibinfo{author}{Notsu\xfnm[ H.]}, \bibinfo{author}{Tabata\xfnm[ M.]}.
\newblock \bibinfo{title}{Error estimates of a pressure-stabilized
  characteristics finite element scheme for the oseen equations}.
\newblock \bibinfo{journal}{Journal of Scientific Computing}
  \bibinfo{year}{2015};\bibinfo{volume}{65}(\bibinfo{number}{3}):\bibinfo{pages}{940--955}.
%Type = Article
\bibitem{notsu2016error}
\bibinfo{author}{Notsu\xfnm[ H.]}, \bibinfo{author}{Tabata\xfnm[ M.]}.
\newblock \bibinfo{title}{Error estimates of a stabilized
  {L}agrange--{G}alerkin scheme for the {N}avier--{S}tokes equations}.
\newblock \bibinfo{journal}{ESAIM: Mathematical Modelling and Numerical
  Analysis}
  \bibinfo{year}{2016}{\natexlab{a}};\bibinfo{volume}{50}(\bibinfo{number}{2}):\bibinfo{pages}{361--380}.
%Type = Article
\bibitem{notsutabata2norder2016}
\bibinfo{author}{Notsu\xfnm[ H.]}, \bibinfo{author}{Tabata\xfnm[ M.]}.
\newblock \bibinfo{title}{Error estimates of a stabilized
  {L}agrange--{G}alerkin scheme of second-order in time for the
  {N}avier--{S}tokes equations}.
\newblock \bibinfo{journal}{Mathematical Fluid Dynamics, Present and Future
  Springer Proceedings in Mathematics \& Statistics}
  \bibinfo{year}{2016}{\natexlab{b}};:\bibinfo{pages}{497--530}.
%Type = Article
\bibitem{paul2012tide}
\bibinfo{author}{Paul\xfnm[ G.C.]}, \bibinfo{author}{Ismail\xfnm[ A.I.M.]}.
\newblock \bibinfo{title}{Tide--surge interaction model including air bubble
  effects for the coast of {B}angladesh}.
\newblock \bibinfo{journal}{Journal of the Franklin Institute}
  \bibinfo{year}{2012};\bibinfo{volume}{349}(\bibinfo{number}{8}):\bibinfo{pages}{2530--2546}.
%Type = Article
\bibitem{paul2013contribution}
\bibinfo{author}{Paul\xfnm[ G.C.]}, \bibinfo{author}{Ismail\xfnm[ A.I.M.]}.
\newblock \bibinfo{title}{Contribution of offshore islands in the prediction of
  water levels due to tide--surge interaction for the coastal region of
  {B}angladesh}.
\newblock \bibinfo{journal}{Natural Hazards}
  \bibinfo{year}{2013};\bibinfo{volume}{65}(\bibinfo{number}{1}):\bibinfo{pages}{13--25}.
%Type = Article
\bibitem{paul2018storm}
\bibinfo{author}{Paul\xfnm[ G.C.]}, \bibinfo{author}{Senthilkumar\xfnm[ S.]},
  \bibinfo{author}{Pria\xfnm[ R.]}.
\newblock \bibinfo{title}{Storm surge simulation along the {M}eghna estuarine
  area: an alternative approach}.
\newblock \bibinfo{journal}{Acta Oceanologica Sinica}
  \bibinfo{year}{2018};\bibinfo{volume}{37}(\bibinfo{number}{1}):\bibinfo{pages}{40--49}.
%Type = Article
\bibitem{pironneau1982}
\bibinfo{author}{Pironneau\xfnm[ O.]}.
\newblock \bibinfo{title}{On the transport-diffusion algorithm and its
  applications to the {N}avier-{S}tokes equations}.
\newblock \bibinfo{journal}{Numerische Mathematik}
  \bibinfo{year}{1982};\bibinfo{volume}{38}(\bibinfo{number}{3}):\bibinfo{pages}{309--332}.
%Type = Book
\bibitem{Pir-1989}
\bibinfo{author}{Pironneau\xfnm[ O.]}.
\newblock \bibinfo{title}{{F}inite {E}lement {M}ethods for {F}luids}.
\newblock \bibinfo{address}{Chichester}: \bibinfo{publisher}{John Wiley \&
  Sons}, \bibinfo{year}{1989}.
%Type = Article
\bibitem{PirTab-2010}
\bibinfo{author}{Pironneau\xfnm[ O.]}, \bibinfo{author}{Tabata\xfnm[ M.]}.
\newblock \bibinfo{title}{Stability and convergence of a
  {G}alerkin-characteristics finite element scheme of lumped mass type}.
\newblock \bibinfo{journal}{International Journal for Numerical Methods in
  Fluids}
  \bibinfo{year}{2010};\bibinfo{volume}{64}:\bibinfo{pages}{1240--1253}.
%Type = Article
\bibitem{roy1999polar}
\bibinfo{author}{Roy\xfnm[ G.]}, \bibinfo{author}{Kabir\xfnm[ A.H.]},
  \bibinfo{author}{Mandal\xfnm[ M.]}, \bibinfo{author}{Haque\xfnm[ M.]}.
\newblock \bibinfo{title}{Polar coordinates shallow water storm surge model for
  the coast of {B}angladesh}.
\newblock \bibinfo{journal}{Dynamics of Atmospheres and Oceans}
  \bibinfo{year}{1999};\bibinfo{volume}{29}(\bibinfo{number}{2-4}):\bibinfo{pages}{397--413}.
%Type = Article
\bibitem{ruitabata2002}
\bibinfo{author}{Rui\xfnm[ H.]}, \bibinfo{author}{Tabata\xfnm[ M.]}.
\newblock \bibinfo{title}{A second order characteristic finite element scheme
  for convection-diffusion problems}.
\newblock \bibinfo{journal}{Numerische Mathematik}
  \bibinfo{year}{2002};\bibinfo{volume}{92}(\bibinfo{number}{1}):\bibinfo{pages}{161--177}.
%Type = Article
\bibitem{RuiTab-2010}
\bibinfo{author}{Rui\xfnm[ H.]}, \bibinfo{author}{Tabata\xfnm[ M.]}.
\newblock \bibinfo{title}{A mass-conservative characteristic finite element
  scheme for convection-diffusion problems}.
\newblock \bibinfo{journal}{Journal of Scientific Computing}
  \bibinfo{year}{2010};\bibinfo{volume}{43}:\bibinfo{pages}{416--432}.
%Type = Article
\bibitem{suli_1988}
\bibinfo{author}{S{\" u}li\xfnm[ E.]}.
\newblock \bibinfo{title}{Convergence and nonlinear stability of the
  {L}agrange-{G}alerkin method for the {N}avier-{S}tokes equations}.
\newblock \bibinfo{journal}{Numerische Mathematik}
  \bibinfo{year}{1988};\bibinfo{volume}{53}(\bibinfo{number}{4}):\bibinfo{pages}{459--483}.
%Type = Article
\bibitem{TabUch-2016-CD}
\bibinfo{author}{Tabata\xfnm[ M.]}, \bibinfo{author}{Uchiumi\xfnm[ S.]}.
\newblock \bibinfo{title}{A genuinely stable {L}agrange--{G}alerkin scheme for
  convection-diffusion problems}.
\newblock \bibinfo{journal}{Japan Journal of Industrial and Applied
  Mathematics}
  \bibinfo{year}{2016};\bibinfo{volume}{33}:\bibinfo{pages}{121--143}.
%Type = Article
\bibitem{TabUch-2018-NS}
\bibinfo{author}{Tabata\xfnm[ M.]}, \bibinfo{author}{Uchiumi\xfnm[ S.]}.
\newblock \bibinfo{title}{An exactly computable {L}agrange--{G}alerkin scheme
  for the {N}avier--{S}tokes equations and its error estimates}.
\newblock \bibinfo{journal}{Mathematics of Computation}
  \bibinfo{year}{2018};\bibinfo{volume}{87}:\bibinfo{pages}{39--67}.
%Type = Article
\bibitem{Uch-2019}
\bibinfo{author}{Uchiumi\xfnm[ S.]}.
\newblock \bibinfo{title}{A viscosity-independent error estimate of a
  pressure-stabilized {L}agrange--{G}alerkin scheme for the {O}seen problem}.
\newblock \bibinfo{journal}{Journal of Scientific Computing}
  \bibinfo{year}{2019};\bibinfo{volume}{80}:\bibinfo{pages}{834--858}.
%
%
%
\end{thebibliography}
%\biboptions{sort&compress}
%

%
%
%
%
\end{document}